\documentclass{amsart}

\title[Morphisms of Cartan connections]{Morphisms of Cartan connections}%
\author{Benjamin McKay}%
\address{University College Cork}%
\email{b.mckay@ucc.ie}%

\usepackage{graphicx}
\usepackage{amssymb}
\usepackage[all]{xy}
\usepackage{varioref}
\usepackage{booktabs}
\usepackage{hyperref}
\usepackage{euscript}
\usepackage[latin1]{inputenc}
\usepackage{array}
\usepackage{amsfonts}
\usepackage{amsmath}
\usepackage{amssymb}
\usepackage{amsthm}
\usepackage{fontenc}
\usepackage{graphicx}
\usepackage{lscape}
\usepackage{ifthen}

\theoremstyle{definition}

\theoremstyle{remark}

\numberwithin{equation}{section}
\newtheorem{theorem}{Theorem}
\newtheorem{corollary}{Corollary}
\newtheorem{lemma}{Lemma}
\newtheorem{proposition}{Proposition}
\theoremstyle{remark}
\newtheorem{conjecture}{Conjecture}
\newtheorem{definition}{Definition}
\newtheorem{example}{Example}
\newtheorem{remark}{Remark}

\newcounter{remarkCounter}
\newlength{\setBracketHeight}

\newcommand{\LieDer}{\ensuremath{\EuScript L}}
\newcommand{\hook}{\ensuremath{\mathbin{ \hbox{\vrule height1.4pt
        width4pt depth-1pt \vrule height4pt width0.4pt depth-1pt}}}}

\newcommand{\R}[1]{\ensuremath{\mathbb{R}^{#1}}}
\newcommand{\C}[1]{\ensuremath{\mathbb{C}^{#1}}}

\newcommand{\CP}[1]{\ensuremath{\mathbb{CP}^{#1}}}
\newcommand{\RP}[1]{\ensuremath{\mathbb{RP}^{#1}}}

\newcommand{\Gr}[2]{\ensuremath{\operatorname{Gr}\left({#1},{#2}\right)}}

\newcommand{\SO}[1]{\operatorname{SO}\left({#1}\right)}
\newcommand{\so}[1]{\mathfrak{so}\left({#1}\right)}
\newcommand{\Or}[1]{\operatorname{O}\left({#1}\right)}

\newcommand{\GL}[1]{\operatorname{GL}\left({#1}\right)}

\newcommand{\gl}[1]{\mathfrak{gl}\left({#1}\right)}
\newcommand{\SL}[1]{\operatorname{SL}\left({#1}\right)}
\newcommand{\PSL}[1]{\mathbb{P}\SL{#1}}

\newcommand{\slLie}[1]{\mathfrak{sl}\left({#1}\right)}

\newcommand{\SU}[1]{\operatorname{SU}\left({#1}\right)}

\newcommand{\Sym}[2]{\ensuremath{\operatorname{Sym}^{#1}\left({#2}\right)}}
\newcommand{\Lm}[2]{\ensuremath{\Lambda^{#1} \left ( {#2} \right )}}
\newcommand{\nForms}[2]{\ensuremath{\Omega^{#1} \left ( {#2} \right)}}

\DeclareMathOperator{\Ad}{Ad}

 \DeclareMathOperator{\ad}{ad}

\newcommand{\Proj}[1]{\mathbb{P}^{#1}}

\newcommand{\OO}[1]{
  \ensuremath{
    \mathcal{O}
    \ifthenelse{\equal{#1}{0}}
      {}
      {\left({#1}\right)}
  }
}
\newcommand{\covOO}[1]{
  \ensuremath{
    \tilde{\mathcal{O}}
    \ifthenelse{\equal{#1}{0}}
      {}
      {\left({#1}\right)}
  }
}

\newcommand{\OOp}[2]{
  \ensuremath{
    \mathcal{O}
    \ifthenelse{\equal{#1}{0}}
      {}
      {\left({#1}\right)}
    \ifthenelse{\equal{#2}{1}}
      {}
      {^{\oplus{#2}}}
  }
}

\newcommand{\Ob}{\ensuremath{\operatorname{Ob}}}

\begin{document}

\thanks{Thanks to Charles Frances for pointing out a serious mistake in a prior version of this paper.}%
\subjclass{MA53A40}%
%\keywords{}%

\date{\today}%
%\dedicatory{}%
%\commby{}%
% ----------------------------------------------------------------
\begin{abstract}
We define what we call \emph{morphisms} of Cartan connections. We generalize
the main theorems on Cartan connections to theorems on
morphisms. Many of the known constructions involving
Cartan connections turn out to be examples of morphisms.
We prove some basic results concerning
completeness of Cartan connections. We provide a new method
to prove completeness of Cartan connections using
families of morphisms.
\end{abstract}
\maketitle  
\tableofcontents
% ----------------------------------------------------------------
\section{Introduction}

\begin{definition}\label{def:CartanConnection}
Let $H \subset G$ be a closed subgroup of a Lie group, with Lie
algebras $\mathfrak{h} \subset \mathfrak{g}$. A $G/H$-geometry
(also known as a \emph{Cartan geometry} modelled on $G/H$) on a manifold
$M$ is a principal right $H$-bundle $E \to M$ and a 1-form $\omega
\in \nForms{1}{E} \otimes \mathfrak{g}$ called the \emph{Cartan
connection}, which satisfies the following conditions:
\begin{enumerate}
\item
Denote the right action of $h \in H$ on $e \in E$ by $r_h e$. The
Cartan connection transforms in the adjoint representation:
\[
r_h^* \omega = \Ad_h^{-1} \omega.
\]
\item
$\omega_e : T_e E \to \mathfrak{g}$ is a linear isomorphism at each
point $e \in E$.
\item
For each $A \in \mathfrak{g}$, define a vector field $\vec{A}$ on
$E$ 
(called a \emph{constant vector field})
by the equation $\vec{A} \hook \omega = A$. For $A \in
\mathfrak{h}$, the vector fields $\vec{A}$ generate the right
$H$-action:
\[
\vec{A}(e) = \left. \frac{d}{dt} r_{e^{tA}} e \right|_{t=0}, \text{ for all $e \in E$}.
\]
\end{enumerate}
\end{definition}

All statements in this paper hold equally true for real or
holomorphic Cartan geometries with obvious modifications.

\begin{example}
The bundle $G \to G/H$ is a $G/H$-geometry, called the \emph{model},
with Cartan connection $\omega=g^{-1} \, dg$ the left invariant
Maurer--Cartan 1-form on $G$.
\end{example}

\begin{example}\label{example:StandardFlatH}
Suppose that $\rho : H \to \GL{n,\R{}}$ is a representation,
and let $G = H \rtimes \R{n}$. 
Write elements of $G$ as ``matrices'',
\[
g = 
\begin{pmatrix}
h & x \\
0 & 1
\end{pmatrix},
\]
for $h \in H$ and $x \in \R{n}$.
Then the Maurer--Cartan 1-form is
\[
g^{-1} \, dg =
\begin{pmatrix}
h^{-1} \, dh & h^{-1} \, dx \\
0 & 0 
\end{pmatrix}.
\]
\end{example}

\begin{definition}
The \emph{universal constant vector field} of a Cartan geometry $E \to M$
is the vector field $X$ on $E \times \mathfrak{g}$ given by
$X(e,A)=\left(\vec{A},0\right)$. The \emph{constant flow}
of a Cartan geometry is the flow of the 
universal constant vector field.
\end{definition}

\begin{example}
Returning to example~\vref{example:StandardFlatH},
the constant flow is 
\[
e^{t X}
\left(
\begin{pmatrix}
h & x \\
0 & 1
\end{pmatrix},
\begin{pmatrix}
A & v \\
0 & 1
\end{pmatrix}
\right)
=
\left(
\begin{pmatrix}
h \, e^{tA} & x+h \left(\int e^{tA} \, dt\right) v  \\
0 & 1
\end{pmatrix},
\begin{pmatrix}
A & v \\
0 & 1
\end{pmatrix}
\right).
\]
The projected curves in $G/H=\R{n}$ are the flow lines of linear vector fields, up to translation.
In particular, if we take $G/H=\R{2}$, for various choices of $H$ we get:
\[
\begin{array}{cl}
H & \text{curves in $\R{2}$} \\
\midrule
\left\{1\right\} & \text{lines} \\
\SO{2} & \text{circles and lines} \\
\R{+} \SO{2} & \text{circles, lines and exponential spirals} \\
\end{array}
\]
\end{example}

\begin{definition}
We will say that a Cartan geometry is \emph{complete} if all of the constant vector fields are complete i.e.
if the flow of the universal constant vector field is complete. 
\end{definition}

\begin{example}
The model is complete, as the constant vector fields are the right invariant vector fields on $G$.
\end{example}

\begin{lemma}[Sharpe \cite{Sharpe:1997} p. 188, theorem 3.15]\label{lemma:TgtBundle}
The Cartan connection of any Cartan geometry $\pi : E \to M$
determines isomorphisms
\[
\xymatrix{%
0 \ar[r] & \ker \pi'(e) \ar[r] \ar[d] & T_e E \ar[r] \ar[d] & T_m M \ar[r] \ar[d] & 0 \\
0 \ar[r] & \mathfrak{h} \ar[r] & \mathfrak{g} \ar[r] &
\mathfrak{g}/\mathfrak{h} \ar[r] & 0 }
\]
for any points $m \in M$ and $e \in E_m$; thus
\[
TM = E \times_H \left( \mathfrak{g}/\mathfrak{h} \right ).
\]
\end{lemma}

\begin{definition}
Suppose that $E_0 \to M_0$ and $E_1 \to M_1$ are two
$G/H$-geometries, with Cartan connections $\omega_0$ and $\omega_1$,
and $X$ is manifold, perhaps with boundary and corners. A smooth map
$\phi_1 : X \to M_1$ is a \emph{development} of a smooth map $\phi_0
: X \to M_0$ if there exists a smooth isomorphism $\Phi : \phi_0^*
E_0 \to \phi_1^* E_1$ of principal $H$-bundles identifying the
1-forms $\omega_0$ with $\omega_1$.
\end{definition}

Development is an equivalence relation. For example, a development
of an open set is precisely a local isomorphism. The graph of $\Phi$
is an integral manifold of the Pfaffian system $\omega_0=\omega_1$
on $\phi_0^* E_0 \times E_1$, and so $\Phi$ is the solution of a
system of (determined or overdetermined) differential equations, and
conversely solutions to those equations determine developments. By
lemma~\vref{lemma:TgtBundle} the developing map $\phi_1$ has
differential $\phi_1'$ of the same rank as $\phi_0'$ at each point
of $X$.

\begin{definition}
Suppose that $E_0 \to M_0$ and $E_1 \to M_1$ are $G/H$-geometries.
Suppose that $\phi_1 : X \to M_1$ is a development of a smooth map
$\phi_0 : X \to M_0$ with isomorphism $\Phi : \phi_0^* E_0 \to
\phi_1^* E_1$. By analogy with Cartan's method of the moving frame,
we will call $e_0$ and $e_1$ \emph{frames} of the development if
$\Phi\left(e_0\right)=e_1$.
\end{definition}
\begin{definition}
Suppose that $E_0 \to M_0$ and $E_1 \to M_1$ are $G/H$-geometries.
We will say that $M_1$ \emph{rolls freely} on curves in $M_0$ to
mean that every curve $\phi_0 : C \to M_0$ has an unramified
covering which admits a development with any chosen frames $e_0 \in
\phi_0^* E_0$ and $e_1 \in E_1$.
\end{definition}
Note that we get to pick $e_1 \in E_1$ arbitrarily, and then we can construct
a development $\phi_1 : C \to M_1$ with $e_1$ lying in $\phi_1^* E_1$.
\begin{definition}
We can similarly define the notion of rolling freely on some family of
curves, for example immersed curves.
\end{definition}
\begin{remark}
We obtain the same development if we replace the frames by $e_0 h$
and $e_1 h$, for any $h \in H$. The curve mentioned could be a real
curve, or (if the manifolds and Cartan geometries are complex
analytic) a complex curve, with boundary and corners and arbitrary
complex analytic singularities. The frames determine the
development, and any frames yield a local development, because the
equation $\omega_0=\omega_1$ is a system of first order ordinary
differential equations in local coordinates.
\end{remark}
\begin{remark}
Ehresmann \cite{Ehresmann:1938} called a Cartan geometry \emph{normal} if it rolls freely on all curves in
its model. Kobayashi \cite{Kobayashi:1954} called such a Cartan geometry \emph{complete}. 
Sharpe \cite{Sharpe:1997} called a Cartan geometry \emph{complete} if it
has complete flow. Kobayashi \cite{Kobayashi:1954} claimed that normalcy and complete flow
are equivalent. (This would explain Sharpe's use of the term \emph{complete},
since he was aware of Kobayashi's paper.) Clifton \cite{Clifton:1966} proved
that normalcy implies complete flow, but he also gave explicit examples to show
that complete flow need not imply normalcy---these concepts are \emph{inequivalent},
i.e. Kobayashi's claim is wrong. Unfortunately this is not well known
in the Cartan geometry literature. I will use the word \emph{complete}
to mean complete flow, following Sharpe, since Sharpe's book is now
the standard reference work on Cartan connections.
There is no way to avoid the unfortunate confusion that fogs the literature.
\end{remark}

\begin{example}
A Riemannian geometry is a Cartan geometry modelled on Euclidean
space. As a Riemannian manifold, the unit sphere rolls freely on
curves in the plane. Rolling in this context has its obvious
intuitive meaning (see Sharpe \cite{Sharpe:1997} pp. 375--390). The
development of a geodesic on the plane is a geodesic on the sphere,
so a portion of a great circle. The upper half of the unit sphere
does \emph{not} roll freely on the plane (by our definition),
because if we draw any straight line segment in the plane of length
more than $\pi$, we can't develop all of it to the upper half of the
sphere. As we will see, every Riemannian manifold rolls
on curves in a given Riemannian manifold $M$ just when
$M$ is complete as a metric space.
\end{example}

\begin{theorem}\label{thm:Main}
Let $M$ be a manifold with a Cartan geometry. The following are equivalent:
\begin{enumerate}
%\item $M_1$ is complete. \label{item:complete}
\item $M$ rolls freely on curves in its model.
\label{item:rollsOnModel}
\item $M$ rolls freely on curves in any Cartan geometry with the same
model. \label{item:rollsOnAnything}
\item $M$ rolls freely on curves in some Cartan geometry with the same
model. \label{item:rollsOnSomething}
\end{enumerate}
\end{theorem}

\begin{remark}
Clifton \cite{Clifton:1966} gives examples of Cartan geometries rolling freely on 
all smooth immersed curves in the model, but not on all smooth curves.
\end{remark}

\begin{theorem}\label{theorem:Immersed}
Let $M_1$ be a manifold with a Cartan geometry. The following are equivalent:
\begin{enumerate}
%\item $M_1$ is complete. \label{item:complete}
\item $M_1$ rolls freely on immersed curves in its model.
\label{item:Immersed:rollsOnModel}
\item $M_1$ rolls freely on immersed curves in any Cartan geometry with the same
model. \label{item:Immersed:rollsOnAnything}
\item $M_1$ rolls freely on immersed curves in some Cartan geometry with the same
model. \label{item:Immersed:rollsOnSomething}
\end{enumerate}
\end{theorem}

This article begins by studying coframings; starting in section~\vref{section:MorphismsOfCartanGeometries},
we move on to study Cartan geometries. We define a notion of morphism of Cartan
geometries, and show that the main theorems on Cartan geometries have analogues
for morphisms. For example, the Frobenius--Gromov theorem is generalized to
the category of Cartan geometry morphisms. We also develop a theory
of morphism deformations, and of infinitesimal morphism deformations.
We employ this theory to prove completeness of some Cartan geometries.

This material is based upon works supported by the Science Foundation Ireland under Grant No. MATF634.

\section{Coframings}

Lets start with a simpler concept than a Cartan geometry.

\subsection{Definitions}

\begin{definition}
A \emph{coframing} 
of a manifold $M$ is a 1-form $\omega \in \nForms{1}{M} \otimes V$,
where $V$ is a vector space of the same dimension as $M$, so that
at each point $m \in M$, the linear map $\omega_m : T_m M \to V$
is a linear isomorphism.
A coframing is also called a \emph{parallelism}
or an \emph{absolute parallelism}.
\end{definition}

Clearly a coframing is a Cartan geometry modelled
on $G/H=V/0$.

\begin{remark}
There are obvious complex analytic analogues
of coframings,
which we can call \emph{holomorphic coframings}.
Many of the results we prove on coframings
have obvious complex analytic analogues.
\end{remark}

\begin{definition}
If two manifolds $M_0$ and $M_1$ have coframings
$\omega_0 \in \nForms{1}{M_0} \otimes V_0$ 
and $\omega_1 \in \nForms{1}{M_1} \otimes V_1$
respectively, and $A : V_0 \to V_1$ is a linear map,
then a smooth map $f_1 : X \to M_1$
is an $A$-\emph{development} of 
a smooth map $f_0 : X \to M_0$
if $f_1^* \omega_1 = A f_0^* \omega_0$.
We will also say that $f_0$ $A$-\emph{develops} to $f_1$,
or just say that $f_0$ develops to $f_1$
if $A$ is understood. If $A$ is not specified,
we will assume that $A=I$ is the identity map.
\end{definition}

\begin{definition}
We will say that \emph{curves $A$-develop} 
from a manifold $M_0$ with a coframing to a 
manifold $M_1$ with a coframing
if every rectifiable map $\R{} \to M_0$ $A$-develops to a map $\R{} \to M_1$.
\end{definition}

\begin{definition}
We will say that \emph{complex curves $A$-develop}
from a complex manifold $M_0$ with a holomorphic
coframing to a complex manifold $M_1$ with a holomorphic
coframing if every holomorphic map $f : C \to M_0$ from a simply
connected Riemann surface develops to a holomorphic map to $M_1$.
\end{definition}

\begin{definition}
If $\omega \in \nForms{1}{M} \otimes V$
is a coframing, and $v \in V$ is a vector,
we denote by $\vec{v}$ the unique
vector field on $M$ so that
$\vec{v} \hook \omega = v$,
called a \emph{constant} vector field.
\end{definition}

\begin{definition}
We will say that \emph{lines develop} to a manifold with coframing,
or that the coframing is \emph{complete}, or has \emph{complete flow} if each 
constant vector field complete, i.e. its flow is defined for all time.
\end{definition}

\begin{definition}
We say similarly that \emph{complex lines develop} to a complex manifold 
with holomorphic coframing, of that the coframing is holomorphically
complete, or has \emph{complete holomorphic flow} if each holomorphic constant
vector field is holomorphically complete, i.e. its flow
is defined for all complex time.
\end{definition}

\begin{definition}
Pick a positive definite inner product on a vector space $V$.
The associated \emph{canonical metric} of a coframing $\omega$
is the Riemannian metric for which the $\omega_m : T_m M \to V$ is an
orthogonal linear map for all $m \in M$.
\end{definition}

\begin{definition}
Pick a Hermitian inner product on a complex vector space $V$.
The \emph{canonical Hermitian metric} of a 
holomorphic coframing is the Hermitian metric for which the
$\omega_m : T_m M \to V$ is unitary for each $m \in M$.
\end{definition}

\begin{example}
Pick a Lie group $G$. The left (right) invariant Maurer--Cartan 1-form
is a left (right) invariant coframing. Lines develop to left translates of 1-parameter
subgroups. Curves develop from any Lie group to any other, by solving a finite dimensional Lie
equation, as we will see. We will also see that every canonical metric on
any Lie group is complete.
\end{example}

\begin{remark}
 The study of completeness of canonical metrics is a small
part of the larger picture of Charles Frances \cite{Frances:2008},
in which the completion of the metric space in the canonical
metric plays an essential role.
\end{remark}

\begin{definition}
Given a coframing $\omega \in \nForms{1}{M} \otimes V$,
the \emph{torsion} of $\omega$ is the quantity $T : M \to V \otimes \Lm{2}{V}^*$
given by $d \omega = T \omega \wedge \omega$.
\end{definition}

It is easy to prove:
\begin{lemma}
A coframing is locally isomorphic to the left invariant Maurer--Cartan
coframing of a Lie group just when its torsion is constant.
A coframing is isomorphic to the left invariant Maurer--Cartan coframing of a Lie group
just when lines develop to the coframing and the torsion is constant.
\end{lemma}

\subsection{Flows and coframings}

\begin{remark}
Let's picture a coframing to which lines develop but curves do not.
Take the plane with a coframing of the form 
\[
\omega= 
\frac{1}{f}
\begin{pmatrix}
dx \\
dy
\end{pmatrix}
\]
for some function $f > 0$.
Draw a parabola. Choose the function
$f$ to grow rapidly along the parabola,
but rapidly approach $1$ away
from the parabola. Lines develop
to this coframing to become lines, but reparameterized
to move more quickly as they approach
the parabola. As we travel along any line, we eventually
move away from the parabola, so
eventually go back to a slow steady speed. 
However, moving along
the parabola at unit speed in the
framing means at speed $f$ in the
usual coframing $\binom{dx}{dy}$,
so in some finite time we leave any compact set.
As we will see, the canonical metric is incomplete.
\end{remark}

\begin{lemma}[Clifton \cite{Clifton:1966}]\label{lemma:MetricIFFdevelopRn}
Let $M$ be an $n$-dimensional manifold with coframing. 
Every canonical metric of $M$ is complete just when curves develop from $\R{n}$ to $M$.
\end{lemma}
\begin{proof}
Take a coframing $\omega \in \nForms{1}{M} \otimes V$ on a manifold $M$.
Take a rectifiable curve $x(t) \in V$.
Let $v(t)=x'(t)$. Suppose that the canonical metric is
complete. Locally integrate the flow of $\vec{v}(t)$ through any
point. In the canonical metric, the flow lines are rectifiable, with
locally bounded velocity, so if they are defined on an open time interval,
then they extend uniquely to the closure of that interval. 
The flow line through each point is defined on a maximal
open interval, which must therefore be the entire real line.

Take a coframing $\omega \in \nForms{1}{M} \otimes V$
on a manifold $M$. Suppose that curves develop from $V$ to $M$.
Take a geodesic $x(t) \in M$ of a canonical metric
on $M$, parameterized by arc length. Suppose that $x(t)$ does not extend beyond
some open interval of time. Then the development to $V$ is 
\[
y(t) = \int x'(t) \hook \omega \, dt.
\]
Because $x(t)$ is parameterized by arc length, $\left|dy/dt\right|=1$ for all $t$.
This curve $y(t)$ is rectifiable, since $dy/dt$ is bounded
and integrable. Extend $y(t)$ to a rectifiable curve on the
closure of the open interval on which $x(t)$ is defined, by continuity. Develop $y(t)$ back
to $M$ to extend $x(t)$. 
\end{proof}

\begin{example}
Curves develop to a coframing on a manifold just when they develop to every other coframing which
agrees with it outside of some compact set.
\end{example}

\begin{example}
On the connected sum of two manifolds with coframings, pick a coframing which
agrees with the original coframings on the summands away from a compact
set where the gluing takes place. Curves develop to the connected sum 
just when they develop to both summands. The same idea works for developing lines.
\end{example}

\begin{lemma}
Let $M$ be a manifold with coframing $\omega \in \nForms{1}{M} \otimes V$.
If a canonical metric on $M$ is complete then, for any integrable function 
$f : [a,b] \to V$, the time varying vector field $\vec{f}(t)$ on $M$ is complete.
\end{lemma}
\begin{proof}
Clearly any flow line of $\vec{f}(t)$ has velocity at most $|f(t)|$.
Let $x(t)$ be a flow line. The distance of $x(t)$ from $x_0$ increases
at a rate of at most $\int |f(t)| \, dt$, and therefore for each time $t$, the
integral curve remains in a closed ball. Because the metric is complete,
closed balls are compact. Suppose that the flow line is defined for time $0 \le t < T$.
The limit $\lim_{t \to T} x(t)$ must exist because of the compactness
of the ball. Therefore we can extend $x(t)$ continuously
to $0 \le t \le T$. But then $\dot{x}(t)$ also extends
to be locally integrable on $0 \le t \le T$, as $\vec{f}(t)$.
Therefore $x(t)$ is rectifiable on $0 \le t \le T$ and a flow line of $\vec{f}(t)$.
But $x(t)$ is defined on a maximal open interval, so must be
defined on $[a,b]$.
\end{proof}

\begin{lemma}\label{lemma:RnImpliesN}
Let $M_1$ be a manifold with coframing $\omega_1 \in \nForms{1}{M_1} \otimes V_1$.
If curves develop from $\R{n}$ to $M_1$ then, for any manifold $M_0$ with 
coframing $\omega_0 \in \nForms{1}{M_0} \otimes V_0$ 
and any linear map $A : V_0 \to V_1$, curves $A$-develop from $M_0$ to $M_1$. 
\end{lemma}
\begin{proof}
Suppose that curves develop from $V_1$ to $M_1$. Then $A$-develop
curves from $M_0$ to $V_1$, and then from $V_1$ to $M_1$. 
\end{proof}

\begin{lemma}\label{lemma:TVVF}
Take manifolds $M_0$ and $M_1$ with coframings $\omega_0 \in \nForms{1}{M_0} \otimes V_0$ and
$\omega_1 \in \nForms{1}{M_1} \otimes V_1$. Take a linear map $A : V_0 \to V_1$.
Let $W_1 \subset V_1$ be the image of $A$. Curves $A$-develop from $M_0$ to $M_1$ just when,
for any integrable function $f : I \to W_1$ ($I \subset \R{}$ any interval),
the time varying vector field $\vec{f}(t)$ on $M_1$ is complete.
\end{lemma}
\begin{proof}
Suppose that all of these time varying vector fields $\vec{f}(t)$ are complete.
Take any rectifiable curve $m_0(t)$ on $M_0$, and let $f(t) = m_0'(t) \hook \omega_0$.
Construct the time varying vector field $\overrightarrow{Af}(t)$;
its flow line is our development.

Conversely, suppose that curves $A$-develop from $M_0$ to $M_1$.
Pick an integrable function $f_1 : I \to W_1$. We can pick a subspace $W_0 \subset V_0$
complementary to the kernel of $A$, so that $A : W_0 \to W_1$ is a linear isomorphism. We can 
let $f_0 = A^{-1} f_1 : I \to W_0$. Pick points $m_0 \in M_0$ and $m_1 \in M_1$
and a time $t_0 \in I$. The time-varying vector field $\overrightarrow{f_0}(t)$ on $M_0$
has a flow line $x_0(t)$ through with $x_0\left(t_0\right)=m_0$, and the time-varying vector
field $\overrightarrow{f_1}(t)$ on $M_1$ has a flow line $x_1(t)$ through 
with $x_1\left(t_0\right)=m_1$. Clearly these are developments 
of one another, over the time interval during which both curves are defined.

There is a time (possibly infinite) $T$, $t_0 < T \le \infty$, 
when the flow line $x_0(t)$ stops being defined. Lets write
this time $T$ as $T=t_0+\Delta t\left(m_0,t_0\right)$.
Similarly define $\Delta t\left(m_1,t_0\right)$ for any point $m_1 \in M_1$.

Since $x_1(t)$ is a development of $x_0(t)$, where both are defined,
we can extend $x_1(t)$ to be defined for as long a time as $x_0(t)$ is defined:
extend it to be the development of $x_0(t)$.
Then $x_1(t)$ is still a flow line of $\overrightarrow{f_1}(t)$.
So flow lines on $M_1$ are defined
for at least as long a time interval as those on $M_0$:
$\Delta t\left(m_0,t_0\right) \ge \Delta t\left(m_1,t_0\right)$
for any points $m_0 \in M_0$ and $m_1 \in M_1$. In particular, 
$\inf_{m_1 \in M_1} \Delta t\left(m_1,t_0\right) \ge 
\sup_{m_0 \in M_0} \Delta t\left(m_0,t_0\right) > 0$.
Therefore $\Delta t\left(m_1,t_0\right)$ is bounded
from below by a positive function of $t_0$.

Pick inner products on $V_0$ and $V_1$ so that 
$A : W_0 \to W_1$ is an isometry. Equip
$M_0$ and $M_1$ with the canonical
metrics. Fix a point $m_0 \in M_0$. Pick a number $r>0$
small enough so that the ball of radius
$r$ about $m_0$ (in the canonical metric) is compact.
The function
\[
\int_{t_0}^{t_1} \left|f_1(t)\right| \, dt
\]
is a continuous function of a variable $t_1$,
where defined. It vanishes at $t_1=t_0$.
Pick $t_1 > t_0$ small enough so that
\[
\int_{t_0}^{t_1} \left|f_1(t)\right| \, dt < r.
\]
The value of $t_1$ you choose may depend on the choice of point $m_0$ and time $t_1$.

At time $t_1>t_0$, the flow line of $\overrightarrow{f_0}(t)$ from 
$m_0$ and time $t_0$ stays in the ball
of radius at most $\int_{t_0}^{t_1} \left|f_0\right| \, dt$
around $m_0$ in $M_0$. So the flow line stays inside this
compact ball of radius $r$.
Therefore the flow line $m_0(t)$ must continue to be defined up to time $t=t_1$. So if
\[
\int_{t_0}^{t_1} \left|f_1(t)\right| \, dt \le r,
\]
then 
\[
\inf_{m_1 \in M_1} \Delta t\left(m_1,t_0\right) \ge t_1-t_0.
\]

Divide up the time axis into small increments $t_0, t_1, t_2, \dots, t_N$ so that 
\[
\int_{t_i}^{t_{i+1}} |f(t)| \, dt < r
\]
on each increment. The flow line $m_1(t)$ of $\overrightarrow{f_1}(t)$
will survive for at least enough time to get from each time $t_i$ to the 
time $t_{i+1}$, and therefore by induction the flow survives for all time.
\end{proof}

Summing up:
\begin{proposition}\label{proposition:Main}
On an $n$-dimensional manifold $M$ with coframing, the following
are equivalent:
\begin{enumerate}
 \item some canonical metric is complete,
 \item every canonical metric is complete,
 \item curves develop from $\R{n}$,
 \item curves $I$-develop from some manifold with coframing,
 \item curves $A$-develop from any manifold with coframing for any $A$.
\end{enumerate}
\end{proposition}
\begin{proof}
Clearly (1) is equivalent to (2), since the metrics are Lipschitz equivalent.
By lemma~\vref{lemma:MetricIFFdevelopRn}, (2) is equivalent to (3).
Clearly (3) implies (4) and (5) implies (4).
By lemma~\vref{lemma:RnImpliesN}, (3) implies (5).
Suppose that curves $I$-develop to $M$ from some manifold $M_0$
with coframing. By lemma~\vref{lemma:TVVF}, (4) implies that bounded time
varying vector fields of a suitable form on $M$ are complete,
which is independent of choice of manifold $M_0$ from
which we were able to develop curves, and so (4) implies (5).
\end{proof}

\begin{example}
In the plane, the coframing
\[
\omega =
\begin{pmatrix}
 \frac{1}{2 + \sin\left(x^2 y\right)} dx \\
 dy
\end{pmatrix}
\]
%\begin{align*}
%X &= \left(2 + \sin\left(x^2 y\right) \right) \partial_x \\
%Y &= \partial_y
%\end{align*}
clearly has complete canonical metric, since the metric is within bounded
dilation of the Euclidean metric. The dual framing is
\begin{align*}
 \vec{\partial}_1 &= \left(2+\sin\left(x^2 y\right)\right) \partial_x, \\
 \vec{\partial}_2 &= \partial_y. \\
\end{align*}
Note that
$\left[\vec{\partial}_1,\vec{\partial}_2\right]=-\cos\left(x^2 y\right) x^2 \, \partial_1$.
Along $y=0$, $\left[\vec{\partial}_1,\vec{\partial}_2\right]=-x^2 \, \partial_x$ is
incomplete. Therefore completeness of the canonical metric does not ensure
that the brackets of vector fields from a coframing are complete.
\end{example}

\begin{conjecture}
If the Lie algebra of vector fields generated by a coframing consists
entirely of complete vector fields, then the canonical metric of the coframing is complete.
\end{conjecture}

\section{Immersed curves}

\begin{definition}
A \emph{ordinary vector field} on a manifold $M$ with coframing $\omega \in \nForms{1}{M} \otimes V$ is a time varying vector field
$\vec{f}(t)$ on $M$ with $f : I \to V$ bounded, continuous and never vanishing on an interval $I \subset \R{}$.
\end{definition}

\begin{lemma}\label{lemma:TVVFNZ}
Take two manifolds $M_0$ and $M_1$ with coframings $\omega_0 \in \nForms{1}{M_0} \otimes V_0$ and
$\omega_1 \in \nForms{1}{M_1} \otimes V_1$. Take a linear map $A: V_0 \to V_1$, with image $W_1 \subset V_1$.
Immersed curves $A$-develop from $M_0$ to $M_1$ just when all ordinary vector fields 
$\vec{f}(t)$ on $M_1$ are complete, for any $f : I \to W_1$.
\end{lemma}
\begin{proof}
The proof is identical to that of lemma~\vref{lemma:TVVF}, but with $f$ nowhere $0$. 
\end{proof}

\begin{corollary}
On an $n$-dimensional manifold $M$ with coframing, the following
are equivalent:
\begin{enumerate}
 \item immersed curves develop from $\R{n}$,
 \item immersed curves develop from some $n$-dimensional manifold with coframing,
 \item immersed curves $A$-develop from any manifold with coframing,
where $A$ is any injective linear map.
\end{enumerate}
\end{corollary}

\begin{remark}
Clifton \cite{Clifton:1966} has an example of a coframing so that lines develop from $\R{n}$ but some immersed curves don't.
\end{remark}

\section{Holomorphic coframings}

\begin{lemma}\label{lemma:MetricIFFdevelopCn}
Let $M_1$ be an $n$-complex-dimensional complex manifold with 
holomorphic coframing $\omega_1 \in \nForms{1}{M_1} \otimes V_1$.
Take a linear map $A : V_0 \to V_1$. If the canonical metric of $M_1$
is complete then complex curves $A$-develop from $V_0$ to $M_1$.
\end{lemma}
\begin{remark}
 It is not known whether there are holomorphic
coframings with incomplete canonical metric to which all 
complex curves develop from a complex vector space.
\end{remark}
\begin{proof}
Suppose that the canonical metric is complete.
Then every real curve in $V_0$ $A$-develops
to $M_1$ by lemma~\vref{lemma:MetricIFFdevelopRn}.
Suppose that $C$ is a simply connected Riemann
surface and that $f : C \to V_0$ is a holomorphic
map. The equations for $A$-development form a complex
analytic system of determined ordinary differential equations
of first order, so there is a unique local solution
near any point. Take any real curve on $C$, and develop
it. The resulting development extends to a holomorphic
development of a neighborhood of the real curve.
Because $C$ is simply connected, there is no monodromy,
and these local solutions extend to all of $C$.
\end{proof}

\begin{lemma}\label{lemma:HolVecFieldsAndDevelopments}
Let $M_0$ and $M_1$ be complex manifolds with 
holomorphic coframings $\omega_0 \in \nForms{1}{M_0} \otimes V_0$
and
$\omega_1 \in \nForms{1}{M_1} \otimes V_1$.
Let $A : V_0 \to V_1$ be a linear map.
Let $W_1$ be the image of $A$.
Let $U$ be a simply connected open set in $\C{}$.
If complex curves $A$-develop from $M_0$ to $M_1$, then
every time varying vector field $\vec{f}(t)$ on $M_1$ is complete,
with $f : U \to W_1$ any holomorphic function.
\end{lemma}
\begin{proof}
Suppose that complex
curves $A$-develop from $M_0$ to $M_1$.  
Let $\Delta_r(t)$ be the disk of radius $r$ about any point
$t \in \C{}$. Take a holomorphic 
function $f_1 : U \to W_1$. Let
$W_0 \subset V_0$ be any complex linear complement to the
kernel of $A$. Let $f_0 = A^{-1} f_1$.
For every point $m_0 \in M$ and initial time $t_0 \in U$,
there is some radius $r=r\left(m_0, t_0\right)$
so that the flow of $\vec{f}_0(t)$ starting at $m_0$ at time 
$t=t_0$ is defined for all $t \in \Delta_r\left(t_0\right)$,
and $0 < r \le \infty$.
Similarly define $r\left(m_1,t_0\right)$ for $m_1 \in M_1$.

Since we can develop curves from $M_0$, we can
develop the flow of $\vec{f}_0(t)$ through
any point $m_0$ into a flow line of $\vec{f}_1(t)$
through any point $m_1$. 
Therefore $r\left(m_1,t_0\right) \ge r\left(m_0,t_0\right)$.
Since this holds for any points $m_0 \in M_0$ and
$m_1 \in M_1$,
\[
\inf_{m_1 \in M_1} r\left(m_1,t_0\right) \ge \sup_{m_0 \in M_0} r\left(m_0,t_0\right).
\]

% Let $x_1(t)$ be the flow line of $\overrightarrow{Af}(t)$ so
% that $x_1\left(t_0\right)=m_1$.
% Consider a point $t_1 \in \C{}$ within distance
% $r\left(m_1,t_0\right)$ of $t_0$. 
% Since $x_1(t)$ is defined for any $t$ for which 
% $\left|t-t_0\right|<r\left(m_1,t_0\right)$,
% we see that $x_1(t)$ is defined for 
% any $t$ for which
% $\left|t-t_1\right|< r\left(m_1,t_0\right)-\left|t_1-t_0\right|$.
% Therefore
% \[
% r\left(x_1\left(t_1\right),t_1\right) \ge r\left(m_1,t_0\right)-\left|t_1-t_0\right|.
% \]
% So for any time $t$,
% \[
% r\left(x_1\left(t\right),t\right) \ge r\left(m_1,t_0\right)-\left|t-t_0\right|.
% \]

Pick some point $t'$ on the boundary of $\Delta_{r\left(m_1,t_0\right)}$,
so that the flow line $x_1(t)$ does not extend to $t=t'$,
and let $t \to t'$. Then
\[
\lim_{t \to t'} r\left(x_1\left(t\right),t\right) \to 0.
\]
Therefore
\begin{equation}\label{equation:limsup}
\lim_{t \to t'} \sup_{m_0 \in M_0} r\left(m_0,t\right) \to 0.
\end{equation}

If $t' \in U$ then $\vec{f}_0(t)$ is holomorphic
at times $t$ near $t'$ so there is a unique flow line $m_0(t)$
of $\vec{f}_0(t)$ with $m_0\left(t'\right)=m_0$,
defined in some disk of positive radius
$r\left(m_0,t'\right)>0$. Clearly
this contradicts equation~\ref{equation:limsup} above.
Therefore $t'$ does not belong to $U$.

Take any real curve $t(s)$ inside $U$ for some real variable $a \le s \le b$.
Let $x_1(t)$ be the integral curve of $\vec{f}_1(t)$ with 
$x_1\left(t\left(a\right)\right)=m_1$.
Then $x_1(t)$ is defined for 
all $t$ in a disk of radius $r\left(m_1,t(a)\right)$.
But then as we vary $s$, at the first moment 
$s=s_1$ that $t(s)$ 
leaves this disk, it doesn't leave $U$, 
so we can extend this flow line $m_1(t)$
to values of $t$ near $t(s)$. We
can then extend to the disk of radius
$r\left(m_1,t\left(s_1\right)\right)$,
etc. In this way we produce a sequence
of extensions of the flow line
to various open intervals, and then
extend to their closures, so
clearly to all $s$ from $s=a$ to $s=b$.

Therefore the flow line of $\vec{f}_1(t)$ is defined for all $t$ along
any real curve inside $U$. Since $U$ is simply connected, and the
local flow lines are unique, there is a unique flow line defined
on all of $U$.
\end{proof}

\begin{lemma}\label{lemma:SomeImpliesCn}
Let $M_1$ be a complex manifold with 
holomorphic coframing $\omega_0 \in \nForms{1}{M_1} \otimes V$. If
complex curves $I$-develop from \emph{some} 
holomorphically coframed manifold to $M_1$ 
then complex curves develop from $V$ to $M_1$.
\end{lemma}
\begin{proof}
Suppose that
$M_0$ is a complex manifold with a holomorphic coframing
$\omega_0 \in \nForms{1}{M_0} \otimes V$. Suppose that complex
curves develop from $M_0$ to $M_1$.  
Then by lemma~\vref{lemma:HolVecFieldsAndDevelopments},
all holomorphic time varying vector fields $\vec{f}(t)$
on $M$ are complete, for $f : U \to V$ a holomorphic
function, where $U \subset \C{}$ is any
simply connected open set.

If we take a complex curve
\[
f : C \to \C{n},
\]
with $C$ a disk, or $C=\C{}$, then clearly
the developments are the flow
lines of $\vec{f}'(t)$.
Therefore we only need to consider the
problem for $C=\CP{1}$, the Riemann sphere.
But for $\CP{1}$, any holomorphic
map $f$ is constant, and so the development
is just a constant map.
\end{proof}

\begin{proposition}\label{proposition:ComplexRolling}
Let $M$ be an $n$-complex-dimensional complex manifold with 
holomorphic coframing. Then the following are equivalent:
\begin{enumerate}
 \item complex curves $I$-develop to $M$ from $\C{n}$,
\item complex curves $I$-develop to $M$ from \emph{some} $n$-dimensional
complex manifold with holomorphic coframing,
\item complex curves $A$-develop to $M$ from \emph{any} holomorphically framed
complex manifold, for any linear map $A$.
\item
time varying vector fields $\vec{f}(t)$ on $M$
are complete, when $f(t)$ is any holomorphic function of $t$,
and $t$ lives in some simply connected open set in $\C{}$.
\end{enumerate}
\end{proposition}
\begin{proof}
Lemma~\vref{lemma:HolVecFieldsAndDevelopments} shows that (2) implies (4).
Lemma~\vref{lemma:SomeImpliesCn} shows that (2) implies (1).
Clearly (1) implies (2) and (3) implies (2). 
Clearly (1) implies (4), by developing the curve $\int f(t) \, dt$ from $\C{n}$.
If we assume (1), we can $A$-develop curves from any complex manifold $M_0$ 
to $\C{n}$ first and then develop to $M$. Therefore (1) implies (3).
To finish, lets see why (4) implies (1). It is clear
that for any simply connected complex curve $C$ biholomorphic to the disk or $\C{}$,
any map $f : C \to \C{n} $
develops to $M$, by integrating $X_t = \vec{f}'(t)$.
But then if $C$ is biholomorphic to $\CP{1}$, $f$
must be constant, so develops trivially from $\C{n}$
to $M$.
\end{proof}

\section{Symmetries of coframings}

\begin{definition}
 A \emph{symmetry} of a coframing
$\omega$ on a manifold
$M$ is a diffeomorphism $F : M \to M$
so that $F^* \omega=\omega$.
An \emph{infinitesimal symmetry}
of a coframing $\omega$ is a vector field
$Y$ so that
$0=\LieDer_Y \omega=0$.
\end{definition}

\begin{proposition}\label{proposition:CompleteInfinitesimalSymmetries}
If lines develop to a coframing, then any infinitesimal symmetry of that coframing
is a complete vector field. 
\end{proposition}
\begin{proof}
Suppose that $Y$ is an infinitesimal symmetry of a coframing
$\omega \in \nForms{1}{M} \otimes V$ on a manifold $M$. Then
$0=\left[Y,\vec{v}\right]$ for each $v \in V$.
The flow lines of $Y$ are permuted
by the flows of the various $\vec{v}$. We can therefore
take any flow line of $Y$, say starting
at a point $m$ and defined for a time $T>0$,
and flow it around, moving the point $m$
to any other nearby point. Suppose that the constant
vector fields $\vec{v}$ of $\omega$ are complete.
We can move the entire flow line to start at any
other nearby point, and still be defined
for time at least $T$. Therefore 
all flow lines through all nearby points 
are defined for  at least time $T$.
This means that as we follow any
flow line of $Y$, at each point
of that flow line, there is still
at least $T$ units of time left
to continue flowing. Since $0<T\le \infty$,
we never run out of time, i.e.
the flow lines are defined for all time.
\end{proof}

\begin{theorem}[Kobayashi \cite{Kobayashi:1995} p. 15, theorem 3.2]
The symmetries of a coframing
form a finite dimensional Lie group.
The Lie algebra of this Lie group
is precisely the set of complete
infinitesimal symmetries.
\end{theorem}

\begin{lemma}
A coframing on a connected manifold is homogeneous just when it has 
a finite collection of complete vector fields which are local infinitesimal symmetries,
and are linearly independent at every point.
\end{lemma}
\begin{proof}
By Kobayashi's theorem, the infinitesimal symmetries form a finite dimensional Lie algebra.
This Lie algebra will be spanned by the finite collection of vector fields 
By a theorem of Palais \cite{Palais:1957}, a finite dimensional Lie algebra of vector
fields, generated by complete vector fields, integrates to a Lie algebra action.
The action is locally transitive, because the vector fields are linearly independent.
Every orbit is open. The manifold is connected. Therefore the manifold is homogeneous
under the group action.
\end{proof}

\section{Morphisms of coframings}

\begin{definition}
Suppose that $M_0$ and $M_1$ are manifolds with coframings $\omega_0 \in \nForms{1}{M_0} \otimes V_0$
and $\omega_1 \in \nForms{1}{M_1} \otimes V_1$. Pick a linear map $A : V_0 \to V_1$.
An $A$-\emph{morphism} is a map $f : M_0 \to M_1$ so that $f^* \omega_1=A \omega_0$.
We will also refer to an $A$-morphism as just a morphism to avoid reference to a specific map $A$.
\end{definition}

\begin{example}
Local isomorphisms are precisely $I$-morphisms.
\end{example}

\begin{example}
 The flow lines of the vector fields $\vec{v}$ are morphisms: take
any manifold $M_1$ with coframing $\omega_1 \in \nForms{1}{M_1} \otimes V_1$. 
Pick an element $v_1 \in V_1$. Let $A : \R{} \to V_1$, $A(x)=x \, v_1$.
Pick a point $m_1 \in M_1$. Let $M_0$ be the open interval of
$\R{}$ of times $t$ for which $e^{\vec{v}_1} m_1$ is defined.
Let $\omega_0=dt$. Map $f : M_0 \to M_1$, $f(t)=e^{t \vec{v}} m_1$.
\end{example}

The sense in which development generalizes morphism is obvious.
\begin{proposition}\label{proposition:DevelopmentToMorphism}
Suppose that $M_0$ and $M_1$ are manifolds bearing coframings $\omega_0 \in \nForms{1}{M_0} \otimes V_0$ and
$\omega_1 \in \nForms{1}{M_1} \otimes V_1$. Suppose that $M_0$ is path connected. Take a linear map $A : V_0 \to V_1$.
Pick points $m_0 \in M_0$ and $m_1 \in M_1$. There is an $A$-morphism $f : M_0 \to M_1$ 
with $f\left(m_0\right)=m_1$ just when every loop in $M_0$ starting and ending at $m_0$
develops to a loop in $M_1$ starting and ending at $m_1$.
\end{proposition}
\begin{proof}
Clearly if there is an $A$-morphism $f : M_0 \to M_1$,
then every map $X \to M_0$ composes with $f$
to produce a development to a map $X \to M_1$.

On the other hand, suppose that all loops in $M_0$,
starting and ending at $m_0$, 
develop to loops in $M_1$
starting and ending at $m_1$. 
It is immediately clear (as in the 
theory of fundamental groups)
that there is a unique
smooth map $f : M_0 \to M_1$
so that $f\left(m_0\right)=m_1$
and so that development
of paths from $M_0$ starting
at $m_0$ is carried out by
composition with $f$.
It is then clear that $f$
is a morphism, since there is a path
on $M_0$ leaving $m_0$ with any
specified tangent vector at
its other end.
\end{proof}

\begin{lemma}\label{lemma:LocalMorphismOfFramings}
Suppose that $M_0$ and $M_1$ are manifolds bearing coframings $\omega_0 \in \nForms{1}{M_0} \otimes V_0$ and
$\omega_1 \in \nForms{1}{M_1} \otimes V_1$. Take a linear map $A : V_0 \to V_1$.
A submanifold $Z$ of $M_0 \times M_1$ is locally the graph of a local $A$-morphism just when 
\begin{enumerate}
\item $Z$ is an integral manifold of the Pfaffian system $\omega_1 = A \omega_0$ and
\item $\omega_0$, pulled back to $M_0 \times M_1$, is a coframing on $Z$.
\end{enumerate}
\end{lemma}
\begin{proof}
Because the 1-form $\omega_0$ is a coframing on $M_0$, its pullback to $Z$
is a coframing just when the composition $Z \subset M_0 \times M_1 \to M_0$
is a local diffeomorphism. Since the problem is local, we can assume
that this composition is a diffeomorphism. So we can assume that $Z$ is the graph
of a map $F : M_0 \to M_1$. Because $\omega_1=A \omega_0$ on $Z$, we see that
$F^*\omega_1=A \omega_0$. 
\end{proof}

\begin{lemma}
If $f : \text{open } \subset M_0 \to M_1$ is an $A$-morphism, 
and $m_0 \in M_0$ and $m_1=f\left(m_0\right)$,
then in logarithmic coordinates, 
$x=\log_{m_0}$ and $y=\log_{m_1}$,
the map $f$ is expressed as the linear map $y=Ax$.
\end{lemma}
\begin{proof}
Clearly if there were going to be a morphism
$f : U_0 \to M_1$, it would have to satisfy
\[
f\left(e^{\vec{v}} m_0\right) = e^{\overrightarrow{Av}}m_1.
\]
But this equation determines a map $f$ completely, near $m_0$:
\[
f(m)=e^{\overrightarrow{A \log_{m_0} m}} m_1.
\]
Clearly in logarithmic coordinates, $f$ is the linear map $A$. 
\end{proof}

\begin{corollary}
Suppose that $f_n : U_n \text{ open } \subset M_0 \to M_1$ is a sequence of $A$-morphisms,
say taking points $m_0(n) \in M_0$ to $m_1(n) \in M_1$.
If the sequences of points converge, say $m_0(n) \to m_0$
and $m_1(n) \to m_1$, then $f_n$ converges with all derivatives to a 
morphism $f : \text{ open } \subset M_0 \to M_1$ in some
neighborhood of $m_0$.
\end{corollary}
\begin{proof}
Clearly we can define a map $f$ by making $f$ be $A$ in logarithmic coordinates.
By the continuity in Picard's theorem, $f$ is the limit of the various
$f_n$ on some open set. Since the various $f_n$ all satisfy $f_n^* \omega_1 = A \omega_0$,
so does $f$.
\end{proof}

\begin{theorem}\label{theorem:GlobalMorphismOfFramings}
Suppose that $M_0$ and $M_1$ are real analytic manifolds bearing real analytic coframings 
$\omega_0 \in \nForms{1}{M_0} \otimes V_0$ and $\omega_1 \in \nForms{1}{M_1} \otimes V_1$.
Suppose that $M_0$ is connected. Take a linear map $A : V_0 \to V_1$.
Suppose that for each $v_0 \in V_0$, the vector field $\overrightarrow{Av}$
on $M_1$ is complete. 

Then every $A$-morphism of coframings $F : U_0 \to M_1$ 
defined on a connected open subset $U_0 \subset M_0$
factors uniquely as $F = \tilde{F} \iota$,
$\tilde{F} : \tilde{M}_0 \to M_1$ 
is an $A$-morphism of coframings from some covering space of $M_0$,
and $\iota : U_0 \to \tilde{M}_0$ is an injective local diffeomorphism
covering the inclusion $U_0 \to M_0$.
The covering map $p : \tilde{M}_0 \to M_0$
yields a map $\left(p,F\right) : \tilde{M}_0 \to M_0 \times M_1$,
making $\tilde{M}_0$ a immersed submanifold of $M_0 \times M_1$.
\end{theorem}
\begin{proof}
Suppose that $M'_0 \subset M_0 \times M_1$ is the graph
of an $A$-morphism on a connected open set of $M_0$. Denote by $\vec{v}_0$ 
the constant vector fields
on $M_0$ for $v_0 \in V_0$, and 
by $\vec{v}_1$ the constant
vector fields on $M_1$ for $v_1 \in V_1$.
On $M_0 \times M_1$,
the flow of the vector field $\hat{v}_0=\vec{v}_0 + \overrightarrow{Av_0}$ is tangent
to $M_0'$ for any vector $v_0 \in V_0$. Moreover, each tangent
plane of $M_0'$ is precisely the span
of all of these vector fields $\hat{v}_0$.

At each point of $M_0'$, for any vectors $v_0, w_0 \in V_0$ the bracket
\(
\left[
\hat{v}_0,\hat{w}_0
\right]
\)
lies tangent to $M_0'$, since the vector fields do. So we can express this bracket on $M_0'$ as a
linear multiple of the various vectors $\hat{v}_0$. Take a basis $e_i, i=1,2,\dots, m$ for $V_0$, and then
$\hat{e}_1 \wedge \hat{e}_2 \wedge \dots \wedge \hat{e}_m \wedge 
\left[\hat{e}_i, \hat{e}_j\right]=0$
for any $i$ and $j$, at every point of $M_0'$.

As we move along the flow of any one of these vector fields
$\hat{v}_0$, starting at any point of $M_0'$, we must find by analyticity that
$\hat{e}_1 \wedge \hat{e}_2 \wedge \dots \wedge \hat{e}_m \wedge 
\left[\hat{e}_i, \hat{e}_j\right]=0$
all along the flow line. We can assume that $M'_0$ is connected. The tangent spaces of $M'_0$ are the
span of the vector fields $\hat{v}_0$. Therefore the flows of these $\hat{v}_0$ vector fields act 
locally transitively on $M_0'$. But then since $M_0'$ is connected,
these vector fields act transitively on $M_0'$.

Let $\tilde{M}_0$ be the orbit of all of the vector fields
$\hat{v}_0$ through some point of $M'_0$.
Then $\tilde{M}_0$  is a submanifold of $M_0 \times M_1$ 
(perhaps not embedded; see Sussmann \cite{Sussmann:1973a,Sussmann:1973b}),
and $\tilde{M}_0$ contains $M'_0$.
At every point of $\tilde{M}_0$, we find 
$\hat{e}_1 \wedge \hat{e}_2 \wedge \dots \wedge \hat{e}_m \wedge 
\left[\hat{e}_i, \hat{e}_j\right]=0$ by analyticity.
% Restricting the Pfaffian system $\omega_1 = A \omega_0$ to $\tilde{M}_0$, 
% it therefore satisfies the conditions of the Frobenius theorem.
% Therefore $\tilde{M}_0$ is foliated by the
% integral manifolds of the Pfaffian system. The vector fields 
% span the tangent spaces of the leaves of the integral manifolds.
% The vector fields act transitively
% on $\tilde{M}_0$, and span the tangent spaces
% of $\tilde{M}_0$, so there is a single leaf
% of this foliation: all of $\tilde{M}_0$.
So $\tilde{M}_0$ is an integral manifold of the Pfaffian
system $\omega_1=A \omega_0$, invariant under the flows of the various $\hat{v}_0$.

Now consider the map $p : \tilde{M}_0 \to M_0$ given
by composing the inclusion $\tilde{M}_0 \subset M_0 \times M_1$
with the forgetful submersion $M_0 \times M_1 \to M_0$.
Clearly $p_* \hat{v}_0 = \vec{v}_0$.
Therefore $p$ is a local diffeomorphism. 
The coframing on $\tilde{M}_0$ is $p^* \omega_0$.

We need to see why $p$ is a covering map. Pick
a point $m_0 \in M_0$, say with logarithm
$\log_{m_0}$ a diffeomorphism in a neighborhood $U_0$ of $m_0$.
Then clearly logarithmic coordinates
around each point in $p^{-1}m_0$ are defined throughout
$p^{-1} U_0$. Let $Z=p^{-1} m_0$. Given any point $m_0' \in M_0$,
take a smooth path $x(t)$ with $x(0)=m_0$ and $x(1)=m_0'$
and $x'(t) \ne 0$. Then the time varying vector field $\vec{x}'(t)$
on $\tilde{M}_0$ will flow the fiber $p^{-1}m_0$ to $p^{-1}m_0'$.
Clearly $p$ is a covering map.

If we started off with a morphism defined on
an open set $U_0 \subset M_0$, then
we can map $U_0 \to \tilde{M}_0 \subset M_0 \times M_1$,
as the graph of that morphism. So this lifts
the inclusion $U_0 \to M_0$ to a local diffeomorphism
$U_0 \to \tilde{M}_0$. 
\end{proof}

\begin{corollary}
Under the hypotheses of theorem~\ref{theorem:GlobalMorphismOfFramings},
if $A$ is injective, then the map $\tilde{M}_0 \to M_1$ 
determined by the $A$-morphism is an immersion. If moreover $M_0$ is complete, then
every $A$-morphism $F' : M_0' \to M_1$ striking
a point in the image of $F$ must factor $F' = \tilde{F} \iota$
where $\iota : \tilde{M}_0' \to \tilde{M}_0$ is a local isomorphism from
a covering space $\tilde{M}_0'$ of $M_0'$.
\end{corollary}

\begin{remark}
 If $M_0$ is itself simply connected and complete,
then clearly $\tilde{M}_0=M_0$, i.e.each morphism extends uniquely to be defined on all of $M_0$.
\end{remark}

\begin{proof}
Completeness of the coframing on $M_0$ clearly implies completeness
of the pullback coframing on any covering space, so on $\tilde{M}_0$.
Pick manifold $M_0'$ with coframing $\omega_0' \in \nForms{1}{M_0'} \otimes V_0$. 
Pick any $A$-morphism $f' : M_0' \to M_1$. In logarithmic coordinates, 
$y=f'\left(x'\right)=Ax$ and $y=f(x)=Ax$, so that $x' \mapsto x=x'$
is a local diffeomorphism. Since $f$ and $f$ are morphisms,
they must identify constant vector fields on $M_0'$ and on $M_0$
respectively with those on $M_1$. Therefore $x' \mapsto x=x'$
is a local isomorphism. By theorem~\vref{theorem:GlobalMorphismOfFramings},
this local isomorphism extends to a local isomorphism
of a covering space.
\end{proof}

\begin{lemma}
Suppose that $A : V_0 \to V_1$ is a linear surjection. Suppose that there is
an $A$-morphism of coframings $M_0 \to M_1$. If $M_0$ is complete,
then $M_1$ is complete. If every canonical metric of $M_0$ is complete,
then every canonical metric of $M_1$ is complete.
\end{lemma}
\begin{proof}
 Take a splitting $V_0 = V_1 \oplus V_1^{\perp}$, and then
develop lines/curves from $V_1$ to $M_1$ by first including these
lines/curves into $V_0$ and then developing into $M_0$ and
then composing with the morphism.
\end{proof}

\begin{lemma}
Suppose that $A : V_0 \to V_1$ is any linear map,
with image $W_1 \subset V_1$. Suppose that
$q : V_1 \to V_1/W_1$ is the obvious linear projection. Take
a manifold $M_1$ with coframing $\omega_1$.
The integral manifolds of the Pfaffian system
$q \omega_1=0$ on $M_1$ on which $\omega_1$ has rank
equal to $\dim W_1$ are locally the images of $A$-morphisms
to $M_1$. The image of any $A$-morphism $M_0 \to M_1$, with $M_0$ connected, 
is an open subset of precisely one of these integral manifolds.
\end{lemma}
\begin{proof}
Without loss of generality, we can assume that as matrices,
\[
A=
\begin{pmatrix} 
I & 0 \\
0 & 0
\end{pmatrix},
\]
and
\[
q=
\begin{pmatrix}
0 & I
\end{pmatrix},
\]
and write out the Pfaffian equations and the result is clear.
\end{proof}

\begin{lemma}
Suppose that $A$ is an injective linear map. 
Any $A$-morphism of coframings is completely
determined by its value at a single point.
\end{lemma}
\begin{proof}
 The morphism intertwines the flows of the constant vector fields, and this determines
the graph of the morphism locally.
\end{proof}

\begin{definition}
Suppose that $\omega \in \nForms{1}{M} \otimes V$ is a coframing on a manifold $M$.
Denote the torsion by $T : M \to V \otimes \Lm{2}{V}^*$, so that $d \omega = T \omega \wedge \omega$.
\end{definition}

\begin{remark}
Suppose that $A : V_0 \to V_1$ is a linear injection, with image $W_1$, and that
$q : V_1 \to V_1/W_1$ is the obvious linear projection.
Denote by $\Lambda^2 A : \Lm{2}{V_0} \to \Lm{2}{V_1}$
the obvious map $\Lambda^2 A(u \wedge v)=A(u) \wedge A(v)$.
Suppose that $\omega_1 \in \nForms{1}{M_1} \otimes V_1$
is a coframing on a manifold $M_1$.
The torsion of the Pfaffian system for
$A$-morphisms (i.e. the system $q\omega_1 = 0$)
is $q T_1 \Lambda^2 A = 0$. Therefore no morphisms
map to points where $q T_1 \Lambda^2 A \ne 0$.

Let $Z$ be the set of points of $M_1$
at which $q T_1 \Lambda^2 A=0$. 
If the framing is real analytic, then
$Z$ is a real analytic variety, and stratified
by smooth subvarieties. Each stratum
is foliated by integral manifolds. If
$\omega_1$ has full rank, i.e. rank equal to $\dim V_0$, on the tangent spaces
of a leaf on one of these strata, then
the leaf has $\omega_1$ as coframing. The
inclusion of the leaf into $M_1$ is
a morphism.
\end{remark}

\begin{remark}
Let's consider a noninjective linear 
map $A : V_0 \to V_1$. 
Suppose that $M_0$ and $M_1$ are
manifolds bearing coframings 
$\omega_0 \in \nForms{1}{M_0} \otimes V_0$
and
$\omega_1 \in \nForms{1}{M_1} \otimes V_1$.
Pick a particular integral manifold of the Pfaffian
system $q \omega_1=0$. Replace $M_1$ by this integral
manifold. So then we can effectively
assume that $A$ is surjective.

If there is an $A$-morphism, then the fibers of $M_0 \to M_1$
will be integral manifolds of the Pfaffian system $A \omega_0=0$.
If we try to construct an $A$-morphism $M_0 \to M_1$, we first have to ask whether
the Pfaffian system $A \omega_0=0$ satisfies the conditions of the
Frobenius theorem. If we 
write $d \omega_0 = T_0 \omega_0 \wedge \omega_0$, and we write an exact sequence
\[
\xymatrix{
0 \ar[r] &  \ker A \ar[r]^{\iota} & V_0 \ar[r]^{A} & V_1 \ar[r]^{q} & V_1/W_1 \ar[r] & 0,
}
\]
then in order that there could be a morphism $M_0 \to M_1$, 
the Frobenius theorem requires $T_0 \Lambda^2 \iota=0$ at all points
of $M_0$. If this holds, then $M_0$ is foliated by the leaves of this Pfaffian
system. If quotient space $M_{\bar{0}}$ of the foliation
is Hausdorff, then it is a manifold with a coframing $\omega_{\bar{0}}$,
so that under the projection map $p : M_0 \to M_{\bar{0}}$,
$p^* \omega_{\bar{0}} = q \omega_0$. We then
try to construct a $\bar{A}$-morphism $M_{\bar{0}} \to M_1$.
Note that $\bar{A}$ is injective.
\end{remark}

\begin{remark}
Suppose that $A$ is an injective linear map.
Lets consider the local problem of constructing manifolds
$M_1$ which admit $A$-morphisms from a given manifold $M_0$
with a given coframing. Without loss of generality 
we can assume that $A=(I \ 0)$. Take any product $M_0 = M_1 \times U$
with any manifold $U$ with trivial tangent bundle. Take $\omega_U$
any coframing on $U$. Let $\omega_0 = \omega_1 + f \, \omega_U$
where $f : M_0 \to V_U^* \otimes \left( V_1 \oplus V_U \right)$
is any linear map. Of course, any $A$-morphism to $M_1$ will be locally isomorphic to this
one for a suitable choice of $f$.
\end{remark}

\begin{corollary}
Suppose that $A : V_0 \to V_1$ is an 
injective linear map with image $W_1$. Suppose that
$q : V_1 \to V_1/W_1$ is
the obvious linear projection. Take
a manifold $M_1$ with coframing $\omega_1 \in \nForms{1}{M_1} \otimes V_1$.
Suppose that $d \omega_1 = T_1 \omega_1 \wedge \omega_1$.
Suppose that $q T_1 \Lambda^2 A=0$.

Then through each point $m_1 \in M_1$, there is a unique
maximal immersed submanifold $M_0=M_0\left(m_1\right)$ 
containing $m_1$ on which $q \omega_N=0$, with coframing $\omega_0=A^{-1} \omega_1$.
The inclusion $M_0 \to M_1$ is an $A$-morphism. 

Every $A$-morphism $F : M'_0 \to M_1$, with $M'_0$ connected and
$m_1$ in the image of $F$
factors as $M'_0 \to M_0 \to M_1$, where
$M_0' \to M_0$ is a local isomorphism of coframings and $M_0=M_0\left(m_1\right)$
for any point $m_1$ in the image of $F$.
\end{corollary}
\begin{proof}
Clearly we must take $M_0$ to be the leaf of the Pfaffian system 
through a point of $M_1$. We can take $M_0$ to be a maximal leaf.
\end{proof}

\section{The Frobenius--Gromov theorem for coframing morphisms}

\begin{definition}
If $\omega \in \nForms{1}{M} \otimes V$ is a coframing
on a manifold $M$, then the \emph{1-torsion}
of $\omega$ is the function $T : M \to V \otimes \Lm{2}{V}^*$
defined by
\[
d \omega = \frac{1}{2} T \omega \wedge \omega.
\]
We also write $T$ as $T^{(0)}$.
For any integer $p>1$, the $p$-torsion
is defined recursively
$T^{(p)}$ to be the function
\[
 T^{(p)} : M \to V \otimes \Lm{2}{V}^* \otimes V^{* \otimes p}
\]
given by
\[
dT^{(p-1)} = T^{(p)} \omega.
\]
\end{definition}

Clearly $T= T^{(1)}, T^{(2)}, \dots, T^{(p)}$ at a point $m \in M$ 
determine the $p$-jet of $T$ at $m$.

\begin{lemma}\label{lemma:BracketsOfAllOrders}
Suppose that $\omega \in \nForms{1}{M} \otimes V$ is a coframing.
Pick any $v_1, v_2, \dots, v_p \in V$. Write out the brackets as
\[
\ad\left(\vec{v}_1\right)
\ad\left(\vec{v}_2\right)
\dots
\ad\left(\vec{v}_{p-1}\right)
\vec{v}_p
=
\overrightarrow{\tau_p\left(v_1,v_2,\dots,v_p\right)},
\]
for some function $\tau_p : M \to \otimes^p V^* \otimes V$.
Then $\tau_1(v)=v$ and inductively
\[
\tau_{p+1}
\left(
v_1, v_2, \dots, v_{p+1}
\right)
=
\LieDer_{\vec{v}_1}  \tau_p\left(v_2, v_3, \dots, v_p\right)
-
T
\left(
v_1 \wedge
\tau_p\left(v_2, v_3, \dots, v_p\right)
\right).
\]
In particular, we can compute the values of $\tau_1, \tau_2, \dots, \tau_p$ at a point
$m \in M$ from the values of $T^{(1)}, T^{(2)}, \dots, T^{(p-1)}$ at $m$.
\end{lemma}
\begin{proof}
Clearly $\tau_1(m)(v)=v$.
Moreover, from the definition of $T$,
clearly
\[
\tau_2(m)\left(v_1,v_2\right)=T(m)\left(v_1,v_2\right).
\]

Suppose by induction that $\tau_p(m)$ is some
function of $T(m), T^{(1)}(m), \dots, T^{(p)}(m)$.
Pick some vectors $v_1, v_2, \dots, v_{p+1} \in V$.
Let $X$ be the vector field
\[
\ad\left(\vec{v}_2\right)
\ad\left(\vec{v}_3\right)
\dots
\ad\left(\vec{v}_{p}\right)
\vec{v}_{p+1}.
\]
Then
\begin{align*}
\tau_{p+1}
\left(
v_1, v_2, \dots, v_{p+1}
\right)
&=
\omega\left(
\left[\vec{v}_1, X\right]
\right)
\\
&=
\LieDer_{\vec{v}_1} \left(X \hook \omega\right)
-
\LieDer_{X} \left(\vec{v}_1 \hook \omega\right)
-
d \omega
\left(
\vec{v}_1,
X
\right)
\\
&=
\LieDer_{\vec{v}_1}  \tau_p\left(v_2, v_3, \dots, v_p\right)
-
T
\left(
v_1 \wedge
\tau_p\left(v_2, v_3, \dots, v_p\right)
\right).
\end{align*}
\end{proof}

\begin{corollary}
 Given the value of a coframing $\omega$ at a point $m$, and the value of the $j$-torsion of the coframing at $m$,
for $j=1,\dots,p-1$, we can compute the bracket at $m$: $\ad\left(\vec{v}_1\right)
\ad\left(\vec{v}_2\right)
\dots
\ad\left(\vec{v}_{p-1}\right)
\vec{v}_p$ of any
collection of up to $p+1$ constant vector fields 
\end{corollary}

\begin{definition}
Suppose that $M_0$ and $M_1$ are manifolds with coframings
$\omega_0 \in \nForms{1}{M_0} \otimes V_0$ and 
$\omega_1 \in \nForms{1}{M_1} \otimes V_1$
and that $A : V_0 \to V_1$ is a linear map. 
For any nonnegative integer $j$, the \emph{$j$-th order obstruction} $\Ob^{(j)}$
of the pair $\left(\omega_0,\omega_1\right)$ is the function
\[
\Ob^{(j)}
=
T^{(j)}_1\left(m_1\right) = A T^{(j)}_0\left(m_0\right)
\Lambda^{2} A \otimes A^{\otimes (p-1)},
\]
on $M_0 \times M_1$, $j=0, 1, 2, \dots$.
The \emph{$j$-th order morphism locus} of the pair $\left(\omega_0,\omega_1\right)$
is the set $V_j \subset M_0 \times M_1$ where the obstructions
of order not more than $j$ vanish.
The \emph{morphism locus} is the intersection of all
of the $j$-order morphism loci, $j=0, 1, 2, \dots$.
\end{definition}

\begin{definition}
Suppose that $M_0$ and $M_1$ are manifolds with coframings
$\omega_0 \in \nForms{1}{M_0} \otimes V_0$ and 
$\omega_1 \in \nForms{1}{M_1} \otimes V_1$
and that $A : V_0 \to V_1$ is a linear map. 
We will say that a point $m_0 \in M_0$
\emph{hits} a point $m_1 \in M_1$ to order $p$ if
$\left(m_0,m_1\right)$ lies in the $p$-th order morphism locus.
We will say that a point $m_0 \in M_0$ \emph{hits}
a point $m_1 \in M_1$ if there is a morphism $F : U_0 \to M_1$
of an open neighborhood $U_0$ of $m_0$ so that $F\left(m_0\right)=m_1$.
\end{definition}

\begin{theorem}\label{theorem:FGMcoframings}
Suppose that $M_0$ and $M_1$ 
are smooth manifolds with smooth coframings
$\omega_0 \in \nForms{1}{M_0} \otimes V_0$ and 
$\omega_1 \in \nForms{1}{M_1} \otimes V_1$
and that $A : V_0 \to V_1$ is a linear map. 
A point $m_0$ hits a point $m_1$ just when
$\left(m_0,m_1\right)$ belongs to the morphism locus.
\end{theorem}
\begin{proof}
Clearly if a point $\left(m_0,m_1\right)$ does not belong
to the morphism locus then $m_0$ does not hit $m_1$.
So we only need to prove that if $\left(m_0,m_1\right)$
belongs to the morphism locus, then $m_0$ hits $m_1$.
Let $Z$ be the morphism locus.
For each vector $v \in V_0$, let $\bar{v}=\vec{v}+\overrightarrow{Av}$,
a vector field on $M_0 \times M_1$. These vector fields
precisely span the plane field given by $\omega_1=A \omega_0$.
We leave the reader to calculate that
\[
\bar{v} \hook d 
\left(
  \Ob^{(j)}
  \left(
    v_0, v_1, w_1, w_2, \dots, w_j
  \right)
\right)
=
  \Ob^{(j+1)}
  \left(
    v_0, v_1, w_1, w_2, \dots, w_j, v
  \right).
\]
It follows (by the Picard uniqueness
theorem) that $Z$ is invariant under $\bar{v}$
for any vector $v \in V_0$. So the orbit
of the family of all vector fields 
$\bar{v}$ through each point of $Z$
lies entirely inside $Z$.

Compute
\begin{align*}
\LieDer_{\bar{v}} \left( \omega_1 - A \omega_0 \right)
&=
\bar{v} 
\hook
d \left( \omega_1 - A \omega_0 \right)
\\
&=
T_1\left(Av,\omega_1\right)
-
A \, T_0\left(v,\omega_0\right).
\end{align*}
In particular, at any point of $Z$, for any vector $w \in V_0$,
\[
\bar{w} \hook \LieDer_{\bar{v}} \left( \omega_1 - A \omega_0 \right) = 0.
\]
Therefore at each point of $Z$, 
the 1-form $\LieDer_{\bar{v}} \left( \omega_1 - A \omega_0 \right)$
vanishes on $\omega_1=A \, \omega_0$, and is therefore a multiple
of $\omega_1-A \, \omega_0$. Integrating, we find that
the ideal of 1-forms generated by the components of $\omega_1-A \omega_0$
(in any basis for $V_1$) is invariant under the flow of
$\bar{v}$ through any point of $Z$. Therefore every 
vector field $\bar{w}$, for $w \in V_0$, is carried
by the flow of $\bar{v}$ into a linear combination
of vector fields of the form $\bar{u}$ for various $u \in V_0$.
Differentiating, the bracket $\left[\bar{v},\bar{w}\right]$
must also lie tangent to $\omega_1=A \, \omega_0$ through
each point of $Z$. 

Indeed we can see this result more directly
from lemma~\vref{lemma:BracketsOfAllOrders}:
compute out
\[
\left[\vec{v},\vec{w}\right]\left(m_0\right)
=
-\overrightarrow{T_0\left(m_0\right)(v,w)}\left(m_0\right),
\]
and
\[
\left[\overrightarrow{Av},\overrightarrow{Aw}\right]\left(m_1\right)
=
-\overrightarrow{T_1\left(m_1\right)(Av,Aw)}\left(m_0\right),
\]
so that at points of $Z$,
\[
\left[\bar{v},\bar{w}\right]\left(m_0,m_1\right)
=
-\overline{T_0\left(m_0\right)(v,w)}\left(m_0,m_1\right).
\]

The orbit through any point of $M_0 \times M_1$ of the family
of all vector fields $\bar{v}$ for $v \in V_0$ is a smooth
submanifold. On the orbit through any point of $Z$, the Frobenius theorem applies
to the Pfaffian system $\omega_1 = A \, \omega_0$.
The orbit is therefore foliated by maximal integral manifolds,
whose tangent planes are precisely spanned by the
various $\bar{v}$ vector fields. But then the
maximal integral manifolds are themselves orbits,
and therefore the entire orbit must be a single maximal integral manifold.
Its tangent space is precisely the span
of the various $\bar{v}$ vector fields.
In particular, it projects locally diffeomorphically
to $M_0$.
\end{proof}

\begin{example}
Suppose that $M_0$ and $M_1$ are manifolds with coframings
\[
\omega_0 \in \nForms{1}{M_0} \otimes V_0
\] 
and 
\[
\omega_1 \in \nForms{1}{M_1} \otimes V_1 
\]
and that $A : V_0 \to V_1$ is a linear map. 
The $j$-order morphism locus is a real analytic
set, and therefore so is the morphism locus.
By the descending chain condition, we can
cover $M_0 \times M_1$ by open sets so
that on each of these open sets, the morphism locus
is equal to the $j$-th order morphism
locus for some $j$. So we can test for local morphisms
with purely local calculations at some order.
\end{example}

\begin{corollary}
Two real analytic coframings agree at a point and have the same torsions of all
orders at that point just when they are equal in a neighborhood of that point.
\end{corollary}
\begin{proof}
Let $A=I$.
% ??? Careful: here we use BCH. I will have to prove it properly
% or give up on it.
% 
% Their iterated brackets of any finite set of constant vector fields  
% must agree. By the Baker--Campbell--Hausdorff formula for
% coframings (lemma~\vref{lemma:BCH}), the two coframings must have
% the same values of $\log_m e^{t \vec{a}} e^{t \vec{b}}m$
% for sufficiently small $a, b$. Replacing $b$ by $sb$ 
% for a small number $s$, and taking $\frac{d}{ds}$, we see
% that the two coframings have the same expression for
% $\vec{b}\left(e^{t\vec{a}}m\right)$ in logarithmic coordinates.
% Since $e^{t\vec{a}}m$ can be any point near $m$, for
% suitable choice of $a$, this ensures that constant vector
% fields are identical throughout a region in which logarithmic coordinates
% are defined. Therefore the coframings are identical in these
% regions, and therefore identical everywhere.
\end{proof}

\section{Hartogs extension of coframing morphisms}

\begin{definition}
 A \emph{Stein manifold} is a complex manifold $X$ which is biholomorphic to a complex 
submanifold $Y \subset \C{N}$, for some integer $N \ge 0$, so that $Y$ is closed as a subset.
\end{definition}

\begin{definition}
An \emph{envelope of holomorphy} of a complex manifold $M$ is a Stein manifold $\hat{M}$ so that
\begin{enumerate}
\item $M \subset \hat{M}$ is an open subset,
\item every holomorphic function on $M$ extends to $\hat{M}$.
\end{enumerate}
\end{definition}

\begin{theorem}[Docquier and Grauert \cite{Docquier/Grauert:1960}]
If $D$ is a domain in a Stein manifold, then $D$ has a unique envelope of holomorphy $\hat{D}$, up
to a unique biholomorphism which is the identity on $D$.
\end{theorem}

\begin{lemma}\label{lemma:ExtendCoframing}
 Suppose that $D$ is a domain in a Stein manifold, with envelope of holomorphy
$\hat{D}$. Then every holomorphic coframing of $D$ extends uniquely to a holomorphic coframing of $\hat{D}$.
\end{lemma}
\begin{proof}
Take a coframing $\omega \in \nForms{1}{D} \otimes V$ on $D$.
Every holomorphic differential form extends uniquely from $D$ to $\hat{D}$
(see McKay \cite{McKay:2009} p. 10 proposition 2), so extend $\omega$ to $\hat{D}$.
If we take a basis for $V$, we can write 
\[
 \omega=
\begin{pmatrix}
 \omega^1 \\
 \omega^2 \\
 \vdots \\
 \omega^n
\end{pmatrix}.
\]
The set of points of $\hat{D}$ at which $\omega$ fails to be a coframing is the complex hypersurface
\[
 \omega^1 \wedge \omega^2 \wedge \dots \wedge \omega^n=0.
\]
Every nonempty complex hypersurface in $\hat{D}$ intersects $D$ (see McKay \cite{McKay:2009} p. 8 lemma 11).
Therefore this complex hypersurface is empty.
\end{proof}

\begin{theorem}\label{theorem:CoframingExtension}
Suppose that $M$ is a complex manifold with complete holomorphic coframing.
Suppose that $D$ is a domain in a Stein manifold, and $D$ has a holomorphic coframing, and that
$f : D \to M$ is a holomorphic morphism of coframings. Then $f$ extends uniquely to a 
holomorphic morphism of coframings $f : \hat{D} \to M$ on the envelope of 
holomorphy $\hat{D}$ of $D$.
\end{theorem}
\begin{proof}
The coframing on $D$ extends uniquely to a holomorphic coframing on $\hat{D}$
by lemma~\vref{lemma:ExtendCoframing}.
The morphism $f$ factors uniquely through a holomorphic morphism on a covering space
of $\hat{D}$, by theorem~\vref{theorem:GlobalMorphismOfFramings},
say $f=f' \iota$ with $f' : \hat{D}' \to M$ a morphism from a covering
space, $p : \hat{D}' \to \hat{D}$, and $\iota' : D \to \hat{D}'$ 
a lift of the inclusion $\iota : D \to \hat{D}$, so $p \iota' = \iota$. On fundamental groups,
$p_* \iota' = \iota$. Note that 
$\iota_* : \pi_1\left(D\right) \to \pi_1\left(\hat{D}\right)$. 
It is well known (see McKay \cite{McKay:2009} p. 7 lemma 6) 
that $\iota_* : $ is surjective. Therefore $p_* \iota'_*$ is surjective,
and so $p_*$ is surjective. Because $p$ is a covering map, $p_*$ is injective. 
Therefore $p_*$ is an isomorphism, i.e. $\hat{D}=\hat{D}'$. 
\end{proof}

\section{Morphisms of Cartan geometries}\label{section:MorphismsOfCartanGeometries}

We now move on from studying coframings to studying Cartan geometries.
We will allow Sharpe's generalization of the concept of Cartan
geometry, but we will use our own terminology.

\begin{definition}[Sharpe \cite{Sharpe:1997} p. 174]
A \emph{local model} $\left(H,\mathfrak{g}\right)$ is a Lie group $H$ and an $H$-module $\mathfrak{g}$,
with the structure of a Lie algebra, on which $H$ acts as automorphisms,
with an injection $\mathfrak{h} \to \mathfrak{g}$ of Lie algebras, where $\mathfrak{h}$
is the Lie algebra of $H$.
\end{definition}

\begin{remark}
 Sharpe \cite{Sharpe:1997} uses the term ``model'' for what we
call the ``local model''. We need the distinction because some
of our results require a homogeneous space as model. Keep
in mind that a local model might not arise from any
homogeneous space, but merely from a Lie subgroup $H \subset G$
which need not be closed. Sharpe's terminology is confusing
and its use should be discouraged.
\end{remark}

\begin{definition}
 The \emph{local model} of a homogeneous space $G/H$ is $\left(H,\mathfrak{g}\right)$
where $\mathfrak{g}$ is the Lie algebra of $G$.
\end{definition}

\begin{remark}
In all of our examples, the only local models we will consider will be 
the local models of homogeneous spaces.
\end{remark}

\begin{definition}\label{def:CartanConnectionSharpe}
Let $\left(H,\mathfrak{g}\right)$ be a local model. A $\left(H,\mathfrak{g}\right)$-geometry
(also known as a \emph{Cartan geometry} modelled on $\left(H,\mathfrak{g}\right)$) on a manifold
$M$ is a principal right $H$-bundle $E \to M$ and a 1-form $\omega
\in \nForms{1}{E} \otimes \mathfrak{g}$ called the \emph{Cartan
connection}, which satisfies the following conditions:
\begin{enumerate}
\item
Denote the right action of $h \in H$ on $e \in E$ by $r_h e$. The
Cartan connection transforms in the adjoint representation:
\[
r_h^* \omega = \Ad_h^{-1} \omega.
\]
\item
$\omega_e : T_e E \to \mathfrak{g}$ is a linear isomorphism at each
point $e \in E$.
\item
For each $A \in \mathfrak{g}$, define a vector field $\vec{A}$ on
$E$ 
(called a \emph{constant vector field})
by the equation $\vec{A} \hook \omega = A$. For $A \in
\mathfrak{h}$, the vector fields $\vec{A}$ generate the right
$H$-action:
\[
\vec{A}(e) = \left. \frac{d}{dt} r_{e^{tA}} e \right|_{t=0}, \text{ for all $e \in E$}.
\]
\end{enumerate}
\end{definition}

\begin{definition}
A \emph{morphism of local models} $\Phi : \left(H_0,\mathfrak{g}_0\right) \to \left(H_1,\mathfrak{g}_1\right)$
is a pair of maps: (1) a Lie group morphism which we denote $\Phi : H_0 \to H_1$
(and we also denote by $\Phi$ the induced morphism of Lie algebras 
$\Phi : \mathfrak{h}_0 \to \mathfrak{h}_1$)
and (2) a linear map, which we denote $\Phi : \mathfrak{g}_0 \to \mathfrak{g}_1$,
extending the Lie algebra morphism $\Phi : \mathfrak{h}_0 \to \mathfrak{h}_1$, so that
\[
\Phi\left(h_0 A_0\right)=\Phi\left(h_0\right)\Phi\left(A_0\right),
\]
for $h_0 \in H_0$ and $A_0 \in \mathfrak{g}_0$.
\end{definition}

\begin{definition}
Suppose that $G_0/H_0$ and $G_1/H_1$ are homogeneous spaces.
A \emph{morphism} $\Phi : G_0/H_0 \to G_1/H_1$ of homogeneous spaces
is a Lie group morphism $\Phi : G_0 \to G_1$ which takes $H_0$ to $H_1$. 
\end{definition}

\begin{example}\label{example:MorphismOfHomSpaces}
A morphism of homogeneous spaces $\Phi : G_0/H_0 \to G_1/H_1$
induces the obvious local model morphism
$\Phi : \left(H_0,\mathfrak{g}_0\right) \to \left(H_1,\mathfrak{g}_1\right)$.
\end{example}

\begin{definition}
Suppose that $\Phi : \left(H_0,\mathfrak{g}_0\right) \to \left(H_1,\mathfrak{g}_1\right)$
is a local model morphism. Suppose that $E_0 \to M_0$ is a $\left(H_0,\mathfrak{g}_0\right)$-geometry with Cartan
connection $\omega_0$ and that $E_1 \to M_1$ is a $\left(H_1,\mathfrak{g}_1\right)$-geometry
with Cartan connection $\omega_1$. A \emph{morphism of Cartan geometries} modelled on $\Phi$
(also called a $\Phi$-morphism) is an $H_0$-equivariant map $F : E_0 \to E_1$
so that $F^* \omega_1 = \Phi \omega_0$.
\end{definition}

\begin{remark}
Clearly under a morphism $F$, the flow lines of each vector field $\vec{A}$ on $E_0$ are carried
to the flow lines of $\overrightarrow{\Phi(A)}$. Indeed (in the notation of the previous definition)
\[
F'\left(e_0\right) \vec{A}\left(e_0\right)=\overrightarrow{\Phi(A)}\left(F\left(e_0\right)\right)
\]
for every point $e_0 \in E_0$.
\end{remark}

\begin{definition}
 If a morphism of Cartan geometries is modelled on a model morphism,
and the model morphism is induced from a morphism of homogeneous spaces (as in 
example~\vref{example:MorphismOfHomSpaces}), we will say that the morphism
of Cartan geometries is modelled on the morphism of homogeneous spaces.
\end{definition}

\begin{remark}
I do not know of a serious example of a Cartan geometry
which is modelled on a local model without being modelled
on a homogeneous space. Nevertheless, there are serious
examples, as we will see, of morphisms of Cartan geometries
which are modelled on morphisms of local models, 
even though the morphisms of local models do not
arise from morphisms of homogeneous spaces. It appears
that this phenomenon was behind Sharpe's use of 
what we have called local models, even though Sharpe
did not have morphisms specifically described in his work.
\end{remark}

\begin{proposition}
Suppose that $\Phi : \left(H_0,\mathfrak{g}_0\right) \to \left(H_1,\mathfrak{g}_1\right)$ 
is a local model morphism.
Suppose that $\pi_0 : E_0 \to M_0$ is an $\left(H_0,\mathfrak{g}_0\right)$-geometry.
Suppose that $\pi_1 : E_1 \to M_1$ is an $\left(H_1,\mathfrak{g}_1\right)$-geometry.
Suppose that $M_0$ is connected.
Pick a point $e_0 \in E_0$ and a point $e_1 \in E_1$.
Let $m_0 = \pi_0 \left(e_0\right) \in M_0$ and 
$m_1 = \pi_1 \left(e_1\right) \in M_1$. 
There is a morphism $F : E_0 \to E_1$ modelled on $\Phi$
with $F\left(e_0\right)=e_1$ if and only if every loop in $M_0$ starting and ending at 
$m_0$ develops (with frames $e_0$ and $e_1$) to a loop in $M_1$ starting and ending at $m_1$.
\end{proposition}
\begin{proof}
Suppose that we have a morphism $F$, so that $F\left(e_0\right)=e_1$.
Since $F$ is $H_0$-equivariant, it induces a unique smooth map
$f : M_0 \to M_1$ by $f \circ \pi_0 = \pi_1 \circ F$.
For any map $h : X \to M_0$ from any manifold $X$ (with boundary,
corners, etc.) the map $H : h^* E_0 \to \left(f \circ h\right)^* E_1$
given by $H\left(m,e\right)=\left(f(m),F(e)\right)$
is clearly a $\Phi$-development.

On the other hand, suppose that 
every loop in $M_0$ starting and ending at 
$m_0$ develops to a loop in $M_1$ starting and ending at $m_1$,
with frames $e_0$ and $e_1$. 
We want to apply our previous results on development of coframings.
Pick a smooth loop
$\Gamma_0 : [0,1] \to E_0$. Then let $\gamma_0 = \pi_0 \circ \Gamma_0$
so $\gamma_0 : [0,1] \to M_0$ is a smooth loop. Develop
$\gamma_0$ to a smooth loop $\gamma_1 : [0,1] \to M_1$,
with some development $H : \gamma_0^* E_0 \to \gamma_1^* E_1$
so that $H^* \omega_1 = \Phi \omega_0$. 
Now let $\Gamma_1 = H \circ \Gamma_0$. Since $\Gamma_0$
is a smooth loop, $\Gamma_1$ is also a smooth loop.
Therefore we can $\Phi$-develop loops to loops for the
coframings $\omega_0$ and $\omega_1$, so by
proposition~\vref{proposition:DevelopmentToMorphism},
there is a unique $\Phi$-coframing morphism
$F : E_0 \to E_1$ so that $\Phi$-coframing development is 
composition with $F$. By uniqueness, $F$ is $H_0$-equivariant,
and therefore is a Cartan geometry morphism.
\end{proof}

\begin{example}
 If $E \to M$ is a $G/H$-geometry, then the bundle map $E \to M$ is
itself a morphism, modelled on the morphism $G \to G/H$.
\end{example}

\begin{example}
Suppose that $E \to M$ is an $\left(H_1,\mathfrak{g}\right)$-geometry, 
and that $H_0 \subset H_1$ is a closed
subgroup. Then $E \to E/H_0$ is a $\left(H_0,\mathfrak{g}\right)$-geometry (with the same
Cartan connection) called a \emph{lift} of $E \to M$.
The obvious map $E/H_0 \to M$ is a morphism, modelled on 
the obvious inclusion $\left(H_0,\mathfrak{g}\right) \to \left(H,\mathfrak{g}\right)$.
It is intuitively clear what the lift of a morphism should mean.
\end{example}

The lift is ``canonical'' or ``universal'' in the following sense.

\begin{lemma}\label{lemma:liftIsCanonical}
Suppose that $H_0$ is a closed subgroup of $H_1$. 
Suppose that $\left(H_1,\mathfrak{g}\right)$ is a local model.
Take the obvious inclusion map of local models 
$\Phi : \left(H_0,\mathfrak{g}\right) \to \left(H_1,\mathfrak{g}\right)$.
Suppose that $E_0 \to M_0$ and $E_1 \to M_1$ are Cartan geometries modelled on $\left(H_0,\mathfrak{g}\right)$ 
and $\left(H_1,\mathfrak{g}\right)$.
Every morphism of Cartan geometries $F : E_0 \to E_1$ modelled
on $\Phi$ factors uniquely as $F=F_{\text{lift}} F_0$ where $F_{\text{lift}} : E_1/H_0 \to E_1/H_1$
is the lift, and $F_0$ is a uniquely determined local isomorphism.
\end{lemma}
\begin{proof}
Just write $F$ as $F=\operatorname{id} F=F_{\text{lift}} F_0$.
\end{proof}

\begin{example}
A \emph{mutation} is a local model morphism $\Phi : \left(H_0,\mathfrak{g}_0\right) \to \left(H_1,\mathfrak{g}_1\right)$
with $\Phi : H_0 \to H_1$ an isomorphism of Lie groups and $\Phi : \mathfrak{g}_0 \to \mathfrak{g}_1$ an
$H_0$-module isomorphism; see
Sharpe \cite{Sharpe:2002} p. 154 for examples. A \emph{mutation} of Cartan geometries $E_0 \to M_0$ and $E_1 \to M_1$
is a morphism modelled on a mutation and for which the underlying map $M_0 \to M_1$ is a diffeomorphism.
\end{example}

\begin{definition}
Say that a model morphism $\Phi : \left(H_0,\mathfrak{g}_0\right) \to \left(H_1,\mathfrak{g}_1\right)$
is a \emph{base epimorphism} if $\mathfrak{g}_0/\mathfrak{h}_0 \to \mathfrak{g}_1/\mathfrak{h}_1$ is onto, etc.,
and a \emph{fiber epimorphism} if $H_0 \to H_1$ is onto, etc. and
an \emph{epimorphism} if both a base and fiber epimorphism, etc.
\end{definition}

We will need some topological lemmas before we can study the
global structure of Cartan geometry epimorphisms.

\begin{lemma}
 Suppose that $\pi : E \to M$ is a fiber bundle mapping
and $X$ is a complete vector field on $M$. Then there
is a complete vector field $Y$ on $E$ so that 
$Y$ lifts $X$ ,i.e. 
$\pi'(e)Y(e)=X(\pi(e))$ for each $e \in E$. Moreover,
the set of complete vector fields on $E$ which lift
complete vector fields on $M$ acts locally transitively
on $E$. 
\end{lemma}
\begin{proof}
Take a locally finite cover of $M$ by relatively
compact open sets $U_{\alpha}$ and a subordinate partition
of unity $f_{\alpha}$. We can pick the open sets
$U_{\alpha}$ so that above each one, $\pi : E \to M$
is a trivial bundle, say $F_{\alpha} : \pi^{-1}U_{\alpha} \to U_{\alpha} \times Q$,
where $Q$ is a typical fiber. Any lift $Y_{\alpha}$ of $f_{\alpha} X$
will have the form
\[
Y_{\alpha}\left(m,q\right)=\left(f_{\alpha}\left(m\right) X\left(m\right),Z_{\alpha}\left(m,q\right)\right),
\]
for some $Z\left(m,q\right) \in T_q Q$. We can pick
$Z_{\alpha}$ to have compact support in $U_{\alpha} \times Q$. Let $Y=\sum_{\alpha} Y_{\alpha}$,
which is well defined because the partition of unity
is locally finite. 

Suppose that the flow of $Y$ through some point $e_0$ exists
only for a time $T$. Let $m_0 = \pi\left(e_0\right)$.
The flow of $X$ exists for all time,
and $e^{TX}m_0$ belongs to a finite number of open sets
$U_{\alpha}$, let's say to $U_1, U_2, \dots, U_N$.
So for time $t$ near $T$, $e^{tX}m_0$
stays inside $U_1 \cap U_2 \cap \dots \cap U_N$.
Expressed in terms of the trivialization
over $U_1$, we will find that inside 
$\pi^{-1} \left(U_1 \cap U_2 \cap \dots \cap U_N\right)$,
the vector field $Y$ has the form
\[
Y\left(m,q\right)
=
\left(X(m),
\sum_{j=1}^N Z_j\left(m,q\right)
\right),
\]
for a finite set of compactly supported functions
$Z_j$. Therefore the flow of $m$ is $\dot{m}=X(m)$,
extending for time past $T$, and the flow
of $q$ is $\dot{q}=\sum_{j=1}^{N} Z_j\left(e^{tX}m_0,q\right)$,
which has compact support in $q$, so is complete.
Clearly we can alter the choice of $Z_1$ as we
like to ensure that $Y$ has required value at
a single point, so that the tangent space of
$E$ at each point is spanned by complete lifts. 
\end{proof}

\begin{corollary}
 The composition of fiber bundle maps is a fiber bundle map.
\end{corollary}
\begin{proof}
Suppose that $P \to Q$ and $Q \to R$ are fiber bundle maps.
 Take the set $V_R$ of complete vector fields on $R$.
 Take the set $V_Q$ of complete lifts of those vector fields on $Q$.
Take the set $V_P$ of complete lifts of those vector fields on $P$.
Each set of vector fields acts transitively, and lifts the
previous, so each vector field in $V_P$ lifts a vector field
in $V_R$. By a theorem of Ehresmann \cite{Ehresmann:1961} or McKay \cite{McKay:2004a},
$P \to R$ is a fiber bundle mapping.
\end{proof}

\begin{theorem}\label{theorem:CompleteImage}
Suppose that $E_0 \to M_0$ and $E_1 \to M_1$
are Cartan geometries, and that $E_0 \to E_1$ is a morphism of Cartan geometries, modelled
on a model epimorphism. Suppose that $M_0$ is complete, and that $M_1$ is connected. Then $M_1$ is complete
and $E_0 \to E_1$ is a fiber bundle map.
\end{theorem}
\begin{remark}
We conjecture that $M_0 \to M_1$ is a fiber bundle map.
\end{remark}
\begin{proof}
Suppose that $F : E_0 \to E_1$ is a morphism modelled on an epimorphism.
Flow lines are carried to flow lines, so the constant flow on $E_1$ is complete on
the image of $F$. Therefore the image of $F$ is an $H_1$-invariant union
of path components of $E_1$, and so is all of $E_1$: $F$ is 
a surjective submersion and $E_1$ has complete flow.
Because the map $E_0 \to E_1$ takes the complete vector field
$\vec{A}$ to $\overrightarrow{\Phi(A)}$, the map $E_0 \to E_1$ is a fiber bundle map; see
Ehresmann \cite{Ehresmann:1961} or McKay \cite{McKay:2004a}. 
\end{proof}

\begin{lemma}\label{lemma:Uniqueness}
Pick Cartan geometries $E_0 \to M_0$ and $E_1 \to M_1$ with local models $\left(H_0,\mathfrak{g}_0\right)$ 
and $\left(H_1,\mathfrak{g}_1\right)$. Suppose that $M_0$ is connected.
Pick a local model morphism $\Phi : \left(H_0,\mathfrak{g}_0\right) \to \left(H_1,\mathfrak{g}_1\right)$.

A $\Phi$-morphism $E_0 \to E_1$, if one exists, is completely
determined by its value at a single point, i.e. if $F_0, F_1 : E_0 \to E_1$
are two morphisms, which agree at a single point of $E_0$, then $F_0=F_1$ 
everywhere on $E_0$.
\end{lemma}
\begin{remark}
So, roughly speaking, the space of morphisms is of dimension at most the dimension of $\mathfrak{g}_1$.
\end{remark}
\begin{proof}
On the graph of a morphism, for each $A \in \mathfrak{g}_0$,
the flow of the vector field $\vec{A}$ on $E_0$ intertwines with the flow of 
$\overrightarrow{\Phi(A)}$ on $E_1$. Therefore on $E_0 \times E_1$, the
graph of a morphism is acted on locally transitively by these
vector fields. This determines the morphism in a maximal
connected open subset of $E_0$. Then $H_0$-equivariance determines
the morphism at all points of $E_0$. 
\end{proof}

\begin{lemma}\label{lemma:LocalDescriptionOfCartanGeometry}
Suppose that $\Phi : \left(H_0, \mathfrak{g}_0\right) \to \left(H_1, \mathfrak{g}_1\right)$ 
is a local model morphism.
Suppose that $E_0 \to M_0$ and $E_1 \to M_1$ are Cartan geometries with those local models
and with Cartan connections $\omega_0$ and $\omega_1$. A submanifold $Z$ of $E_0 \times E_1$ 
is locally the graph of a local morphism just when 
\begin{enumerate}
\item $Z$ is an integral manifold of the Pfaffian system $\omega_1 = \Phi \omega_0$ and
\item
the obvious map $E_0 \times E_1 \to E_0$ restricts to a local diffeomorphism $Z \to E_0$.
\end{enumerate}
\end{lemma}
\begin{proof}
Since the problem is local, we can assume that $Z \subset E_0 \times E_1 \to E_0$
is a diffeomorphism to an open set $Z_0 \subset E_0$. So we can assume that $Z$ is the graph
of a map $F : E_0 \to E_1$. Because $\omega_1= \Phi \omega_0$ on $Z$, we see that
$F^*\omega_1= \Phi \omega_0$. Therefore $F_* \vec{A} = \vec{A}$.
Because the infinitesimal generators of the right $H_0$-action on $E_0$
are the vector fields $\vec{A}$ for $A \in \mathfrak{h}_0$,
$F$ is equivariant under the identity component of $H_0$. 

Since our problem is local, we can assume that $E_0$ and $E_1$ are trivial,
$E_0 = M_0 \times H_0$ and $E_1 = M_1 \times H_1$.
Denote the identity component of $H_0$ as $H_0^0$.
We can assume that $Z$ is contained
inside $M_0 \times H^0_0 \times M_1 \times H_1$. Therefore
we can extend the map $F$ by $H_0$-equivariance to be defined on $E_0$.
\end{proof}

\begin{definition}
Suppose that $\Phi : \left(H_0, \mathfrak{g}_0\right) \to \left(H_1,\mathfrak{g}_1\right)$
is a local model morphism.
Suppose that $E_0 \to M_0$ and $E_1 \to M_1$ are Cartan geometries
with those local models.
The \emph{$\Phi$-obstruction} between $E_0$ and $E_1$ is the function 
$K_1 \Lambda^2 \Phi - \Phi K_0 : E_0 \times E_1 \to \Lm{2}{\mathfrak{g}_0/\mathfrak{h}_0} \otimes \mathfrak{g}_1$.
Let $q : \mathfrak{g}_1 \to \mathfrak{g}_1/\Phi \mathfrak{g}_0$ be the obvious
quotient map. The \emph{$\Phi$-obstruction} of $E_1$ is the
function $q K_1 \Lambda^2 \Phi$.
\end{definition}

\begin{corollary}\label{corollary:LocalConstructionOfMorphisms}
Suppose that $\Phi : \left(H_0, \mathfrak{g}_0\right) \to \left(H_1, \mathfrak{g}_1\right)$ 
is a local model morphism. Suppose that $E_0 \to M_0$ and $E_1 \to M_1$ are Cartan geometries with 
those local models and with Cartan connections $\omega_0$ and $\omega_1$.
Suppose that the $\Phi$-obstruction vanishes at every point of $E_0 \times E_1$.
Consider the Pfaffian system $\omega_1 = \Phi \omega_0$.
Through each point of $E_0 \times E_1$ there passes a unique maximal connected 
integral manifold $Z$ of this Pfaffian system. Let $E_0'$ be the union of $H_0$-orbits of the points 
of $Z$. Then $H_0$ acts freely and properly on $E_0'$. Let $M_0'=E_0'/H_0$.
The composition $E_0' \to E_0 \times E_1 \to E_0$ is an $H_0$-equivariant
local diffeomorphism, quotienting to a local diffeomorphism $M_0' \to M_0$.
The 1-form $\omega_0$ pulls back to Cartan connection on $E_0' \to M_0'$,
and the map $E_0' \to E_0$ is a local isomorphism of Cartan geometries.
The composition $E_0' \to E_0 \times E_1 \to E_1$ is a $\Phi$-morphism.
Every $\Phi$-morphism $F : U_0 \to E_1$ from an open subset
$U_0 \subset E_0$ whose graph intersects $E_0'$
factors through $E_0' \to E_1$.
\end{corollary}
\begin{proof}
By lemma~\vref{lemma:LocalDescriptionOfCartanGeometry},
we can see locally that $E_0'$ is just going to be the graph of a
local isomorphism. Factoring of morphisms through $E_0' \to E_1$ 
is clear from lemma~\vref{lemma:Uniqueness}.
\end{proof}

\begin{theorem}\label{theorem:RealAnalyticCoveringExtension}
Suppose that $\Phi : \left(H_0, \mathfrak{g}_0\right) \to \left(H_1, \mathfrak{g}_1\right)$ 
is a local model morphism. 
Suppose that $E_0 \to M_0$ and $E_1 \to M_1$ are
real analytic Cartan geometries with those local models.
Suppose that $M_1$ is complete.

Every real analytic $\Phi$-morphism
from a connected open subset $U_0 \subset E_0$ to $E_1$ extends uniquely to a real analytic
$\Phi$-morphism $\tilde{E}_0 \to E_1$ from the pullback Cartan geometry $\tilde{E}_0 = p^* E_0$
of some covering map $p : \tilde{M}_0 \to M_0$ with $H_0$-invariant lift
\[
\xymatrix{
U_0 \ar[r] \ar[dr] & \tilde{E}_0 \ar[d] \\
                   & E_0.
}
\]
\end{theorem}
\begin{remark}
Ehresmann \cite{Ehresmann:1938} claimed a special case of this theorem, without
proof. Kobayashi \cite{Kobayashi:1954} also claimed without proof the
same special case, as a consequence of his incorrect main theorem
(which also appeared without proof).
\end{remark}
\begin{proof}
Suppose that the bundles of the Cartan geometries are $E_0 \to M_0$ and $E_1 \to M_1$,
with Cartan connections $\omega_0$ and $\omega_1$.
By hypothesis, we have a local morphism on some open sets. By 
theorem~\vref{theorem:GlobalMorphismOfFramings}, we can extend this local morphism to a morphism
of coframings $F : \tilde{E}_0 \to E_1$ on some covering space $\tilde{E}_0$ of 
$E_0$. Moreover, $\tilde{E}_0 \subset E_0 \times E_1$ is a submanifold (by construction
in theorem~\vref{theorem:GlobalMorphismOfFramings}), a union of integral
manifolds of the Pfaffian system $\omega_1=A \omega_0$.
Since the construction is $H_0$-invariant, the resulting coframing will
transform in the adjoint $H_0$-representation on $\tilde{E}_0$.
So the $H_0$-action on $\tilde{E}_0$ is free and proper,
because it is on $E_0 \times E_1$. Let $\tilde{M}_0 = \tilde{E}_0/H_0$.
Clearly $\tilde{M}_0$ is a smooth manifold and $\tilde{E}_0 \to \tilde{M}_0$ 
is a principal right $H_0$-bundle. Pull $\omega_0$ back to $E_0 \times E_1$,
and then restrict to $\tilde{E}_0$ to produce a Cartan connection.
Map $\tilde{E}_0 \to E_1$ by restricting the forgetful map $E_0 \times E_1 \to E_1$.
This map is a covering map, and $H_0$-equivariant, and therefore
descends to a covering map $\tilde{M}_0 \to M_0$.
\end{proof}

\subsection{The Frobenius--Gromov theorem for morphisms of Cartan geometries}

\begin{definition}
Suppose that $E \to M$ is a Cartan geometry
with local model $\left(H, \mathfrak{g}\right)$ and Cartan connection $\omega$.
Any function $f : M \to W$ to any vector space
$W$ has differential $df = \nabla f \omega$, say,
where $\nabla f : E \to W \otimes \mathfrak{g}^*$.
Clearly we can define $\nabla^2 f = \nabla \nabla f$,
and so on to all orders,
\[
 \nabla^j f : E \to W \otimes \mathfrak{g}^{*\otimes j}.
\]
\end{definition}

An easy calculation yields:
\begin{lemma}
 Suppose that $E \to M$ is a Cartan geometry with local model 
$\left(H, \mathfrak{g}\right)$. Suppose that $W$ is an $H$-module,
and $f$ is a smooth section of $E \times_H W$. Suppose that $\rho : H \to \GL{W}$
is the $H$-representation on $W$. In other words, $f : E \to W$ and
\[
r_h^* f = \rho(h)^{-1}f.
\]
Then under right $H$-action,
\[
 r_h^* \nabla^j f = \rho(h)^{-1} \otimes \left(\Ad(h)^{-1}\right)^{\otimes j}.
\]
\end{lemma}

\begin{definition}
The curvature of the Cartan geometry is given by is a function
\[
K : E \to \mathfrak{g} \otimes \Lm{2}{\mathfrak{g}/\mathfrak{h}}^*
\]
given by
\[
d \omega + \frac{1}{2}\left[\omega,\omega\right] = K \omega \wedge \omega.
\]
So we have defined $\nabla^0 K = K, \nabla^1 K, \nabla^2 K, \dots$ and
we see that the functions $\nabla^j K$ for $j<p$ determine the
$p$-jet of $K$. Clearly
\[
\nabla^j K : E \to \mathfrak{g} \otimes \Lm{2}{\mathfrak{g}/\mathfrak{h}}^* \otimes \bigotimes^{j} \mathfrak{g}^*
\]
\end{definition}

\begin{definition}
Suppose that $\Phi : \left(H_0, \mathfrak{g}_0\right) \to \left(H_1, \mathfrak{g}_1\right)$
is a local model morphism.
Suppose that $\pi_0 : E_0 \to M_0$ and $\pi_1 : E_1 \to M_1$ are real analytic
manifolds bearing real analytic Cartan geometries with those local models. 

Say that a point $e_0 \in E_0$ $\Phi$-\emph{hits} a point $e_1 \in E_1$
to order $p$ if 
\[
 \nabla^j K_1\left(e_1\right) \Lambda^{2} \Phi \otimes \bigotimes^{j} \Phi = \Phi \nabla^j K_0\left(e_0\right),
\]
for $j=0, 1, 2, \dots, p$.

Say that a point $e_0 \in E_0$ $\Phi$-\emph{hits} a point $e_1 \in E_1$
if there is an open set $U \subset M$ containing
the point $m_0=\pi_0\left(e_0\right)$ and a morphism
$F : \pi_0^{-1} U \to E_0$ so that $F\left(e_0\right)=e_1$.
We will say \emph{hits} instead of $\Phi$-hits
if $\Phi$ is understood.
\end{definition}

\begin{theorem}\label{theorem:FGMforCartanGeometryMorphisms}
Suppose that $\Phi : \left(H_0, \mathfrak{g}_0\right) \to \left(H_1, \mathfrak{g}_1\right)$
is a local model morphism. Suppose that $E_0 \to M_0$ and $E_1 \to M_1$ are real analytic
manifolds bearing real analytic Cartan geometries with those local models.

For each compact set $K \subset M_0 \times M_1$,
there is an integer $p$ so that if
a point $e_0 \in E_0$ hits a point $e_1 \in E_1$ to order $p$
and if $\left(\pi_0\left(e_0\right),\pi_1\left(e_1\right)\right) \in K$ 
then $e_0$ hits $e_1$.

In particular, the set of points $\left(e_0,e_1\right) \in E_0 \times E_1$
so that $e_0$ hits $e_1$ is an analytic variety.
\end{theorem}
\begin{proof}
Clearly it is enough to prove the result locally,
so we can assume that $E_0 = M_0 \times H_0$
and $E_1 = M_1 \times H_1$. Fix a compact
subset $K_0 \subset H_0$ and a compact
subset $K_1 \subset H_1$, each containing
a neighborhood of the origin. To
each compact set $K \subset M_0 \times M_1$,
associate the set $K'$
of points $\left(m_0,h_0,m_1,h_1\right)$
with $\left(m_0,m_1\right) \in K$ 
and $h_0 \in K_0$ and $h_1 \in K_1$.

We now have two notions of hitting: 
for the Cartan geometry and for
the coframings given by the
Cartan connections.
By lemma~\vref{lemma:LocalDescriptionOfCartanGeometry},
any coframing $\Phi$-morphism on an open
subset of $E_0$ is locally a
Cartan geometry $\Phi$-morphism.
So we can hit with the
coframings just when we can hit with
the Cartan geometries.

Clearly by theorem~\vref{theorem:FGMcoframings},
we can find an integer $p$ so that
the points hit to order $p$ lying
inside $K'$ are the points hit. 
By right invariance of the $p$-hitting and hitting criteria,
the same value of $p$ will work everywhere:
any point $e_0$ that hits a point $e_1$ to $p$-th order 
lies in the domain of a 
a (locally unique) coframing $\Phi$-morphism $F$
taking $e_0$ to $e_1$.
\end{proof}

\begin{remark}
In the notation of the previous theorem, the set of pairs 
$\left(e_0,e_1\right) \in E_0 \times E_1$ such that 
$e_0$ hits $e_1$ is an analytic variety,
invariant under the action
\[
\left(e_0,e_1\right)h_0=\left(e_0 h_0, e_1 \Phi\left(h_0\right)\right). 
\]
Therefore it projects to an analytic variety in the quotient $E_0 \times_{H_0} E_1$.

In the special case where $M_1=G_1/H_1$ is the homogeneous space with
its standard flat geometry, this analytic variety is $G_1$-invariant
too, so actually quotients to $M_0$: i.e. the set of
points of $M_0$ that hit one and hence every point in the model $G_1/H_1$
forms an analytic subvariety of $M_0$. But this is just the set of
points of $M_0$ near which there is a $\Phi$-morphism to $G_1/H_1$.
In particular, if the morphism $\Phi : G_0/H_0 \to G_1/H_1$ is
an immersion, then taking $M_1=G_1/H_1$, we find that 
then either $M_0$ is everywhere locally hitting
$G_0/H_0$ or there are no points of $M_0$ near which
there is a local $\Phi$-morphism to $G_1/H_1$.

In particular, if $\Phi=I$ is the identity, we recover
the fact that local homogeneity on a dense set implies
local homogeneity everywhere. 
\end{remark}

\begin{proposition}
Suppose that $\Phi : \left(H_0, \mathfrak{g}_0\right) \to \left(H_1, \mathfrak{g}_1\right)$
is a local model morphism.
Suppose that $M_0$ and $M_1$ are real analytic
manifolds bearing real analytic Cartan geometries
with those models. Suppose that $M_1$ has constant curvature. Then
the set of points of $M_0$ which hit some point of
$M_1$ (or equivalently hit a given point of $M_1$,
or equivalently hit every point of $M_1$) is a union of path components of $M_0$.
\end{proposition}
\begin{proof}
We can assume that $M_1$ has zero curvature, by replacing
its local model by a mutation. If a $\Phi$-morphism is defined near some point $m_0 \in M_0$,
say taking $m_0$ to some point $m_1 \in G_1/H_1$, then it is defined at all points nearby,
say on some open set $U_0 \subset M_0$. Its various translates under local
isomorphisms of $M_1$ form a family
of $\Phi$-morphisms so that we can get any point of $U_0$ to map to any
point of $M_1$ by some translate. Therefore the set of points of $M_0$
hitting any given point $m_1$ is open
and not empty. So above each point $m_0 \in U_0$, there is some point $e_0 \in E_0$ at which
\[
 0 = \Phi \nabla^j K_0\left(e_0\right),
\]
for $j=0,1,\dots,p$. This equation is $H_0$-invariant,
and therefore holds on an open $H_0$-invariant
subset of $E_0$. But then by analytic
continuation it holds everywhere above a union of path components of $M_0$.
\end{proof}

% \begin{corollary}
% Suppose that $\Phi : \left(H_0, \mathfrak{g}_0\right) \to \left(H_1, \mathfrak{g}_1\right)$
% is a gene morphism of homogeneous spaces $G_0/H_0$ and $G_1/H_1$.
% Suppose that $M_0$ is a real analytic
% manifold bearing a real analytic Cartan geometry
% modelled on $G_0/H_0$ with complete flow.
% Suppose that $M_0$ is connected and simply connected.
% Then either no point of $M_0$ hits any point
% of $G_1/H_1$, or there is a morphism
% $M_0 \to G_1/H_1$ globally defined.
% \end{corollary}

\begin{proposition}
Suppose that $\Phi : \left(H_0, \mathfrak{g}_0\right) \to \left(H_1, \mathfrak{g}_1\right)$
is an epimorphism of local models. Suppose that $M_0$ and $M_1$ are real analytic
manifolds bearing real analytic Cartan geometries
with those models. Suppose that $M_0$ has constant curvature. Then
the set of points of $M_1$ which are hit by some point of
$M_0$ (or equivalently are hit by a given point of $M_0$,
or equivalently are hit by every point of $M_0$) is a union of path components of $M_1$.
\end{proposition}
\begin{proof}
We can assume that $M_0$ has zero curvature, by replacing
its model by a mutation. If a $\Phi$-epimorphism hits some point $m_1 \in M_1$,
then it hits all points nearby,
say on some open set $U_1 \subset M_1$. Its various translates under local
isomorphisms of $M_0$ form a family
of $\Phi$-morphisms so that we can get any point of $M_0$ to map to any
point of $U_1$ by some translate. Therefore the set of points of $M_1$
hit by a given point $m_0$ is open
and not empty. So above each point $m_1 \in U_1$, there is some point $e_1 \in E_1$ at which
\[
0 =  \nabla^j K_1\left(e_1\right) \Lambda^{2} \Phi \otimes \bigotimes^{j} \Phi,
\]
for $j=0,1,\dots,p$. This equation is $H_1$-invariant,
and therefore holds on an open $H_1$-invariant
subset of $E_1$. But then by analytic
continuation it holds everywhere above a union of path components of $M_1$.
\end{proof}

\section{Stratification theorem for morphisms of Cartan geometries}

We recall the stratification theorem for algebraic group actions:
\begin{theorem}[Rosenlicht \cite{Rosenlicht:1963}]\label{theorem:Rosenlicht}
Pick a field $k$. If an algebraic group $G$ over $k$ acts algebraically
on an algebraic variety $V$ over $k$, then 
there is a dense $G$-invariant $k$-open subset $U \subset V$
so that $U/G$ is a algebraic variety over $k$, $U \to U/G$ a
regular morphism, and $k(U)^G=k\left(U/G\right)$.
\end{theorem}

\begin{corollary}\label{corollary:Rosenlicht}
If a real algebraic group $G$ acts algebraically
on a real algebraic variety $V$, then 
there is a collection of 
$G$-invariant algebraic subvarieties
$V= V_0 \supset V_1 \supset V_2 \dots$, 
of successively lower dimensions,
so that, if we let $U_j = V_j \backslash V_{j+1}$,
then 
\begin{enumerate}
\item $U_j/G$ is a smooth algebraic variety over $k$, and 
\item
the map $U_j \to U_j/G$ is a smooth submersion for each $j$ and
\item
the function fields are $\R{}(U)^G=\R{}\left(U/G\right)$ and
\item
each variety $U_j/G$ is affine.
\end{enumerate}
\end{corollary}
\begin{proof}
We pick an open set $U_0 \subset V_0=V$ as in theorem~\vref{theorem:Rosenlicht},
and set $V_1 = V_0 \backslash U_0$. But then take the points where the map
$U_0 \to U_0/G$ is not a submersion, and add them to $V_1$, and thereby remove
them from $U_0$. Repeat inductively. If $U_1/G$ is not affine, then
cut out a proper subvariety from it, say $W_1 \subset U_1/G$, so that
$U_1/G \setminus W_1$ is affine. Then let $W_1'$ be the preimage of
$W_1$ in $U_1$, and let $U_1' = U_1 \setminus W_1'$ and $V_2 = V_1 \cup W_1'$,
etc. to inductively arrange that all of the quotient varieties are affine.
\end{proof}

Recall that a rational map $f : X \to Y$ from a real analytic space
$X$ to an affine variety $Y$ is a map defined away from some
analytic subvariety so that the composition with a rational function is rational,
i.e. locally a quotient of real analytic functions.

\begin{definition}
Suppose that $V$ is a real analytic variety and $X$ is a real algebraic variety.
A \emph{striating map}
is map $I : V \to X$ (not necessarily continuous)
with a choice of closed subvarieties $V=V_0 \supset V_1 \supset V_2 \dots$
(called the \emph{strata}) of successively lower dimension and closed affine
subvarieties $X_0, X_1, X_2 \dots \subset X$ so that
\begin{enumerate}
\item 
the varieties $X_1, X_2, \dots$ are affine and
\item
if we let $U_j = V_j \backslash V_{j+1}$, then $U_j$ is smooth and 
\item
$\left.I\right|_{U_j} : U_j \to X$ is a smooth map of locally constant rank and
\item
$\left.I\right|_{U_j} : U_j \to X$ has image inside $X_j \subset X$ and
\item
$\left.I\right|_{U_j} : V_j \to X_j$ is a rational map for each $j$.
\end{enumerate}
\end{definition}

\begin{remark}
A striating map endows each set $U_j$ with a smooth analytic foliation: the level sets of $I$. 
So $V$ has a smooth foliation (the level sets of $I$) except on a subvariety $V_1$,
but $V_1$ itself has a smooth foliation (again the level sets of $I$) except on a subvariety $V_2$, etc.
\end{remark}

\begin{remark}
We will consider the closed subvarieties $V_j$ and $X_j$ to be part of the data of the striating map.
\end{remark}

\begin{example}
By Rosenlicht's theorem,  if a real algebraic group $G$ acts algebraically
on a real algebraic variety $V$, then there is a $G$-invariant striating map.
On each stratum, this striating map is algebraic and submersive. 
We can take $X$ to be just the formal disjoint union of the
various $U_j/G$, so that each $U_j/G$ is a closed and open subset of $X$.
It is not clear if the hypotheses of Rosenlicht's theorem ensure that there is one 
striating map through which all others factor.
\end{example}

\begin{example}
Suppose that $I : V \to X$ is any regular map of real algebraic varieties.
We can then let $V_1$ be the set of points where $V$ is not smooth,
together with those points where $I$ is not of maximal rank.
Restrict $I$ to $V_1$ and let $V_2$ be the set of points of $V_1$
where $V_1$ is singular or $I$ is not of maximal rank, etc.
We set $U_j = V_j \backslash V_{j+1}$.
In this way $I$ is canonically a
striating map. Clearly not every striating map comes about this way.
\end{example}

\begin{lemma}
 If $I : M \to X$ is a striating map, and $G$ is a 
real analytic Lie group with real analytic action on $M$ with 
smooth quotient $M/G$, and if $I$ is $G$-invariant,
then $I$ descends to a striating map $I : M/G \to X$.
\end{lemma}
\begin{proof}
The strata of $I$ are $G$-invariant so descend
to analytic subvarieties, as does the map $I$. 
Descent doesn't alter rank of invariant maps.
\end{proof}

\begin{definition}
A \emph{refinement} $I' : V \to X$ of a striating map $I : V \to X$
is a striating map for which $I'=I$ at every point of $V$, but
whose strata $V_j'$ are strictly larger: $V_j \subset V_j'$.
\end{definition}

\begin{lemma}\label{lemma:ComposeStriatingMap}
Suppose that $I : V \to X$ is a striating map invariant under a real analytic group action 
of a group $G$ on a real analytic variety $V$ and $F : W \to V$
is a $G$-equivariant real analytic map between real analytic varieties
with real analytic $G$-action. 
Then there is a refinement $J : W \to X$ of $I \circ F$ 
which is a $G$-invariant striating map.
\end{lemma}
\begin{proof}
For the moment, replace $X$ by the algebraic Zariski closure of the image of $I \circ F$,
and $V$ by the preimage of that algebraic Zariski closure under $I$.
So we can assume that $F$ has image striking $U_0$. Therefore $I \circ F : F^{-1} U_0 \to X_0$
is well defined. Then take the set of points of $F^{-1} U_0$ where $I \circ F$ has
rank locally constant, say $W_0$, as our first stratum. Repeat inductively on
the complement of $W_0$. Note that the complement of $W_0$ might have some components
mapping under $I \circ F$ into $X_0$, so we might have to attach a copy
of $X_0$ to $X_1$. This is consistent with the definition of striating map. 
\end{proof}

\begin{definition}
A local model $\left(H,\mathfrak{g}\right)$ is \emph{algebraic} if
$H$ is a linear algebraic group and the representation of $H$ on $\mathfrak{g}$ is algebraic.
\end{definition}

\begin{definition}
A Cartan geometry is of \emph{algebraic type} if its local model is algebraic.
\end{definition}

\begin{remark}
A Cartan geometry of algebraic type need not be in any
sense algebraic. For example, all pseudo-Riemannian
geometries are of algebraic type.
\end{remark}

\begin{theorem}
Suppose that $M_0$ and $M_1$ are real analytic
manifolds bearing real analytic Cartan geometries.
Suppose that the Cartan geometry on  $M_1$ is of algebraic type.
Then the set of pairs $\left(m_0,m_1\right) \in M_0 \times M_1$
so that $m_0$ hits $m_1$ (for a given local model morphism) is an analytic
subvariety of $M_0 \times M_1$.

Suppose that $W_0 \subset M_0$ and $W_1 \subset M_1$
are relatively compact open sets.
Then there is a pair of striating maps
$I_0 : W_0 \to X$ and $I_1 : W_1 \to X$,
so that a point $m_0 \in W_0$ hits a point $m_1 \in W_1$ just when
$I_0\left(m_0\right)=I_1\left(m_1\right)$.
So a point of $W_0$ hits a point of $W_1$ just when 
the points lie on corresponding leaves on corresponding strata.
These striating maps can be chosen to be invariant under the pseudogroups
of local isomorphisms of the Cartan geometries.
\end{theorem}
\begin{proof}
The problem of determining whether $m_0$ hits $m_1$ 
is the same problem as asking if the quantity
\[
 \nabla^j K_1\left(e_1\right) \Lambda^{2} \Phi \otimes \bigotimes^{j} \Phi 
\]
can be made equal to the quantity
\[
\Phi \nabla^j K_0\left(e_0\right),
\]
by $H_0 \times H_1$-action, for $j=1,2,\dots,p$
at some points $e_0$ and $e_1$ over points $m_0$ and $m_1$ respectively,
for some sufficiently large integer $p$.
This integer $p$ will depend on the choice of sets $W_0$
and $W_1$.

% (The reader might wonder why $H_1$ plays a different role than $H_0$
% in this theorem. Each point $e_0$ has to go somewhere under any
% morphism, so to get a morphism we can pick a particular point $e_0$,
% and then try to find a point $e_1$ that it hits. So we only
% need to move around $e_1$ with the $H_1$-action
% to see if we can find a point that $e_0$ hits.)

Let $\kappa^j_1 = \nabla^j K_1\left(e_1\right) \Lambda^{2} \Phi \otimes \bigotimes^{j} \Phi$.
Let $\kappa^j_0 = \Phi \nabla^j K_0\left(e_0\right)$.
Let
\[
J_0 = 
\left( \kappa^0_0, \kappa^1_0, \dots, \kappa^p_0 \right) 
\]
and
\[
J_1 = 
\left( \kappa^0_1, \kappa^1_1, \dots, \kappa^p_1 \right) 
\]
So $J_1$ is a function on $E_1$ valued in some $H_1$-module, say $V$, and is $H_1$-equivariant.
We must determine if the $H_1$-orbit in $V$ of $J_1$
is the same as the $H_1$-orbit in $V$ of $J_0$. 

% (The reader is probably confused about the $H_1$-action
% on $J_0$ since $J_0$ is a function on $E_0$ and
% $H_1$ doesn't act on $E_0$. What we really mean is that we are
% looking at the $H_1$-orbit of the value of $J_0$ at a point. 
% This makes sense because $J_0$ is valued in an $H_1$-module,
% even though $E_0$ is not an $H_1$-bundle.
% This is not a problem, because we have a morphism $H_0 \to H_1$,
% so that $H_0$ actually acts on this 
% $H_1$-module.)

By corollary~\vref{corollary:Rosenlicht}, there is
an $H_1$-invariant striating map $I : V \to X$,
which is a submersion on each stratum, which distinguishes
orbits of the $H_1$-action on $V$. Therefore
$I_1=I \circ J_1 : E_1 \to X$ is an $H_1$-invariant 
map. But it might not be striating, because its
rank could change at various points. However,
it is close to what we need:
while it might not be striating, its level sets 
are precisely the points of $E_1$
at which $J_1$ takes on values in a given $H_1$-orbit.
So our next problem is to refine the maps $I_1$ and $I_0=I \circ J_0$ 
to become striating maps, i.e. locally constant rank on each stratum.

By lemma~\vref{lemma:ComposeStriatingMap},
we can assume that $I_0$ is a striating map by refinement,
and the same for $I_1$. The strata of $X$ could be
altered in each of these processes of refinement, but
(as we see in the proof of lemma~\vref{lemma:ComposeStriatingMap})
only by moving various strata $X_i$ into various other lower strata $X_{i+1}$
repeatedly. So we can assume without loss of generality
that $I_0$ and $I_1$ have the same strata
in $X$. Moreover, we can ensure that $I_0=I_1$
precisely at points at which $I \circ J_0=I \circ J_1$.

So now both $I_0$ and $I_1$ are striating, and $H_1$-invariant,
so descend to striating maps $I_0 : M_0 \to X$ and $I_1 : M_1 \to X$ respectively.
Moreover, these maps have level sets given precisely
by the points of $W_0$ or $W_1$ above which
$J_0$ or $J_1$ reach particular $H_1$-orbits.
Therefore a point $m_0 \in M_0$ hits a
point $m_1 \in M_1$ just when $I_0\left(m_0\right)=I_1\left(m_1\right)$.
\end{proof}

\begin{remark} 
 In the previous theorem, if $H_1$ is compact, then its
space of invariants is a suitable choice of $X$ and $I_0$ and $I_1$
can be chosen to be the quotient mapping to $X$; 
see Procesi \cite{Procesi:2007} p. 556, theorem 2.
In particular, both $I_0$ and $I_1$ will then
be smooth real analytic maps to affine space,
with image inside the affine algebraic variety $X$.
\end{remark}

\begin{corollary}
Suppose that $M_0$ and $M_1$ are real analytic
manifolds bearing real analytic Cartan geometries
with those models.
Suppose that $M_1$ is of algebraic type. Pick a point $m_0 \in M_0$. 
Let $Z$ be the set of points $m_1 \in M_1$ so that $m_0$ hits $m_1$.
Then $Z$ is dense if and only if $Z$ is Zariski dense.
\end{corollary}

\begin{remark}
For simplicity, we stated the results in this section for real manifolds and varieties,
but their complex analytic analogues will have obvious formulations and identical proofs.
\end{remark}

\begin{corollary}[Dumitrescu \cite{Dumitrescu:2001}]
Suppose that $M_0$ and $M_1$ are complex
manifolds bearing holomorphic Cartan geometries.
Suppose that $M_1$ is of algebraic type. 
\begin{enumerate}
 \item 
Suppose that $M_0$ admits no nonconstant meromorphic functions. 
Pick a point $m_1 \in M_1$.  If some point of $M_0$ 
hits $m_1$, then all points of $M_0$ hit $m_1$.
\item
Suppose instead that $M_1$ admits no nonconstant meromorphic functions. 
Pick a point $m_0 \in M_0$.  If some point of $M_1$ 
is hit by $m_0$, then all points of $M_1$ are hit by $m_0$.
\end{enumerate}
\end{corollary}
\begin{proof}
Since there are no nonconstant meromorphic functions 
on $M_0$ (or $M_1$), the striating map $I_0$ (or $I_1$) must be constant.
\end{proof}

\begin{example}
Suppose that $M_1$ is a real analytic manifold with a real analytic projective connection. 
The set of points of $M_1$ through which there are flat totally geodesic surfaces is a closed 
real analytic subvariety of $M_1$. Indeed we let $M_0=\Proj{2}$ imbedded as a 
linear subspace in $\Proj{n}$, let $G_0$ be the subgroup of $\PSL{n+1,\R{}}$
preserving $M_0$, and let $H_0$ be the subgroup of $G_0$ fixing a point.
\end{example}

\begin{example}
 Complexifying the last example: suppose that $M_1$ is a complex manifold
with holomorphic projective connection. Suppose that $M_1$ has no nonconstant
meromorphic functions. Then either every point of $M_1$ lies on a flat
totally geodesic complex surface, or no point does.
\end{example}

\section{Disguises: morphisms covering the identity map}

\begin{remark}
 It is unclear how many different morphisms $E_0 \to E_1$
between two fixed Cartan geometries $E_0 \to M_0$ and $E_1 \to M_1$
could induce the same underlying map $M_0 \to M_1$.
\end{remark}

We want a weaker notion than isomorphism, allowing the Cartan geometry to change on a fixed manifold.

\begin{definition}
Suppose that $\Phi : \left(H_0,\mathfrak{g}_0\right) \to \left(H_1,\mathfrak{g}_1\right)$
is a local model morphism, and that the induced linear map
$\Phi : \mathfrak{g}_0/\mathfrak{h}_0 \to \mathfrak{g}_1/\mathfrak{h}_1$
is an isomorphism of $H_0$-modules. We will call $\Phi$ an \emph{disguise}.
\end{definition}
In our earlier terminology, an disguise is a \emph{base isomorphism}.

\begin{theorem}\label{theorem:ProjMorphism}
Pick an disguise $\Phi : \left(H_0,\mathfrak{g}_0\right) \to \left(H_1,\mathfrak{g}_1\right)$.
Pick a $\left(H_0,\mathfrak{g}_0\right)$-geometry $E_0 \to M$.
Let 
\[
E_1 = E_0 \times_{H_0} H_1.
\]
Let $\omega_1$ be the 1-form on $E_0 \times H_1$ which at a point $\left(e_0,h_1\right)$ is
\[
\omega_1 = h_1^{-1} \, dh_1 + \Ad\left(h_1\right)^{-1} \Phi \omega_0.
\]
Then $\omega_1$ is the pullback of a unique $1$-form, which we also denote by $\omega_1$,
on $E_1$. Moreover $E_1 \to M$ is a $\left(H_1,\mathfrak{g}_1\right)$-geometry. 
Define $F_{\text{disguise}} : E_0 \to E_1$ given by 
composition of $e_0 \in E_0 \mapsto \left(e_0,1\right) \in E_0 \times H_1$
with $E_0 \times H_1 \to E_1$. This map $F_{\text{disguise}}$ is a $\Phi$-morphism,
called the $\Phi$-\emph{disguise} of $E_0$. 

Lets write $\bar{\omega}_0$ to denote the soldering form
$\omega_0 + \mathfrak{h}_0 \in \nForms{1}{E_0} \otimes \mathfrak{g}_0/\mathfrak{h}_0$.
Then the curvature of the disguise is
\[
 K_1 \bar{\omega}_1 \wedge \bar{\omega}_1
=
\Ad\left(h_1\right)^{-1} 
\left(
  K_0 \bar{\omega}_0 \wedge \bar{\omega}_0
+
\frac{1}{2}
\left[
  \Phi \omega_0,
  \Phi \omega_0
\right]
-
\Phi
\left[
  \omega_0,
  \omega_0
\right]
\right).
\]

Every $\Phi$-morphism $F$ of Cartan geometries factors uniquely as $F=F_1 F_{\text{disguise}}$
where $F_1$ is a uniquely determined local isomorphism of $\left(H_1,\mathfrak{g}_1\right)$-geometries.
\end{theorem}
\begin{remark}
 Our terminology is a bit vague, as to whether $\Phi$ or $F$ is the disguise, but it should be clear that there is
essentially a unique Cartan geometry disguise $F_{\text{disguise}}$ modelled on a given model disguise.
\end{remark}
\begin{proof}
The right action of $H_0$ on $E_0 \times H_1$ is 
\[
\left(e_0, h_1\right)h_0
=
\left(e_0 h_0, \Phi\left(h_0\right)^{-1} h_1\right), h_0 \in H_0.
\]
This action commutes with the right $H_1$-action
\[
\left(e_0, h_1\right)k_1
=
\left(e_0 h_0, h_1 k_1\right), k_1 \in H_1.
\]
The 1-form $\omega_1$ is invariant under the right $H_0$-action,
and transforms in the adjoint representation under the right $H_1$-action. Moreover $\omega_1$
vanishes on the orbits of the $H_0$-action, and therefore descends to a 1-form on $E_1$.
Clearly $\omega_1$ is a Cartan connection on $E_1$.
The curvature is a simple calculation.

Define $F_{\text{disguise}} : E_0 \to E_1$ given by 
composition of $e_0 \in E_0 \mapsto \left(e_0,1\right) \in E_0 \times H_1$
with $E_0 \times H_1 \to E_1$. This map $F=F_{\text{disguise}}$ 
is equivariant for the $H_0$-action, and $F^* \omega_1 = \Phi \omega_0$, so
a $\Phi$-morphism. 

Suppose that $F : E_0 \to E_1'$ is a $\Phi$-morphism,
where $E_1' \to M_1'$ is a $G_1/H_1$-geometry.
Define a map (which we also denote by $F$) say $F : E_0 \times H_1 \to E_1'$,
by $F\left(e_0,h_1\right)=F\left(e_0\right)h_1$.
This map is clearly a submersion to $E_1'$;
its fibers are precisely the right $H_0$-orbits.
This map is invariant under the right $H_0$-action,
so descends to a local diffeomorphism $F_1 : E_1 \to E_1'$.
This local diffeomorphism is $H_1$-equivariant.

Pulling back the Cartan connection $\omega'_1$ on $E_1'$ to $E_0$, because $E_0 \to E_1'$
is a morphism,
\[
F^* \omega_1' = \Phi \omega_0.
\]
But then pulling back to $E_0 \times H_1$,
\[
F^* \omega_1' = h_1^{-1} \, dh_1 + 
\Ad\left(h_1\right)^{-1} \, \Phi \omega_0 = \omega_1.
\]
Therefore the map $F_1 : E_1 \to E_1'$ is a local isomorphism of $G_1/H_1$-geometries.
\end{proof}

\begin{example}[From affine connections to projective connections]
Let
\[
G_0 = \GL{n,\R{}} \rtimes \R{n}, \ H_0=\GL{n,\R{}},  
\]
so $G_0/H_0$ is affine $n$-space. 
Write elements of $G_0$ as
matrices
\[
g_0 =
\begin{pmatrix}
 1 & 0 \\
 x_0 & h_0
\end{pmatrix}
\]
and elements of $H_0$ as
\[
\begin{pmatrix}
 1 & 0 \\
 0 & h_0
\end{pmatrix}.
\]

Let $G_1=\PSL{n+1,\R{}}$,
let $H_1$ be the subgroup fixing the point
$\left[0,0,\dots,0,1\right] \in \RP{n}$.  
Let $\Phi$ be the
obvious composition $G_0 \to \GL{n+1,\R{}} \to \PSL{n+1,\R{}}$;
then $\Phi$ is a disguise of homogeneous spaces.
This disguise takes a manifold with affine connection
to the same manifold with the induced projective connection.
This is \emph{not} necessarily the normal projective connection.
By theorem~\vref{theorem:ProjMorphism},
there is a unique morphism associating a projective connection
to an affine connection.
\end{example}

\begin{example}[From Riemannian metrics to affine connections]
Let $G_0=\Or{n} \rtimes \R{n}$ and $H_0=\Or{n}$.
Let $G_1=\GL{n,\C{}} \rtimes \R{n}$ and $H_1 = \GL{n,\C{}}$.
The obvious map $\Phi : G_0/H_0 \to G_1/H_1$
is a disguise, which takes a Riemannian manifold
and forgets the Riemannian
metric, remembering only the Levi--Civita connection.
\end{example}

\begin{example}[Induced effective geometries]\label{example:Effectivedisguise}
Suppose that $N \subset H$ is a closed
subgroup which is a normal subgroup of $G$.
We have an obvious disguise
$\Phi : \left(H,\mathfrak{g}\right) \to \left(H/N,\mathfrak{g}/\mathfrak{n}\right)$.
By theorem~\vref{theorem:ProjMorphism},
there is a unique disguise
between any $G/H$-geometry and some
unique $(G/N)/(H/N)$-geometry
on the same manifold.
In particular, if we take
$N$ to be the largest normal subgroup
of $G$ contained in $H$, then 
the resulting disguise
could be called the \emph{induced effective Cartan
geometry}; see Sharpe \cite{Sharpe:2002}
for the definition of an effective Cartan geometry.
\end{example}

\begin{example}[Fefferman constructions]
 If $G_0 \subset G_1$ and $H_0$ contains
$G_0 \cap H_1 \subset H_0$, and the
composition $\mathfrak{g}_0 \to \mathfrak{g}_1 \to \mathfrak{g}_1/\mathfrak{h}_1$
is surjective, then the associated morphism 
of homogeneous spaces is called
a \emph{Fefferman} morphism of homogeneous
spaces. A morphism of Cartan geometries
modelled on a Fefferman morphism is called
a \emph{Fefferman construction}; see 
\v{C}ap \cite{Cap:2006}.
By lemma~\vref{lemma:liftIsCanonical} and theorem~\ref{theorem:ProjMorphism},
there is a unique Fefferman construction with a given model:
the composition of the lift
modelled on $G_0/\left(G_0 \cap H_1\right) \to G_0/H_0$
with the obvious disguise modelled on 
$G_0/\left(G_0 \cap H_1\right) \to G_1/H_1$.
\end{example}

\begin{remark}
Let's say that a local model morphism 
$\Phi : \left(H_0,\mathfrak{g}_0\right) \to \left(H_1,\mathfrak{g}_1\right)$
is \emph{target natural} if to every manifold with 
$\left(H_0,\mathfrak{g}_0\right)$-geometry there is 
associated a $\Phi$-morphism to a manifold with 
$\left(H_1,\mathfrak{g}_1\right)$-geometry, 
and this association is invariant under local
isomorphisms. Clearly we have examples
of target natural morphisms from Riemannian
geometry and projective connections above.
Our theorem above says that disguises
of homogeneous space genes are target natural.
They are the only known target natural
morphisms, besides the construction
of a normal projective connection
from any projective connection,
and the construction of a torsion-free connection
from any affine connection.

Let's say that a local model morphism 
$\Phi : \left(H_0,\mathfrak{g}_0\right) \to \left(H_1,\mathfrak{g}_1\right)$
is \emph{source natural} if to every manifold with 
$\left(H_1,\mathfrak{g}_1\right)$-geometry, there is 
associated a manifold with $\left(H_0,\mathfrak{g}_0\right)$-geometry,
and a $\Phi$-morphism between these manifolds,
and this association is invariant under local
isomorphisms. 
The lift of a Cartan geometry is modelled on a 
source natural morphism. The classification
of source and target natural morphisms
of local models is not known or even conjectured.
\end{remark}

\subsection{Disguises and reductive structure group}

\begin{definition}
Lets say that a local model $\left(H,\mathfrak{g}\right)$ is
\emph{nearly reductive} if $\mathfrak{g}$ splits as an $H$-module
into a sum $\mathfrak{g}=\mathfrak{h} \oplus \mathfrak{g}/\mathfrak{h}$.
\end{definition}

\begin{remark}
 Sharpe \cite{Sharpe:2002} calls a local model
\emph{reductive} under this condition.
Sharpe's terminology is almost certain to confuse the 
reader into thinking that $H$ (or perhaps
$G$) is supposed to be a reductive linear algebraic group.
Sharpe admits (p. 197 footnote 17) that the terms
could be confused in this way.
\end{remark}

\begin{definition}
If $\pi : E \to M$ is a Cartan
geometry modelled on a nearly reductive model $\left(H,\mathfrak{g}\right)$, then the
splitting $\mathfrak{g}=\mathfrak{h} \oplus \mathfrak{g}/\mathfrak{h}$
splits the Cartan connection $\omega$ into a sum $\omega=\omega^H + \omega^{\perp}$
with $\omega^H$ a 1-form valued in $\mathfrak{h}$ and $\omega^{\perp}$ a 1-form valued in $V$.
\end{definition}

\begin{example}
The splitting of the Cartan connection is an disguise, modelled on the model disguise
\[
\Phi : \left(H, \mathfrak{g}\right) \to \left(H,\mathfrak{h} \oplus \mathfrak{g}/\mathfrak{h}\right).
\]
The 1-form $\omega^H$ is a connection for the bundle $E \to M$.
\end{example}

\section{Morphisms as integral manifolds of a Pfaffian system}

\begin{definition}
If $E \to M$ is a Cartan geometry with Cartan connection $\omega$,
with model $\left(H,\mathfrak{g}\right)$, then the \emph{curvature}
of the Cartan geometry is the function 
$K : E \to \mathfrak{g} \otimes \Lm{2}{\mathfrak{g}/\mathfrak{h}}^*$
for which
\[
d \omega + \frac{1}{2}\left[\omega,\omega\right] = \frac{1}{2} K \omega \wedge \omega.
\]
(In this definition, $\omega \wedge \omega$ is the form
valued in $\Lm{2}{\mathfrak{g}/\mathfrak{h}}$ given
by $\omega \wedge \omega(v,w)=\omega(v) \wedge \omega(w)$.)
\end{definition}

Three natural problems arise in studying morphisms of Cartan geometries:
\begin{enumerate}
\item
 How many morphisms are there from 
a fixed Cartan geometry to a fixed Cartan geometry?
\item
How many morphisms are there from a fixed
Cartan geometry to all possible Cartan geometries?
\item
How many morphisms are there from all possible
Cartan geometries to a fixed Cartan geometry?
\end{enumerate}

Each of these vague questions leads to a different
Pfaffian system, whose integral manifolds are locally
the graphs of morphisms of the required form. 

\begin{example}
Suppose that $\Phi : \left(H_0,\mathfrak{g}_0\right) \to \left(H_1,\mathfrak{g}_1\right)$ 
is a local model morphism. Suppose that $E_0 \to M_0$
and $E_1 \to M_1$ are Cartan geometries with these models,
with Cartan connections $\omega_0$ and $\omega_1$ respectively.
On $E_0 \times E_1$, consider the Pfaffian system $\omega_1 = \Phi \omega_0$.
We denote by 
$\Lambda^2 \Phi : \Lm{2}{\mathfrak{g}_0/\mathfrak{h}_0} \to \Lm{2}{\mathfrak{g}_1/\mathfrak{h}_1}$
the linear map $\Lambda^2 \Phi (v \wedge w)=F(v) \wedge F(w)$. 
The torsion of the Pfaffian system at any integral element
is 
\(
K_1 \Lambda^2 \Phi - \Phi K_0.
\)

Let $V \subset E_0 \times E_1$ be the subset of points on which
$K_1 \Lambda^{2} \Phi=\Phi K_0$. Suppose that $V^{\text{sm}} \subset V$
is some open subset which is locally a submanifold
of $E_0 \times E_1$. Then the Pfaffian system restricts
to $V^{\text{sm}}$ to satisfy the conditions of the Frobenius
theorem. Therefore if $V^{\text{sm}}$ is not empty, then it
is foliated by integral manifolds of the Pfaffian system. The
1-form $\omega_0$ restricts to a coframing on one of these integral manifolds
just when that integral manifold is locally the graph of a local
morphism.
Clearly the integral manifolds living inside $V^{\text{sm}}$
therefore form a finite dimensional family.
Indeed if the local models arise from homogeneous space modesl,
then this family is never larger in dimension than
the associated family of morphisms of the models.

Keep in mind that there might be singular points
of $V$ through which integral manifolds pass. If the Cartan geometry is real analytic, 
then $V$ is stratified into smooth submanifolds, so that
a Zariski open subset of each integral manifold will have to lie inside one of
the strata, and so we could in principal find all integral
manifolds in this way. In the real analytic category, we can therefore
justify the claim that the local morphisms from open subsets
of $M_0$ to open subsets of $M_1$ form a finite dimensional family
of no larger dimension than the associated family of morphisms of
the models.

For example, if the homogeneous space $G_0/H_0$
is one dimensional, then $\Lambda^2 \Phi=0$, 
and $K_0=0$, so there is no torsion and the Pfaffian system
foliates $E_1$ by integral manifolds, which are locally
the graphs of local morphisms from open sets of $G_0/H_0$.
\end{example}

\begin{example}\label{example:AnyMorphisms}
Instead of trying to construct morphisms between specific Cartan geometries, 
we could ask if there are any morphisms from any Cartan geometries
to a fixed Cartan geometry. 
Suppose that $\Phi : \left(H_0,\mathfrak{g}_0\right) \to \left(H_1,\mathfrak{g}_1\right)$ 
is a local model morphism. 
Suppose that 
$\Phi : \mathfrak{g}_0 \to \mathfrak{g}_1$ is an injective linear map.
Suppose that $E_1 \to M_1$ is a $\left(H_1,\mathfrak{g}_1\right)$-geometry.
Let $q : \mathfrak{g}_1 \to \mathfrak{g}_1/\Phi\left(\mathfrak{g}_0\right)$
be the obvious quotient map. On $E_1$,
consider the Pfaffian system $q \omega_1 = 0$.
The torsion of this Pfaffian system is 
\[
q K_1 \Lambda^{2} \Phi 
\]
Again we let $V$ be the set of points on which the torsion vanishes,
and let $V^{\text{sm}}$ be any open subset of $V$ which is locally a
submanifold of $E_1$. Again, on $V^{\text{sm}}$ the Pfaffian
system satisfies the conditions of the Frobenius theorem,
so is foliated by integral manifolds. Integral manifolds of the Pfaffian
system on which $\omega_1$ is a  coframing (valued in $\mathfrak{g}_0$!)
are locally graphs of local morphisms from submanifolds
with $\left(H_0,\mathfrak{g}_0\right)$-geometry.

Similar remarks to those of the previous example apply concerning
integral manifolds through singular points of $V$.
\end{example}

\section{Morphisms and complete flows}

\begin{definition}
 A local model morphism $\Phi : \left(H_0,\mathfrak{g}_0\right) \to \left(H_1,\mathfrak{g}_1\right)$
is called \emph{immersive downstairs} if the induced linear map
$\Phi : \mathfrak{g}_0/\mathfrak{h}_0 \to \mathfrak{g}_1/\mathfrak{h}_1$
is injective.
\end{definition}

\begin{proposition}
Pick a homogeneous space morphism 
\[
\Phi : G_0/H_0 \to G_1/H_1
\]
which is immersive downstairs. Pick a $\left(H_1,\mathfrak{g}_1\right)$-geometry $E_1 \to M_1$.
Let 
\[
q : \mathfrak{g}_1 \to \mathfrak{g}_1/\Phi\left(\mathfrak{g}_0\right)
\]
be the obvious linear quotient map. Suppose that 
\begin{enumerate}
\item $q K_1 \Lambda^2 \Phi=0$ and
\item $\LieDer_{\overrightarrow{\Phi(A)}} K_1 \Lambda^2 \Phi=0$ for each $A \in \mathfrak{g}_0$
and 
\item $\overrightarrow{\Phi(A)}$ is a complete vector field, for each $A \in \mathfrak{g}_0$.
\end{enumerate}
Then each point of $E_1$ lies in the image of a morphism from a
mutation $G_0'/H_0$ of a covering of $G_0/H_0$. This morphism is unique up to $G_0'$-action.
In particular, if $K_1 \Lambda^2 \Phi=\Phi \left[,\right]_{\mathfrak{g}_0}$,
we can take $G_0'$ to be a covering of $G_0$.
\end{proposition}
\begin{proof}
From example~\ref{example:AnyMorphisms}, we can see that (1) ensures
that $E_1$ is foliated by images of $\Phi$-morphisms. Next, (2)
ensures that each of the leaves is a Cartan geometry with constant curvature.
Any Cartan geometry satisfies
\[
\left[\vec{A},\vec{B}\right]=\overrightarrow{\left[A,B\right]}+\overrightarrow{K\left(A,B\right)},
\]
so that a Cartan geometry with constant curvature is a mutation
for the bracket
\[
 \left[A,B\right]'=\left[A,B\right]+K\left(A,B\right).
\]
Condition (3) ensures completeness, so that the
geometry is the image of a local isomorphism from its model.
The model must have the same structure group $H_0$,
being only a mutation. We can replace the model by a covering
space to arrange that the model is a covering of 
a mutation of $G_0/H_0$.
\end{proof}

\begin{remark}
 The image of an immersive downstairs morphism
is an immersed submanifold. If it is complete and flat, then it is a quotient
of its model $G_0/H_0$, up to replacing the model by a 
covering space $G_0'/H_0$. The quotient is by a group
action of a discrete subgroup $\Gamma \subset G_0'$.
\end{remark}

\section{Hartogs extension of Cartan geometry morphisms}

\begin{lemma}[Matsushima and Morimoto \cite{Matsushima/Morimoto:1960}]
A complex Lie group is Stein if and only if the identity
component of its center contains no complex torus.
\end{lemma}

\begin{proposition}
Pick 
\begin{enumerate}
\item
a holomorphic local 
model morphism $\Phi : \left(H_0,\mathfrak{g}_0\right) \to \left(H_1,\mathfrak{g}_1\right)$ and
\item
a holomorphic $\left(H_0,\mathfrak{g}_0\right)$-geometry $E_0 \to M_0$ and
\item
a complete holomorphic $\left(H_1,\mathfrak{g}_1\right)$-geometry $E_1 \to M_1$ and
\item
a $\Phi$-morphism $F : E_0 \to E_1$.
\end{enumerate}
Suppose that $H_1$ is Stein and $M_0$ is a domain in a Stein manifold  with envelope of holomorphy $\hat{M}_0$.

Then the holomorphic principal $H_0$-bundle $E_0 \to M_0$ extends to a holomorphic
principal $H_0$-bundle $E'_0 \to \hat{M}_0$ if and only if
the Cartan geometry on $E_0$ extends to a 
holomorphic Cartan geometry $E'_0 \to \hat{M}_0$.
If either extension exists, then both are unique and
$F$ extends to a unique $\Phi$-morphism $F : E'_0 \to E_1$.
\end{proposition}
\begin{proof}
The Cartan geometry $E_0 \to M_0$ extends to a Cartan geometry $E'_0 \to \hat{M}_0$ just when the
holomorphic principal $H_0$-bundle $E_0 \to M_0$ extends to a holomorphic principal $H_0$-bundle $E'_0 \to \hat{M}_0$;
see McKay \cite{McKay:2009} p. 19 theorem 8. By theorem~\vref{theorem:CoframingExtension}, the morphism $F$ extends to $F : E'_0 \to E_1$.
\end{proof}

\section{Morphism deformations}

\begin{definition}\label{definition:family}
A \emph{family of Cartan geometries} modelled on a local model $\left(H_0,\mathfrak{g}_0\right)$ is a choice of 
\begin{enumerate}
 \item 
a principal $H_0$-bundle $\pi : E \to M$ and
\item
a foliation of $M$, i.e. a bracket-closed vector subbundle $VM \subset TM$
(and we denote by $VE$ the preimage of $VM$ by $\pi'$),
\item
a section $\omega_0$ of the vector bundle $V^*E \otimes \mathfrak{g}_0 \to E$ so that 
\item
for each leaf $M_0 \subset M$ of $VM$, the bundle 
$\pi^{-1} M_0 \to M_0$ is a Cartan geometry with $\omega_0$ as Cartan connection.
\end{enumerate}
\end{definition}

\begin{definition}
If the foliation $VM$ of a family of Cartan geometries is the kernel of a smooth submersion to a manifold, say $S$, we will
call $S$ the \emph{parameter space} of the family.
\end{definition}

\begin{remark}
It is often convenient to assume that the parameter space has a chosen point, say $s_0 \in S$, so that
we imagine that the deformation is deforming a specific Cartan geometry $E_{p_0} \to M_{p_0}$.
\end{remark}

\begin{definition}\label{definition:morphismDeformation}
Pick a local model morphism $\Phi : \left(H_0,\mathfrak{g}_0\right) \to \left(H_1,\mathfrak{g}_1\right)$ 
A $\Phi$-\emph{morphism deformation} is a choice of
\begin{enumerate}
\item
deformation of Cartan geometries $E \to M$ modelled on $\left(H_0,\mathfrak{g}_0\right)$, with notation as in
definition~\vref{definition:family}, and
\item
a Cartan geometry $E_1 \to M_1$ modelled on $\left(H_1,\mathfrak{g}_1\right)$ and
\item
an $H_0$-equivariant map $F : E \to E_1$ so that
\item
above each leaf of the foliation of $M$, $F$ is a $\Phi$-morphism of Cartan geometries. 
\end{enumerate}
\end{definition}

\begin{definition}
Suppose that $F_0 : E_0 \to E_1$ is a morphism of Cartan geometries.
A morphism deformation $F : E \to E_1$ is a morphism deformation
of $F_0$ if there is a map $g : E_0 \to E$ which immerses 
$E_0$ as a maximal integral manifold of $VE$, and so that 
$F \circ f = F_0$.
\end{definition}

\section{Universal morphism deformations}

\begin{definition}
A deformation $F : E \to E_1$ of a morphism $E_0 \to E_1$
is \emph{versal} along $E_0$ if for every deformation $F' : E' \to E_1$ 
of $E_0$, there is an $I$-morphism $G : E' \to E$ so that 
$F' = F \circ G$. 
\end{definition}

\begin{definition}
A versal deformation is \emph{universal} if all of its induced
$I$-morphisms to all other versal deformations are
immersions.
\end{definition}

\begin{lemma}\label{lemma:Universal}
Suppose that
$\Phi : \left(H_0,\mathfrak{g}_0\right) \to \left(H_1,\mathfrak{g}_1\right)$ 
is a immersive local model morphism, i.e. an immersion
on $H_0$ and on $\mathfrak{g}_0$.
Let $q : \mathfrak{g}_1 \to \mathfrak{g}_1/\Phi \mathfrak{g}_0$
be the obvious linear quotient map.
Suppose that $E_1 \to M_1$ is a $\left(H_1,\mathfrak{g}_1\right)$-geometry.
Suppose that the $\Phi$-obstruction of $E_1$ vanishes.
Then there is a unique universal $\Phi$-deformation
$F : E \to E_1$ from a principal $H_0$-bundle $\pi : E \to M$ 
up to isomorphism. The universal $\Phi$-deformation
has $E=E_1$, $F=\text{Id}$, $M=E/H_0$, 
$VE = \left(q\omega_1=0\right)$, $\omega=\omega_1$,
and $VM=\pi' VE$.
\end{lemma}
\begin{proof}
First let's see that $\pi : E \to M$ is a morphism
deformation, with morphism $F=\text{Id}$.
Because $E_1 \to M_1$ is a principal right $H_1$-bundle,
and $H_0 \subset H_1$ is a subgroup, $H_0$ acts
freely and properly on $E=E_1$. Therefore the 
quotient $M=E/H_0$ is a manifold, and the projection 
$\pi : E \to M$ is a principal right $H_0$-bundle.
(Indeed $\pi : E \to M$ is just the lift of the $H_1$-geometry
to an $H_0$-geometry.) The vector bundle $VE$ is 
an $H_0$-invariant subbundle of $TE$, and bracket
closed, hence consists entirely of its own
Cauchy characterstics, so descends to a vector subbundle $VM \subset TM$
and $VE = \left(\pi'\right)^{-1} VM$.
Moreover, $VM$ is bracket closed because $VE$ is.
The rest is clear.

Next, let's see why this morphism deformation is
versal. Pick $F' : E' \to E_1$ any
morphism deformation. Let $G = F' : E' \to E_1$,
viewed as an $I$-morphism deformation. Clearly
we find $F' = F \circ G$.

Next, let's see why this morphism deformation is
universal. Suppose that $F' : E' \to E_1$ 
is another versal $\Phi$-morphism deformation. Then
there is an $I$-morphism $G : E_1 \to E'$ so that
$F=F' \circ G$. But $F=\text{Id}$, so $\text{Id}=F' \circ G$.
Therefore $G$ is an immersion. 

Uniqueness of the universal morphism deformation is clear.
\end{proof}

\begin{example}\label{example:ProjectiveConnectionGeodesics}
 Suppose that $G_1=\PSL{n+1,\R{}}$ and $H_1$ is the subgroup fixing a point $p_0 \in \RP{n}$.
Let $G_0 \subset G_1$ be the subgroup fixing a projective line in $\RP{n}$
containing $p_0$, and let $H_0=G_0 \cap H_1$. Let $\Phi : G_0 \to G_1$ be the obvious
inclusion mapping. Then our lemma says that every projective connection (i.e. $G_1/H_1$-geometry)
has a universal space of pointed geodesics $M$, and a universal deformation of geodesics.
It has a smooth parameter space of geodesics if and only if the $\Phi$-obstruction Pfaffian
system consists in the leaves of a continuous map.
\end{example}

\begin{remark}
It is not clear how to build a universal deformation of a nonimmersive morphism.
\end{remark}

\begin{example}\label{example:ScalarODE}
 It is well known that every scalar ordinary differential equation
of some order, say order $n+1$, induces a Cartan geometry
on its configuration space. To be more precise, 
write $\Sym{n}{\R{2}}^*$ for the space of homogeneous
polynomials of degree $n$ in two variables $x$ and $y$, write
$C_n \subset \GL{2,\R{}}$ for the group $\pm I$ if $n$ is even and $I$ if
$n$ is odd, and let
$G_1=\left(\GL{2,\R{}}/C_n\right) \rtimes \Sym{n}{\R{2}}^*$.

Let $\OO{n}$ be the set of all pairs $(L,q)$
so that $L \subset \R{2}$ is a line through the origin
and $q \in \Sym{n}{L}^*$ is a homogeneous polynomial
on $L$ of degree $n$. 
There is an obvious map $\OO{n} \to \RP{1}$ given by
$\left(L,q\right) \mapsto L$. The manifold $\OO{n}$
is the total space of the usual line bundle also
called $\OO{n}$ over $\RP{1}$.

Let $G_1$ act on $\OO{n}$ by
\[
\left(g,p\right)\left(L,q\right)=
\left(gL,q \circ g^{-1} + \left.p\right|_L\right).
\]
Let $H_1$ be the subgroup preserving the point
$\left(L_0,0\right)$ where $L_0$ is the
horizontal axis in $\R{2}$.

It is well known (see Lagrange \cite{Lagrange:1957a,Lagrange:1957b}, Fels \cite{Fels:1993,Fels:1995}, 
Dunajski and Tod \cite{Dunajski/Tod:2006}, Godli{\'n}ski and Nurowski \cite{Godlinski/Nurowski:2007}, 
Doubrov \cite{Doubrov:2008}) that every scalar ordinary differential equation 
of order $n+1$, 
\[
\frac{d^{n+1} y}{dx^{n+1}}
=
f\left(x,y,\frac{dy}{dx},\frac{d^2 y}{dx^2}, \dots, \frac{d^n y}{dx^n} \right)
\]
imposes a $G_1/H_1$-geometry on the 
surface parameterized by the variables $x,y$,
and that this $G_1/H_1$-geometry is invariant
under fiber-preserving transformations of
the ordinary differential equation.
Therefore we can view any $G_1/H_1$-geometry
on a surface as a natural geometric generalization
of an ordinary differential equation
of order $n+1$.

Let $G_0=\GL{2,\R{}}/C_n$ and $H_0 = G_0 \cap H_1$.
So $G_0/H_0=\RP{1}$. The orbits of $G_0$ in $\OO{n}$ are precisely
the graphs of the global sections of $\OO{n}$
given by taking a homogeneous polynomial $Q(x,y)$ of degree $n$
on $\R{2}$, and mapping $L \in \RP{1} \mapsto \left(L,\left.Q\right|_{L}\right) \in \OO{n}$.
In particular, we have an obvious morphism $\Phi : G_0/H_0 \to G_1/H_1$
of homogeneous spaces: the zero section $\RP{1} \to \OO{n}$ of the bundle map 
$\OO{n} \to \RP{1}$.

If we have a $G_1/H_1$-geometry $E \to M$ which is constructed out of an ordinary
differential equation of degree $n+1$, then the $\Phi$-morphisms of open sets of $\RP{1}$ 
are precisely the local solutions, equipped with a natural projective connection.
The universal morphism deformation is just the family of all solutions of 
the original scalar ordinary differential equation. 
Therefore in all scalar ODE systems of degree at least 2 are examples of universal morphism deformations
of Cartan geometries. 
\end{example}

\begin{definition}
Suppose that
$\Phi : \left(H_0,\mathfrak{g}_0\right) \to \left(H_1,\mathfrak{g}_1\right)$ 
is a immersive local model morphism, i.e. an immersion
on $H_0$ and on $\mathfrak{g}_0$.
Let $q : \mathfrak{g}_1 \to \mathfrak{g}_1/\Phi \mathfrak{g}_0$
be the obvious linear quotient map.
Suppose that $E_1 \to M_1$ is a $G_1/H_1$-geometry.

We will say that the $G_1/H_1$-geometry $E_1 \to M_1$
is \emph{$\Phi$-tame} (or just \emph{tame} if $\Phi$ is understood) if
\begin{enumerate}
\item
the $\Phi$-obstruction of $E_1$ vanishes (so that in particular $E_1 \to M_1$
has a universal morphism deformation)
 \item 
the universal morphism deformation $E \to M$ of $E_1 \to M_1$ admits a parameter space $S$,
and
\item the map $E \to S$ is a smooth fiber bundle morphism with compact fibers.
\end{enumerate}
\end{definition}

\section{The equation of first order deformation}

\begin{definition}
Pick a local model monomorphism $\Phi : \left(H_0,\mathfrak{g}_0\right) \to \left(H_1,\mathfrak{g}_1\right)$.
Let $q : \mathfrak{g}_1 \to \mathfrak{g}_1/\Phi \mathfrak{g}_0$ be the obvious
linear projection map.

Suppose that $F : E \to E_1$ is a $\Phi$-morphism deformation, with bundles 
$H_0 \to E \to M$ and $H_1 \to E_1 \to M_1$. 
Clearly $q F^* \omega_1$, a section of the vector bundle 
$\left(TE/VE\right)^* \otimes \left(\mathfrak{g}_1/\Phi \mathfrak{g}_0\right)$.
Since $q F^* \omega_1$ is clearly $H_0$-equivariant, it is also
a section of the vector bundle 
$Q_E = \left(TE/VE\right)^* \otimes_{H_0} \left(\mathfrak{g}_1/\Phi \mathfrak{g}_0\right) \to M$;
we will write $\delta F$ to mean $q F^* \omega_1$ viewed as a section 
of the vector bundle $Q_E \to M$.
This object $\delta F$ will be called the \emph{first order deformation} of $F$.
The restriction of $\delta F$ to a leaf of $VM$ will be called
the first order deformation of that leaf.
\end{definition}

\begin{definition}
In the notation of the previous definition,
for any vector $A \in \mathfrak{g}_1$, 
write $\bar{A}$ for $A+\mathfrak{h}_1 \in \mathfrak{g}_1/\mathfrak{h}_1$.
\end{definition}

\begin{definition}
 Suppose that $\rho : \mathfrak{g} \to \gl{V}$ is a representation, and that 
$E \to M$ is a $\left(H,\mathfrak{g}\right)$-geometry or a family of $\left(H,\mathfrak{g}\right)$-geometries. The
\emph{covariant derivative} $\nabla$ is the operator on sections of $E \times_H V$
given by
\[
\nabla s = ds + \rho(\omega)s.
\]
\end{definition}

\begin{lemma}
Consider a local model morphism 
$\Phi : \left(H_0,\mathfrak{g}_0\right) \to \left(H_1,\mathfrak{g}_1\right)$.
Consider a morphism deformation, in the notation of definition~\ref{definition:morphismDeformation},
modelled on $\Phi$. Suppose that the parameter space $S$ is an interval $S \subset \R{}$,
and let $t : S \to \R{}$ be the inclusion mapping.
The first order deformation is $\delta F = a \, dt$
where $a : E \to \mathfrak{g}_1/ \Phi \mathfrak{g}_0$
is a uniquely determined $H_0$-equivariant function,
i.e. a section of $E \times_{H_0} \left(\mathfrak{g}_1/ \Phi \mathfrak{g}_0\right)$.
Let $\rho : \mathfrak{g}_0 \to \gl{\mathfrak{g}_1/\Phi \mathfrak{g}_0}$
be the obvious representation of the Lie algebra $\mathfrak{g}_0$.
Let $Q = \mathfrak{g}_1/ \left( \Phi \mathfrak{g}_0 + \mathfrak{h}_1 \right)$.
Let $R : E \to \left(\mathfrak{g}_0/\mathfrak{h}_0\right)^* \otimes Q^* \otimes \left(\mathfrak{g}_1/\Phi \mathfrak{g}_0\right)$
be the map 
\[
R(e)(A,B) = q K_1\left(F\left(e\right)\right) \left( \overline{\Phi A} \wedge \bar{C} \right),
\]
where $qC=B$. This function $R$ is well defined, and called the
\emph{morphism curvature} or \emph{$\Phi$-curvature}. The function $a$ restricted to a leaf $E_0 \subset E$ 
of the morphism deformation satisfies the \emph{equation of first variation}:
\[
 \nabla a = R\left(\overline{\Phi \omega_0}, \bar{a}\right).
\]
\end{lemma}
\begin{remark}
The morphism curvature is analogous to the Riemann curvature tensor appearing in
the Jacobi vector field equation of a geodesic in a pseudo-Riemannian manifold;
see example~\vref{example:ProjectiveConnectionCurvature}.
\end{remark}
\begin{proof}
The first order deformation is $\delta F = q F^* \omega_1$.
It vanishes on $VE$, a corank 1 subbundle of $TE$. Clearly
$VE=\left(dt=0\right)$, so $\delta F = a \, dt$ for a 
unique function $a$, and clearly $a : E \to \mathfrak{g}_1/ \Phi \mathfrak{g}_0$.
The $H_0$-equivariance of $a$ follows immediately from that of $\omega_1$.
\begin{align*}
d \left(a \, dt\right) 
&=
da \wedge dt \\
&=
d \delta F 
\\
&=
d q F^* \omega_1
\\
&=
q F^* \left( -\frac{1}{2}\left[\omega_1,\omega_1\right] + \frac{1}{2}K_1 \bar{\omega}_1 \wedge \bar{\omega}_1 \right).
\end{align*}
Each vector tangent to $E_0$ has the form
$\vec{A}$ for some unique $A \in \mathfrak{g}_0$.
\begin{align*}
 \vec{A} \hook \left( da \wedge dt\right)
&=
\left(\vec{A} \hook da\right) dt 
\\
&=
q \left( - \Ad(A) F^* \omega_1 + K_1\left(\overline{\Phi A},F^* \bar{\omega}_1\right) \right)
\\
&=
- \rho(A) q \omega_1 + q K_1 \left(\overline{\Phi A},F^* \bar{\omega}_1\right)
\\
&=
- \rho(A) a \, dt + q K_1 \left(\overline{\Phi A}, F^* \bar{\omega}_1\right).
\end{align*}
If we pick some vector tangent to $E_0$, say $\vec{B}$ for some
$B \in \mathfrak{g}_0$, then
\begin{align*}
0&=
 \vec{B} \hook \left( \left( \overline{\Phi A} \hook da\right) dt \right)
\\
&=
\vec{B} \hook \left(
- \rho(A) a \, dt + q K_1 \left(\overline{\Phi A}, F^* \bar{\omega}_1\right).
\right)
\\
&=
q K_1 \left(\overline{\Phi A},\overline{\Phi B}\right)
\\
&=
q \Phi K_0 \left(\bar{A}, \bar{B}\right)
\\
&=
0.
\end{align*}
Therefore the expression
\[
q K_1 \left(\overline{\Phi A},F^* \bar{\omega}_1\right)
\]
vanishes on vectors tangent to $E_0$, and so is
defined on vectors in $\left.TE\right|_{E_0}/TE_0$,
and is expressible in terms of $q F^*\omega_1$, i.e.
in terms of $a \, dt$, as
\[
q K_1 \left(\overline{\Phi A},F^* \bar{\omega}_1\right)
=R\left(\overline{\Phi A}, \bar{a}\right) \, dt.
\]
Hence $R$ is well defined, and the equation of first variation follows.
\end{proof}

Next, we would like a notion of infinitesimal deformation
which does not depend on the existence of any actual
deformation.

\begin{definition}
Consider a local model morphism 
$\Phi : \left(H_0,\mathfrak{g}_0\right) \to \left(H_1,\mathfrak{g}_1\right)$.
Suppose that $E_0 \to M_0$ and $E_1 \to M_1$ are Cartan geometries
and that $F : E_0 \to E_1$ is a $\Phi$-morphism.
Let $\rho : \mathfrak{g}_0 \to \gl{\mathfrak{g}_1/\Phi \mathfrak{g}_0}$
be the obvious representation of the Lie algebra $\mathfrak{g}_0$.
Let $Q = \mathfrak{g}_1/ \left( \Phi \mathfrak{g}_0 + \mathfrak{h}_1 \right)$.
Let  $R : E_0 \to \left(\mathfrak{g}_0/\mathfrak{h}_0\right)^* \otimes Q^* \otimes \left(\mathfrak{g}_1/\Phi \mathfrak{g}_0\right)$
be the map 
\[
R(e)(A,B) = q K_1\left(F\left(e\right)\right) \left( \overline{\Phi A} \wedge \bar{C} \right),
\]
where $qC=B$. 
Again we will refer to $R$ as the \emph{$\Phi$-curvature}, or \emph{morphism curvature}.

An \emph{infinitesimal one-parameter morphism deformation} of $F$
is a section $a$ of $E \times_{H_0} \left(\mathfrak{g}_1/ \Phi \mathfrak{g}_0\right)$
satisfying the equation of first variation
\[
 \nabla a = R\left(\overline{\Phi \omega_0}, \bar{a}\right).
\]
\end{definition}

\begin{remark}
 Since the equation of first variation is a total differential
equation, each solution is completely determined by its value
at any chosen point. In particular, the dimension of the
space of global solutions is at most the dimension of
$\mathfrak{g}_1/ \Phi \mathfrak{g}_0$. Formally we can
think of the space of infinitesimal morphism deformations 
as the tangent space to the space of morphism deformations.
\end{remark}

\begin{example}
Consider a morphism of homogeneous spaces $\Phi : G_0/H_0 \to G_1/H_1$.
The infinitesimal deformations of the model are just the global solutions of $\nabla a=0$.
We can take any element $a_0 \in \mathfrak{g}_1/\Phi \mathfrak{g}_0$,
and then let $a\left(g_0\right) = \rho\left(g_0\right)^{-1}a_0$.
It is easy to see that every infinitesimal morphism
deformation has this form, so the space
of infinitesimal morphism deformations of the
model $G_0/H_0 \to G_1/H_1$ is precisely $\mathfrak{g}_1/\Phi \mathfrak{g}_0$.
\end{example}

\begin{remark}\label{remark:SectionOfNormalBundle}
To each infinitesimal deformation $a$ of an immersive morphism $F : E_0 \to E_1$
of Cartan geometries, there is associated a section
of the normal bundle $F^* TE_1 /TE_0$, given by using
the Cartan connection $\omega_1$ on $E_1$ to identify
the normal bundle with $E_0 \times_{H_0} \left(\mathfrak{g}_1/ \Phi \mathfrak{g}_0\right)$.
\end{remark}

\begin{lemma}
 Suppose that $E \to M \to S$ is a family of Cartan geometries,
say with maps $\pi : E \to M$ and $p : M \to S$.
For each point $s \in S$, let $M_s = p^{-1} s$ and $E_s = \pi^{-1} p^{-1} s$.
Suppose that $F : E \to E_1$ is an immersive morphism deformation. Then 
at each point $s \in S$, to each tangent vector $\dot{s} \in T_s S$,
we can associate a unique infinitesimal deformation 
$a=a_{\dot{s}} : E_s \to \mathfrak{g}_1/ \Phi \mathfrak{g}_0$
so that, identifying $a$ with a section $A$ of the
normal bundle $F^* TE_1 /TE_s$ (as in remark~\vref{remark:SectionOfNormalBundle}), and taking any path
$e(t) \in E$ for which $\left(p \circ \pi \circ e\right)'(0)=\dot{s}$,
we have 
\[
A\left(e(0)\right)=F'\left(e(0)\right)e'(0) + F'\left(e(0)\right) V_{e(0)}E.
\]
We call this $A=A_{\dot{s}}$ the \emph{associated section of the normal bundle}.
\end{lemma}
\begin{proof}
We have to show that we can construct an infinitesimal
deformation $a$ with
\[
a\left(e(0)\right)=F'\left(e(0)\right)e'(0) \hook q \omega_1.
\]
Take any path $s(t)$ in $S$ with $s'(0)=\dot{s}$,
and replace $S$ by an interval of this path,
and $E$ and $M$ by their respective pullbacks
via the map $t \mapsto s(t)$. Then we have
to prove that there is some $a$ so that
\[
 a \, dt = q F^* \omega_1.
\]
But clearly since $qF^* \omega_1=a \, dt$ 
for a unique infinitesimal deformation $a$ 
as seen above, this choice of $a$ is precisely
the one we must make.

Now consider returning to a general parameter
space $S$. We can pick out a choice of $a$ as above,
depending on a choice of path $s(t)$.
We need to show that this choice will satisfy
\[
a\left(e(0)\right)=F'\left(e(0)\right)e'(0) \hook q \omega_1
\]
for any path $e(t)$ in $E$ for which $\left(p \circ \phi \circ e\right)'(0)=\dot{s}$.
This certainly works as long as $e(t)$ is a path
inside $E$ lying over the curve $s(t)$.

So now we need only show that this
infinitesimal deformation $a$ depends only
on $\dot{s}$, not on the particular
path $s(t)$. So suppose that we have
two paths $s_1(t)$ and $s_2(t)$ in $S$,
with $s_1'(0)=s_2'(0)=\dot{s}$.
We then make infinitesimal deformations 
$a_1$ and $a_2$ so that if $e_1(t)$
is a path in the pullback of $E$
over $s_1(t)$, then
\[
a_1\left(e_1(0)\right)=F'\left(e_1(0)\right)e_1'(0) \hook q \omega_1,
\]
and similarly for $e_2(t)$. Assume that $e_1(0)=e_2(0)=e_0$
and that $\left(p \circ \pi \circ e_1\right)'(0)=\left(p \circ \pi \circ e_2\right)'(0)=\dot{s}$.
So $e_1'(0)-e_2'(0) \in V_{e_0} E$. So
\begin{align*}
a_1\left(e_1(0)\right)-a_2\left(e_2(0)\right)
&=
\left(e_1'(0)-e_2'(0)\right) \hook qF^* \omega_1(e_0)
\\
&=
0.
\end{align*}
\end{proof}

\begin{lemma}
Suppose that $\Phi : \left(H_0, \mathfrak{g}_0\right) \to \left(H_1, \mathfrak{g}_1\right)$ is an immersive 
local model morphism.
Suppose that $E_1 \to M_1$ is a $\left(H_1,\mathfrak{g}_1\right)$-geometry with vanishing $\Phi$-obstruction.
Suppose that $E \to E_1$ is the universal $\Phi$-morphism 
deformation, and has a parameter space $S$. The map $\dot{s} \mapsto A_{\dot{s}}$ 
taking a tangent vector on the parameter space to its
associated section of the normal bundle is a linear
isomorphism from tangent vectors to $TS$ to infinitesimal
morphism deformations.
\end{lemma}
\begin{proof}
Let's start by finding the kernel. If $A_{\dot{s}}=0$,
then taking any path
$e(t) \in E$ for which $\left(p \circ \pi \circ e\right)'(0)=\dot{s}$,
we have 
\[
0=A\left(e(0)\right)=F'\left(e(0)\right)e'(0) + F'\left(e(0)\right) V_{e(0)}E.
\]
So $F'$ is not an immersion. But $F : E \to E_1$ must
be an immersion by lemma~\vref{lemma:Universal}.

The universal deformation has the same dimension
as the model, and its parameter space is therefore
the same dimension as $\mathfrak{g}_1/\Phi \mathfrak{g}_0$,
and so the linear map $\dot{s} \mapsto A_{\dot{s}}$
is a linear isomorphism.
\end{proof}

\begin{example}\label{example:ProjectiveConnectionCurvature}
We continue example~\vref{example:ProjectiveConnectionGeodesics}: geodesics
of projective connections. We can write any element
$A$ in the various Lie algebras in the manner indicated:
\[
\begin{array}{ccc}
\text{Lie algebra} & \text{Also known as} & \text{Typical element} \\
\mathfrak{g}_1 
& 
\slLie{n+1,\R{}}
& 
\begin{pmatrix}
 A^0_0 & A^0_1 & A^0_J \\
A^1_0 & A^1_1 & A^1_J \\
A^I_0 & A^I_1 & A^I_J
\end{pmatrix}
\\
\mathfrak{h}_1
& 
\text{stabilizer of a point}
&
\begin{pmatrix}
 A^0_0 & A^0_1 & A^0_J \\
0 & A^1_1 & A^1_J \\
0 & A^I_1 & A^I_J
\end{pmatrix}
\\
\mathfrak{g}_0
&
\text{stabilizer of a line}
&
\begin{pmatrix}
 A^0_0 & A^0_1 & A^0_J \\
A^1_0 & A^1_1 & A^1_J \\
0 & 0 & A^I_J
\end{pmatrix}
\\
\mathfrak{h}_0 
& 
\text{stabilizer of a pointed line}
&
\begin{pmatrix}
 A^0_0 & A^0_1 & A^0_J \\
0 & A^1_1 & A^1_J \\
0 & 0 & A^I_J
\end{pmatrix}
\end{array}
\]
with indices $I,J = 2, \dots, n$.

The Cartan geometry of a projective connection looks like
\[
\omega_1 =
\begin{pmatrix}
 \omega^0_0 & \omega^0_1 & \omega^0_J \\
 \omega^1_0 & \omega^1_1 & \omega^1_J \\
 \omega^I_0 & \omega^I_1 & \omega^I_J \\
\end{pmatrix}.
\]
The quotient $\bar{\omega}_1 = \omega_1 + \mathfrak{h}_1$
looks like
\[
\bar{\omega}_1 
=
\begin{pmatrix}
\omega^1_0 \\
\omega^I_0
\end{pmatrix}.
\]
The curvature of a projective connection looks like
\[
K_1 \bar{\omega}_1 \wedge \bar{\omega}_1= 
\begin{pmatrix}
 \nabla \omega^0_0 
 & 
  \nabla \omega^0_1 
 & 
  \nabla \omega^0_J \\
  \nabla \omega^1_0 
 & 
  \nabla \omega^1_1 
 &  \nabla \omega^1_J \\
  \nabla \omega^I_0 
 & 
  \nabla \omega^I_1 
 & 
  \nabla \omega^I_J \\
\end{pmatrix}
\]
where
\[
\nabla \omega^{\mu}_{\nu} = 
 2 K^{\mu}_{\nu 1I} \omega^1_0 \wedge \omega^I_0  
 +
 K^{\mu}_{\nu IJ} \omega^I_0 \wedge \omega^J_0
\]
for $\mu, \nu=0,1,2,\dots,n$.
The first order variation $\delta F$ 
of the universal deformation is 
\begin{align*}
\delta F
&=
q K_1 \bar{\omega}_1 \wedge \bar{\omega}_1
\\
&=
\begin{pmatrix}
\nabla \omega^I_0 & \nabla \omega^I_1 \\
\end{pmatrix}
\\
&=
\begin{pmatrix}
 2 K^{I}_{0 1I} \omega^1_0 \wedge \omega^I_0  
 +
 K^{I}_{0 IJ} \omega^I_0 \wedge \omega^J_0
& 
 2 K^{I}_{1 1I} \omega^1_0 \wedge \omega^I_0  
 +
 K^{I}_{1 IJ} \omega^I_0 \wedge \omega^J_0
\end{pmatrix}.
\end{align*}

The morphism curvature of the universal deformation is
\[
R\left(A^1_0,B^I_0\right)
=
\begin{pmatrix}
K^I_{01J} A^1_0 B^J_0
&
K^I_{11J} A^1_0 B^J_0
\end{pmatrix}.
\]
An infinitesimal deformation along one of the morphisms
(i.e. leaves) of the universal deformation has the form
\[
a = 
\begin{pmatrix}
 a^I_0 & a^I_1
\end{pmatrix}.
\]
The equation of first variation is
\[
\nabla
\begin{pmatrix}
a^I_0 & a^I_1
\end{pmatrix}
=
\begin{pmatrix}
K^I_{01J} \omega^1_0 a^J_0
&
K^I_{11J} \omega^1_0 a^J_0
\end{pmatrix}.
\]

The representation $\rho$ is
\[
\rho(A)
\begin{pmatrix}
 B^I_0 & B^I_1
\end{pmatrix}
=
\begin{pmatrix}
A^I_J B^J_0 - B^I_0 A^0_0 - B^I_1 A^1_0 
&
A^I_J B^J_1 - B^I_0 A^0_1 - B^I_1 A^1_1
\end{pmatrix}.
\]
Therefore the covariant derivative $\nabla a$
of any infinitesimal deformation $a$ is
\[
\nabla 
\begin{pmatrix}
a^I_0 & a^I_1
\end{pmatrix}
=
\begin{pmatrix}
\omega^I_J a^J_0 - a^I_0 \omega^0_0 - a^I_1 \omega^1_0 
&
\omega^I_J a^J_1 - a^I_0 \omega^0_1 - a^I_1 \omega^1_1
\end{pmatrix}.
\]
Finally, we obtain the equation of first variation
of the universal deformation:
\begin{align*}
da^I_0 &= a^I_0 \omega^0_0 + a^I_1 \omega^1_0 - \omega^I_J a^J_0 + K^I_{01J} \omega^1_0 a^J_0 \\
da^I_1 &= a^I_0 \omega^0_1 + a^I_1 \omega^1_1 - \omega^I_J a^J_1 + K^I_{11J} \omega^1_0 a^J_0.
\end{align*}

The projective connection is torsion-free if and only if $K^I_{01J}=0$. Therefore
a torsion-free projective connection has equation of first variation
of its universal deformation:
\begin{align*}
da^I_0 &= a^I_0 \omega^0_0 + a^I_1 \omega^1_0 - \omega^I_J a^J_0 \\
da^I_1 &= a^I_0 \omega^0_1 + a^I_1 \omega^1_1 - \omega^I_J a^J_1 + K^I_{11J} \omega^1_0 a^J_0,
\end{align*}
which we can see is naturally a 2nd order system of ordinary differential
equations along each geodesic: precisely the well known Jacobi vector field equations.
\end{example}

\begin{example}
We leave the reader to generalize the previous example to totally geodesic projectively flat subspaces
inside a manifold with projective connection. This is really just a slight change of notation from
the previous example; see Cartan \cite{Cartan:1992} \emph{Le\,cons sur la
th\'eorie des espaces \'a connexion projective}, chapter VI.
\end{example}

\begin{example}\label{example:WilczynskiInvariants}
Let's continue with example~\vref{example:ScalarODE}: $G_1/H_1=\OO{n}$,
$G_1=\left(\GL{2,\R{}}/C_n\right) \rtimes \Sym{n}{\R{2}}^*$.
It is convenient to write elements $(g,p) \in G_1$ formally as matrices
\[
\begin{pmatrix}
g & p \\
0 & 1
\end{pmatrix}
\]
to make the semidirect product structure more apparent. In this setting, we write $gp$ to mean
$p \circ g^{-1}$, as a homogeneous polynomial. Clearly
\[
\begin{pmatrix}
g & p \\
0 & 1
\end{pmatrix}^{-1}
=
\begin{pmatrix}
g^{-1} & -g^{-1}p \\
0 & 1
\end{pmatrix}.
\]
We will write any homogeneous polynomial $p \in \Sym{n}{\C{2}}^*$ 
as 
\[
p\left(x,y\right)=\sum_{i+j=n} p_{ij} x^i y^j. 
\]

Suppose that 
\[
\xymatrix{%
H_1 \ar[r] & E_1 \ar[d]^{\pi} \\
         & M_1,
}%
\]
is an $\OO{n}$-geometry on a surface $M_1$, with Cartan 
connection $\omega_1 \in \nForms{1}{E_1} \otimes \mathfrak{g}$.
We can write
\[
\omega=
\begin{pmatrix}
 \gamma & \varpi \\
0 & 0
\end{pmatrix},
\]
with $\gamma \in \nForms{1}{E} \otimes \gl{2,\C{}}$ and
$\varpi \in \nForms{1}{E} \otimes \Sym{n}{\C{2}}^*$,
say
\[
\gamma = 
\begin{pmatrix}
 \gamma^1_1 & \gamma^1_2 \\
\gamma^2_1 & \gamma^2_2
\end{pmatrix}
\]
and
\[
 \varpi = \left( \varpi_{ij} \right),
\]
where we write each homogeneous polynomial $p\left(x,y\right)$
as
\[
p\left(x,y\right)=\sum_{i+j=n} p_{ij} x^i y^j.
\]

Under $H_1$-action on $E_1$, $\omega_1$ varies in the adjoint representation,
$r_h^* \omega=\Ad_h^{-1} \omega_1$. Therefore, if we write
\[
h=
\begin{pmatrix}
 g & p \\
0 & 1
\end{pmatrix}
\]
then
\[
r_h^* 
\begin{pmatrix}
\gamma & \varpi \\
0 & 0
\end{pmatrix}
=
\begin{pmatrix}
 g^{-1} \gamma g &  g^{-1} \gamma p + g^{-1} \varpi \\
0 & 0
\end{pmatrix}.
\]
The action on $\varpi$, expanded out, when we think of $\varpi$
as a 1-form valued in symmetric $n$-linear forms, is
\begin{align*}
\left(g^{-1} \gamma p + g^{-1} \varpi\right)\left(v_1, v_2, \dots, v_n\right)
=&
-
\sum_k
p\left(v_1, v_2, \dots, v_{k-1},
g^{-1} \gamma v_k, v_{k+1}, \dots, v_n\right)
\\
&+
\varpi\left(gv_1, gv_2, \dots, gv_n\right).
\end{align*}

The subalgebra $\mathfrak{h}_1 \subset \mathfrak{g}_1$ consists in the elements for
which $\gamma^2_1=0$ and $\varpi_{n0}=0$. Therefore we can treat the pair
\[
\left(\gamma^2_1, \varpi_{n0} \right)
\]
as being a 1-form valued in $\mathfrak{g}_1/\mathfrak{h}_1$. In particular,
we can write our structure equations as
\[
d 
\begin{pmatrix}
\gamma & \varpi \\
0 & 0 
\end{pmatrix}
+
\begin{pmatrix}
\gamma & \varpi \\
0 & 0 
\end{pmatrix}
\wedge 
\begin{pmatrix}
\gamma & \varpi \\
0 & 0 
\end{pmatrix}
=
\begin{pmatrix}
K & L \\
0 & 0 
\end{pmatrix}
\gamma^2_1 \wedge \varpi_{n0},
\]
where 
$K : E_1 \to \Lm{2}{\mathfrak{g}_1/\mathfrak{h}_1}^* \otimes \gl{2,\C{}}$ 
and 
$W : E_1 \to \Lm{2}{\mathfrak{g}_1/\mathfrak{h}_1}^* \otimes \Sym{n}{\C{2}}^*$.
The invariants $W=\left(W_{ij}\right)$ are called the
\emph{Wilczynski invariants}. Expanding this out gives
\begin{align*}
d \gamma + \gamma \wedge \gamma &= K \, \gamma^2_1 \wedge \varpi_{n0}, \\
d \varpi + \gamma \wedge \varpi &= W \, \gamma^2_1 \wedge \varpi_{n0}. 
\end{align*}
It is convenient to write $\varpi = \left(\varpi_{ij}\right)$ corresponding
to the expression for polynomials $p\left(x,y\right)=\sum_{i+j=n} p_{ij} x^i y^j$.
In terms of such an expression,
\[
d \varpi_{ij} = \gamma^1_1 \wedge \varpi_{ij} + \gamma^2_1 \wedge \varpi_{i-1,j+1}
+ L_{ij} \gamma^2_1 \wedge \varpi_{n0}.
\]

% We calculate that under $H_1$-action, if 
% \[
% h=
% \begin{pmatrix}
%  g & p \\
% 0 & 1
% \end{pmatrix}
% \]
% then
% \[
% r_h^* \gamma^2_1=\frac{g^1_1}{g^2_2} \gamma^2_1, 
% \]
% since $\left(g^{-1}\right)^2_1=0$. 
% 
% We immediately see that 
% \[
% d \gamma^2_1 = - \left(\gamma^2_2 - \gamma^1_1 + K^2_1 \varpi_{n0} \right) \wedge \gamma^2_1, 
% \]
% so that the Pfaffian system $\gamma^2_1=0$ satisfies the hypotheses of the Frobenius theorem.
% Moreover this Pfaffian system descends (modulo Cauchy characteristics) to 
% the surface $M_1$, giving a foliation of $M_1$ by curves. In the
% model $\OO{n}$, this foliation is the foliation by the leaves of the map $\OO{n} \to \RP{1}$.

Consider the morphism of homogeneous spaces $\Phi : G_0/H_0 \to G_1/H_1$, where
$G_0=\GL{2,\R{}}/C_n$ and $H_0=G_0 \cap G_1$.
The $\Phi$-torsion vanishes on any $G_1/H_1$-geometry, because $G_0/H_0$ is 1-dimensional.
The $\Phi$-morphism curvature is precisely $R=L$:
the morphism curvature is the collection of Wilczynski invariants.
In particular, the scalar ordinary differential equations which satisfy $L=0$ are of obvious
interest, since their first order deformations are the same as those
of the model equation $d^{n+1}y/dx^{n+1}=0$. The scalar
ordinary differential equations with $L=0$ (vanishing Wilczynski invariants) are well known and geometrically
characterized by Dunajski and Tod \cite{Dunajski/Tod:2006} and Doubrov \cite{Doubrov:2008}.
\end{example}

\section{The canonical linear system of a morphism}

\begin{definition}
 Suppose that $\left(H,\mathfrak{g}\right)$ is a local model, $V$ a $\mathfrak{g}$-module,
$V'$ an $H$-module, $\lambda : V \to V'$ a morphism of $H$-modules,
and $E \to M$ a $\left(H,\mathfrak{g}\right)$-geometry. The \emph{locus map} $f_{\lambda}$ is the
map associating to each point $m \in M$ the set
of parallel sections of $E \times_H V$ whose image
under $E \times_H V \to E \times_H V'$ vanishes
at $m$, so $f_{\lambda} : M \to \cup_{k=0}^{\dim V} \Gr{k}{Z}$
where $Z$ is the space of parallel sections of $E \times V \to M$.
\end{definition}

\begin{example}
Suppose that $G/H$ is a homogeneous space.
Let $K$ be the kernel of $\lambda : V \to V'$. The model locus map is the map taking
$gH \in G/H \mapsto gK \in \Gr{\dim K}{V}$. The model locus map is an immersion just when
the Lie algebra of the stabilizer of $K$ in $G$ is contained in the Lie
algebra of $H$.
\end{example}

\begin{definition}
Suppose that $\Phi : \left(H_0, \mathfrak{g}_0\right) \to \left(H_1, \mathfrak{g}_1\right)$
is a local model morphism. Suppose that $E_0 \to M_0$
and $E_1 \to M_1$ are Cartan geometries and $F : E_0 \to E_1$ 
is a $\Phi$-morphism. Then the locus map $\Lambda_F$ of the morphism $F$ is the locus map
of the $H_0$-module morphism $\lambda : \mathfrak{g}_1/\Phi \mathfrak{g}_0 \to \mathfrak{g}_1/\left(\mathfrak{h}_1+\Phi \mathfrak{g}_0\right)$.
\end{definition}

The locus map of a morphism associates to each point $m_0$
of the source $M_0$ the set of all infinitesimal
deformations whose underlying normal bundle section
vanishes at $m_0$.

\begin{example}\label{example:LocusMapOfModelAnImmersion}
Suppose that $\Phi : G_0/H_0 \to G_1/H_1$ is a morphism
of homogeneous spaces. Consider the locus map of $\Phi$ as a $\Phi$-morphism
of Cartan geometries. This locus map is a map $\Lambda_{\Phi} : G_0/H_0 \to \Gr{\dim K}{\mathfrak{g}_1/\Phi \mathfrak{g}_0}$
where $K=\mathfrak{h}_1+\Phi \mathfrak{g}_0 \subset \mathfrak{g}_1$.
It is defined as
\[
\Lambda_{\Phi}\left(g_0H_0\right)=g_0 K \subset \mathfrak{g}_1.
\]
Let $\mathfrak{k}$ be the set of all $A \in \mathfrak{g}_1$ for which $\Ad(A)K \subset K$.
The locus map is an immersion just when $\mathfrak{k} \subset \mathfrak{h}_1$.
\end{example}

\begin{definition}
The \emph{expected dimension} of the space of infinitesimal deformations of
a morphism of Cartan geometries is the dimension of the space of infinitesimal deformations of
the model.
\end{definition}

We can apply all of the ideas we developed so far to
prove completeness results.

\begin{proposition}
Suppose that $\Phi : G_0/H_0 \to G_1/H_1$ is a morphism of homogeneous spaces. Suppose that the model locus map
is an immersion. Suppose that $E_0 \to M_0$ is a flat $G_0/H_0$-geometry on a compact manifold $M_0$.
Suppose that $F : E_0 \to E_1$ is a $\Phi$-morphism, whose infinitesimal deformations have expected dimension,
and whose deformation curvature vanishes. Then the locus map of $F$ is a covering map to its image and
the geometry $E_0 \to M_0$ is complete.
\end{proposition}
\begin{proof}
Because $M_0$ is flat, it is locally isomorphic
to the model. Because the deformation curvature
vanishes, each local isomorphism identifies
the infinitesimal deformations and their 
vanishing loci, and therefore near each
point, the locus map is identified locally.
Because $M_0$ is compact, the image $\Lambda_{F} M_0$ of the locus map 
of $M_0$ is identical to that of $G_0/H_0$,
a $G_0$-orbit inside the appropriate Grassmannian;
let's call it $G_0/H_0'$. 
Moreover, the locus map is an immersion, because
it is an immersion of the model. In
particular, $H_0' = H_0 \rtimes \Gamma$,
for some discrete group $\Gamma$. By homogeneity,
the image $G_0/H_0'$ of either locus map 
is an embedded submanifold of the
Grassmannian. By compactness of $M_0$, the
image $G_0/H_0'$ is a compact manifold, and so the
locus map is a finite covering map $\Lambda_F : M_0 \to G_0/H_0'$.
By homogeneity, the locus map of the model
is a covering map $\Lambda_{\Phi} : G_0/H_0 \to G_0/H_0'$.
Let $\tilde{M}_0$ be the pullback:
\[
\xymatrix{
\tilde{M}_0 \ar[d] \ar[r] & G_0/H_0 \ar[d] \\
M_0 \ar[r] & G_0/H_0'.
}
\]
Clearly $\tilde{M}_0$ bears the pullback Cartan
geometry from $M_0$, and the maps from
$\tilde{M}_0$ to $M_0$ and to $G_0/H_0$
are local isomorphisms of Cartan geometries.
Indeed theorem~\vref{theorem:RealAnalyticCoveringExtension}
already furnishes us this same 
manifold $\tilde{M}_0$. 

Take any element $A \in \mathfrak{g}_0$.
Look at the vector field $X_A$ on $G_0/H_0'$
which is the infinitesimal left action
of $A$. This vector field lifts
via the locus map, since the locus map is a covering
map, to a vector field on $M_0$.
It also lifts in the obvious way
to a vector field on $G_0/H_0$.
(Keep in mind that this is \emph{not}
the constant vector field that comes from the
model Cartan connection. Instead it is 
an infinitesimal symmetry of the
flat Cartan geometry.) This vector
field on $M_0$ is locally identified
with the corresponding field on $G_0/H_0$,
via any local isomorphism given
by locally inverting the locus maps.
The flows intertwine in the obvious
way, so that the completeness of flows
on the homogeneous spaces $G_0/H_0$ and
$G_0/H_0'$ ensures the completeness
of $X_A$ on all four manifolds:
\[
\xymatrix{
\tilde{M}_0 \ar[d] \ar[r] & G_0/H_0 \ar[d] \\
M_0 \ar[r] & G_0/H_0'.
}
\]
These vector fields generate a locally transitive
action of $\mathfrak{g}_0$ on all four manifolds.
By a theorem of Ehresmann \cite{Ehresmann:1961} or McKay \cite{McKay:2004a},
each of these maps is a fiber bundle mapping, so
a covering map. Completeness of Cartan geometries 
is preserved and reflected by covering maps so the Cartan geometry 
on $M_0$ is complete.
\end{proof}

\begin{corollary}\label{corollary:completeFromTame}
Suppose that $\Phi : G_0/H_0 \to G_1/H_1$ is an immersive morphism
of homogeneous spaces. Suppose that the locus map 
of $G_0/H_0$ is an immersion (see example~\vref{example:LocusMapOfModelAnImmersion}).
Suppose that $E_1 \to M_1$ is a $G_1/H_1$-geometry
with vanishing $\Phi$-obstruction, vanishing
deformation curvature, and is $\Phi$-tame.
Then the vector fields $\vec{A}$ on $E$, for $A \in \Phi \mathfrak{g}_0$,
are complete.
\end{corollary}

\begin{example}
\begin{theorem}[McKay \cite{McKay2007}]\label{theorem:TameProjectiveConnectionsAreComplete}
Every normal projective connection on a surface,
all of whose geodesics are embedded closed curves,
is tame and complete.
\end{theorem}
\begin{proof}
A projective connection on a surface is a Cartan geometry
modelled on $G/H=\RP{2}$ with $G=\PSL{3,\R{}}$.
It is normal just when certain curvature terms vanish;
see McKay \cite{McKay2007}. The morphism curvature vanishes, as we can easily check. If it
tame then it is complete by corollary~\vref{corollary:completeFromTame}.
By a theorem of LeBrun and Mason \cite{LeBrunMason:2002}, 
a projective connection is tame just precisely when all of the
geodesics are closed and embedded. 
\end{proof}
\begin{remark}
LeBrun and Mason \cite{LeBrunMason:2002} prove that the
only normal projective connection on $\RP{2}$ all
of whose geodesics are closed and embedded curves
is the usual flat connection, i.e. on $\RP{2}$ tame
implies flat and isomorphic to the model.
\end{remark}
\end{example}

\begin{example}
Let's use the previous example to see that
tameness is strictly stronger than completeness.
\begin{theorem}[McKay \cite{McKay2007}]
On any compact manifold, any pseudo-Riemannian metric with positive
Ricci curvature induces a complete normal projective connection.
\end{theorem}
\begin{corollary}
A smooth compact surface of positive Gauss curvature has complete normal projective connection,
but is not tame unless it is a Zoll surface.
\end{corollary}
Since Zoll surfaces are quite rare, we can already see that
tameness is much stronger than completeness.
\end{example}

\begin{example}
\begin{definition}
 A \emph{$G/H$-structure} is a flat $G/H$-geometry.
\end{definition}
\begin{theorem}[Sharpe \cite{Sharpe:1997}]\label{theorem:CompleteImpliesQuasiKleinian}
A $G/H$-structure on a manifold $M$ is complete
if and only if $M=\Gamma\backslash \tilde{G}/H$,
where $\tilde{G}$ is a Lie group containing $H$ and 
$\tilde{G}/H$ has the same local model as 
$G/H$ and $\Gamma \subset \tilde{G}$ is a discrete subgroup 
acting freely and properly on $\tilde{G}/H$.
\end{theorem}
\begin{remark}
 Hence classification of such structures reduces to algebra.
\end{remark}
\begin{corollary}
If a surface $M$ has a flat projective connection,
all of whose geodesics are embedded closed curves,
then the surface is the sphere or real projective
plane with the usual flat projective connection.
\end{corollary}
\begin{proof}
By theorem~\vref{theorem:TameProjectiveConnectionsAreComplete},
the projective connection is complete. Flat
and complete implies isomorphic to the model
up to covering, by theorem~\vref{theorem:CompleteImpliesQuasiKleinian}.
\end{proof}
\end{example}

\begin{example}
A totally geodesic surface is a morphism
of Cartan geometries to any manifold with projective connection, in the obvious way:
in terms of Cartan's structure equations
for projective connections (see Kobayashi
and Nagano \cite{KobayashiNagano:1964}).

\begin{theorem}
If a normal projective connection on a manifold has all geodesics closed and,
for some integer $k>1$, contains
a compact totally geodesic submanifold
of dimension $k$ through each point tangent to each
$k$-plane, then it is a lens space $\Gamma\backslash S^n$,
where $\Gamma$ is a finite group of orthogonal
transformations of $\R{n+1}$, acting freely
on the unit sphere $S^n$.
\end{theorem}
\begin{proof}
It is easy to see that the torsion for such morphisms is just
the curvature of the projective connection.
Therefore the existence of such a
large family of totally geodesic
surfaces immediately implies that
the curvature vanishes.

The induced projective connections
on the totally geodesic surfaces are
normal, because the ambient projective
connection is normal.
The compactness of each surface ensures that
their induced projective connections are all complete. But then
applying corollary~\vref{corollary:completeFromTame},
we find immediately that the ambient projective connection is complete.

Completeness and flatness implies
by theorem~\vref{theorem:CompleteImpliesQuasiKleinian}
that $M$ is the model up to covering.
Since the universal covering of the model
is the sphere, $M$ is covered by the sphere,
say $M=\Gamma\backslash S^n$, where $\Gamma$
is a discrete group acting by projective
automorphisms on $S^n$, freely and properly.
Since $S^n$ is compact, $\Gamma$ is finite,
and therefore $\Gamma$ is, up to projective
conjugacy, a subgroup of the maximal
compact subgroup of the projective
automorphism group of $S^n$, i.e.
$\Gamma$ is a group of orthogonal
transformations.
\end{proof}

A projective Blaschke conjecture:
\begin{conjecture}
If a normal projective connection has all its geodesics
closed embedded curves, then it is the standard normal
projective connection of a compact rank one symmetric space.
\end{conjecture}
\end{example}

\begin{example}
Continuing with example~\vref{example:WilczynskiInvariants},
take a $\OO{n}$-geometry with vanishing Wilczynski invariants.
In that example, we saw that the $\Phi$-torsion for 
the obvious morphism $\Phi$ vanishes and we identified the 
morphism curvature with the Wilczynski invariants.
If the $\OO{n}$-geometry is tame then it is complete 
by corollary~\vref{corollary:completeFromTame}.
It isn't clear how to classify the flat complete or the
flat tame $\OO{n}$-geometries.
\begin{conjecture}
 An $\OO{n}$-geometry is tame if and only if it is isomorphic to a covering
space of the model.
\end{conjecture}
\end{example}

\begin{example}
For certain types of model there is a purely algebraic description
of the flat tame geometries with that model.
\begin{proposition}
Suppose that $\Phi : G_0/H_0 \to G_1/H_1$ is an immersive morphism
of homogeneous spaces. Suppose that the locus map 
of $G_0/H_0$ is an immersion (see example~\vref{example:LocusMapOfModelAnImmersion}).
Suppose that $E_1 \to M_1$ is a $\Phi$-amenable $G_1/H_1$-structure.
Then the vector fields $\vec{A}$ on $E$, for $A \in \Phi \mathfrak{g}_0$,
are complete. If moreover $\mathfrak{h}_1 + \Phi \mathfrak{g}_0$
generates $\mathfrak{g}_1$ as a Lie algebra, then
$M_1$ is complete. In particular, $M_1=\Gamma \backslash \tilde{G}_1/H_1$
as in theorem~\vref{theorem:CompleteImpliesQuasiKleinian}.
\end{proposition}
\begin{proof}
 Corollary~\vref{corollary:completeFromTame} applies
because all of the curvature of $M_1$ vanishes.
If moreover $\mathfrak{h}_1 + \Phi \mathfrak{g}_0$
generates $\mathfrak{g}_1$ as a Lie algebra, then
$M_1$ is flat and complete, so the rest follows
from theorem~\vref{theorem:CompleteImpliesQuasiKleinian}.
\end{proof}
\end{example}

\begin{example}
The obvious inclusion $\SL{2,\R{}} \to \SL{n+1,\R{}}$
as linear transformations of a 2-plane fixing a complementary
$(n-1)$-plane yields a morphism $\Phi : \RP{1} \to \RP{n}$.
We will say that a projective connection is
\emph{tame} if it is $\Phi$-tame.
\begin{theorem}
A manifold with flat complete projective connection is a lens space $\Gamma\backslash S^n$
where $S^n$ is the sphere with its usual flat 
projective connection (the Levi-Civita connection
of the usual round Riemannian metric), and $\Gamma \subset \Or{n+1}$
is a finite subgroup of the orthogonal group $\Or{n+1}$
acting freely on $S^n$.

A manifold with flat tame projective connection is $S^n$ or $\RP{n}$ 
with its usual flat projective connection.
\end{theorem}
\begin{proof}
By theorem~\vref{theorem:CompleteImpliesQuasiKleinian}, if a 
manifold has a flat complete projective connection,
it must be a quotient of the universal covering
space of the model, i.e. of the sphere $S^n$
with its usual flat projective connection,
quotiented by a discrete group $\Gamma$ of automorphisms
of its projective connection. 
Since the sphere is compact, $\Gamma$ must
be a finite group. The automorphisms of the
projective connection are the invertible
linear maps (modulo scaling by positive constants).
We can assume that $\Gamma$ lies in the
maximal compact subgroup of the automorphism
group, the orthogonal group.

To see when such a quotient is tame, we need 
to check that the action of $\Gamma$ on the
Grassmannian of oriented 2-planes in $\R{n+1}$
has no fixed points, or has only fixed points.
This happens just for $\Gamma=\left\{I\right\}$ or $\Gamma=\left\{\pm I\right\}$.
\end{proof}
\end{example}

\begin{example}
A flat conformal geometry has lots of curves which are
locally identified with 2-spheres in the model. (Some
of these have to be thought of as 2-spheres of zero radius,
in order to present them as morphisms.) 
Clearly tameness of these objects implies completeness
of the conformal geometry, and so each path component of
the manifold is a quotient of the model, with
induced conformal geometry, i.e. $\Gamma \backslash S^n$,
where $\Gamma$ is a finite group of rotations. We leave the reader to
work out the relevant Lie algebras, homogeneous spaces and morphism.
\end{example}

\begin{example}
A \emph{$k$-plane field} is a rank $k$ subbundle of the tangent bundle
of a manifold. 
If $P$ is a sheaf of vector fields on a manifold $M$, let $P^{(1)}$ be the 
sheaf of locally defined vector fields spanned
by local sections of $P$ together with Lie brackets $[Y,Z]$ of local
sections of $P$. Define $P^{(2)}, P^{(3)}, \dots$ by induction:
\[
P^{(j+1)}=\left(P^{(j)}\right)^{(1)}.
\]
Let's say that a 2-plane field $P$
on a 5-manifold is \emph{nondegenerate}
if $P^{(1)}$ is a 3-plane field and $P^{(2)}=TM$.

Let $G_2$ be the adjoint form of the split form of the exceptional
simple group of dimension 14.
A nondegenerate 2-plane field 
on a 5-manifold $M$ gives rise to a Cartan geometry
$E \to M$ modelled on $G_2/P_2$ for a certain
parabolic subgroup $P_2 \subset G_2$; see \cite{Cartan:30,Gardner:1989,Sternberg:1983}. 
(The subgroup $P_2$ is unique up to conjugacy.)
We say the plane field is \emph{flat} if this Cartan geometry is flat.
Cartan proved that a nondegenerate 2-plane field 
is flat if and only if it is locally isomorphic
to the model 2-plane field on the model $G_2/P_2$.
(This model 2-plane field is the unique $G_2$-invariant
2-plane field on $G_2/P_2$. The group of diffeomorphisms
of $G_2/P_2$ preserving that 2-plane field 
is a Lie group and has the same Lie algebra
as $G_2$, although it is strictly larger than $G_2$.
It has at least 4 components.) 

Let $B \subset G_2$ be a Borel subgroup 
contained in $P_2$, and let $P_1 \subset G_2$
be the parabolic subgroup of $G_2$ 
containing $B$ and not conjugate to $P_2$
(see \cite{Fulton/Harris:1991} for proof
that $P_1$ exists and is unique and for
an explicit description of $P_1$.)
Then $E/B$ has a foliation by
curves: the leaves of the exterior differential
system $\omega + \mathfrak{p}_1=0$. 
These leaves are curves, known as
\emph{characteristics} \cite{Cartan:30}
or \emph{singular extremals} \cite{Bryant/Hsu:1993}.
Let $W$ be the space of characteristics.
As far as the author knows
$W$ might not be Hausdorff. We will
say that the 2-plane field is \emph{tame}
if $W$ is a smooth manifold and $E/B \to W$
is a smooth circle bundle (not necessarily
a principal circle bundle). Clearly
the model is tame. Clearly tameness
is $\Phi$-tameness for the obvious
choice of $\Phi : P_1/B \to G_2/P_2$. In particular,
by corollary~\vref{corollary:completeFromTame},
tameness implies completeness
of the Cartan geometry.

For example, take any two surfaces
$S$ and $S'$ with Riemannian
metrics. The bundles $OS$ and $OS'$ of orthonormal
bases for the tangent spaces of $S$ and $S'$ are each
equipped with the usual 1-forms (see \cite{BCGGG:1991} p. 89 for
example), say 
$\omega_1, \omega_2, \omega_{12}$
and
$\omega'_1, \omega'_2, \omega_{12}'$
satisfying the structure equations of
Riemannian geometry:
\begin{align*}
 d 
\begin{pmatrix}
 \omega_1 \\
\omega_2
\end{pmatrix}
&=
-
\begin{pmatrix}
 0 & \omega_{12} \\
-\omega_{12} & 0 
\end{pmatrix}
\wedge
\begin{pmatrix}
 \omega_1 \\
\omega_2
\end{pmatrix}
, \
d \omega_{12} = K \omega_1 \wedge \omega_2, \\
 d 
\begin{pmatrix}
 \omega'_1 \\
\omega'_2
\end{pmatrix}
&=
-
\begin{pmatrix}
 0 & \omega'_{12} \\
-\omega'_{12} & 0 
\end{pmatrix}
\wedge
\begin{pmatrix}
 \omega'_1 \\
\omega'_2
\end{pmatrix}
, \
d \omega'_{12} = K' \omega'_1 \wedge \omega'_2.
\end{align*}
The 3-plane field $\mathcal{I} = \left(\omega_1=\omega_1', \omega_2=\omega_2', \omega_{12}=\omega_{12}'\right)$
on $OS \times OS'$ is semibasic for the circle action of rotating
the tangent spaces of $S$ and $S'$ simultaneously at the same rate.
Moreover the paths drawn out in $OS \times OS'$ as orbits of this
circle action are integral curves  of $\mathcal{I}$.
Therefore $\mathcal{I}$ descends to a 2-plane field (also
called $\mathcal{I}$) on the space of orbits. This space of orbits is a 5-manifold, 
which we call $M$, and consists precisely in the
linear isometries between tangent spaces of $S$ and $S'$. 
We can think of the paths of the resulting $2$-plane field 
as ``rolling of $S$ along a path on $S'$,'' or vice versa.
Clearly the isometries of $S$ and those of $S'$ act
on $M$ preserving the 2-plane field.

It is well known (see \cite{Bryant/Hsu:1993}) that $\mathcal{I}$
is a nondegenerate just above those points where $S$ and $S'$
have distinct Gauss curvature. There is a natural family
of $\mathcal{I}$-tangent curves: pick an oriented geodesic on $S$,
and a point $p$ on that geodesic, an oriented geodesic on $S'$
and a point $p'$ on that geodesic, and a linear isometry
$T_p S \to T_{p'} S'$ identifying the unit tangent
tangent vectors of the geodesics, and extend that 
identification by parallel transport along both
geodesics. It is well known that the $\Phi$-morphisms, i.e.
the characteristics, are precisely
these $\mathcal{I}$-tangent curves in $M$.

\begin{conjecture}
Suppose that $S$ and $S'$ are surfaces with 
different Gauss curvatures at any pair of points, 
and $M$ is the space of linear isometries 
of tangent planes of $S$ and $S'$.
Then the $G_2/P_2$-geometry of the nondegenerate 2-plane
field on the space  is complete just when $S$ and $S'$ have
complete Riemannian metrics.
\end{conjecture}

\begin{conjecture}
Suppose that $S$ and $S'$ are surfaces with 
different Gauss curvatures at any pair of points, 
and $M$ is the space of linear isometries 
of tangent planes of $S$ and $S'$.
Then the $G_2/P_2$-geometry of the nondegenerate 2-plane
field on $M$ is tame just when both surfaces 
are Zoll surfaces (i.e. all of their geodesics are closed
and embedded curves) and the length $L$ of all closed
geodesics on $S$ is a rational multiple of the length $L'$
of all closed geodesics on $S'$.
\end{conjecture}

Consider a special case: let $S$ be a sphere
of radius 1 and $S'$ be a sphere of radius 3,
and let $M$ be the space of linear isometries
between tangent spaces of $S$ and $S'$.
It is well known \cite{Agrachev:2007,Sagerschnig:2006}
that the model $G_2/P_2$  is a single component of the space $M$
consisting precisely in the linear isometries which preserve orientation. 
Clearly the ``phase space'' manifold $M$ has two components, so consists of two copies of $G_2/P_2$.
Isometries of the sphere of radius 1 and of the sphere of radius 3 act on $M$ preserving
the 2-plane field, so their orientation preserving components sit in $G_2$, giving
the maximal compact subgroup $\SO{3} \times \SO{3} \subset G_2$.
The automorphism group of the 2-plane field contains $\Or{3} \times \Or{3}$,
and has Lie algebra $\mathfrak{g}_2$ by Cartan's work \cite{Cartan:30}.
The maximal compact subgroup therefore has Lie algebra $\so{3} \times \so{3}$,
so the automorphism group of the 2-plane field has at least 4 components.
It does not appear to be known in the literature what the 
automorphism group of Cartan's example is.
It is known \cite{Agrachev:2007,Sagerschnig:2006} 
that the simply connected double cover of $G_2/P_2$ is $S^2 \times S^3$.

\begin{lemma}
The characteristics of any flat complete nondegenerate 2-plane field are closed and embedded.
\end{lemma}
\begin{proof}
By theorem~\vref{theorem:CompleteImpliesQuasiKleinian} every flat complete $G_2/P_2$-geometry
on a connected manifold has the form $M=\Gamma \backslash \tilde{G}_2/P_2$, where
$\Gamma$ is a discrete subgroup of the covering group $\tilde{G}_2$
of $G_2$ which acts on the universal covering space of $G_2/P_2$. Clearly $\tilde{G}_2/P_2$ is simply connected and compact. 
It is enough to check therefore that the characteristics
of $\tilde{G}_2/P_2$ are closed, which is easy to check. They are embedded because they are homogeneous.
\end{proof}

\begin{theorem}
There are infinitely many nonhomeomorphic compact 5-manifolds admitting flat tame 
nondegenerate 2-plane fields, and there are infinitely many 
nonhomeomorphic compact 5-manifolds admitting flat complete but \emph{not} tame
nondegenerate 2-plane fields.
\end{theorem}
\begin{proof}
By theorem~\vref{theorem:CompleteImpliesQuasiKleinian} every flat complete $G_2/P_2$-geometry
on a connected manifold has the form $M=\Gamma \backslash \tilde{G}_2/P_2$, where
$\Gamma$ is a discrete subgroup of the covering group $\tilde{G}_2$
of $G_2$ which acts on the universal covering space of $G_2/P_2$. Clearly $\tilde{G}_2/P_2$ is simply connected and compact.
Therefore $\Gamma$ must be finite.  Without loss of generality, $\Gamma$ is a finite subgroup of the
maximal compact subgroup of $\tilde{G}_2$,
i.e. of $\left(\SU{2} \times \SU{2}\right)/\left<(-1,-1)\right>$.

% The finite subgroups of $\SO{3}$
% are well known \cite{Burban:2009}: the rotation symmetries
% of the $n$-sided pyramid and $2n$-sided
% double pyramid and of the Platonic
% solids. To be more precise, the
% finite subgroups of $\SO{3}$, up to
% conjugacy, labelled by orbifold notation, are:
% \begin{enumerate}
% \item $1$: the identity group (order 1),
%  \item $nn$: the cyclic group generated
% by rotation by angle $2\pi/n$ around an axis (order $n$),
% \item $22n$: the dihedral group generated
% by a rotation by angle $2\pi/n$ around
% an axis, and by rotation by angle $\pi$
% around a line perpendicular to that axis (order $2n$),
% \item $332$: the group of rotational
% symmetries of the tetrahedron (order $12$),
% \item $432$: the group of rotational
% symmetries of the octahedron (order $24$),
% \item $532$: the group of rotational
% symmetries of the icosahedron (order $60$).
% \end{enumerate}

The finite subgroups of $\SU{2}$ are
well known: essentially they are the preimages
via $\SU{2} \to \SO{3}$ of the finite
subgroups of $\SO{3}$. To be precise,
the finite subgroups of $\SU{2}$, up
to conjugacy, are:
\begin{enumerate}
\item
 The cyclic group generated by 
\[
\begin{pmatrix}
 \varepsilon & 0 \\
0 & \frac{1}{\varepsilon}
\end{pmatrix}, \varepsilon=e^{2 \pi i/n}
\]
\item
the binary dihedral group generated by
\[
\begin{pmatrix}
 \varepsilon & 0 \\
0 & \frac{1}{\varepsilon}
\end{pmatrix},
\begin{pmatrix}
 0 & 1 \\
-1 & 0
\end{pmatrix}, \varepsilon=e^{2 \pi i/n}
\]
\item
the binary tetrahedral group generated by
\[
 \begin{pmatrix}
  i & 0  \\
  0 & -i
 \end{pmatrix},
\begin{pmatrix}
 0 & 1 \\
-1 & 0
\end{pmatrix},
\frac{1}{\sqrt{2}}
\begin{pmatrix}
 \varepsilon^{7}  &  \varepsilon^{7} \\
 \varepsilon^{5}  & \varepsilon  \\
\end{pmatrix}, 
\varepsilon = e^{2 \pi i/8}
\]
\item
the binary octahedral group generated by
\[
 \begin{pmatrix}
  i & 0  \\
  0 & -i
 \end{pmatrix},
\begin{pmatrix}
 0 & 1 \\
-1 & 0
\end{pmatrix},
\frac{1}{\sqrt{2}}
\begin{pmatrix}
 \varepsilon^{7}  &  \varepsilon^{7} \\
 \varepsilon^{5}  & \varepsilon  \\
\end{pmatrix},
 \begin{pmatrix}
  \varepsilon & 0  \\
  0 & \varepsilon^{7}
 \end{pmatrix},
\varepsilon = e^{2 \pi i/8}
\]
\item
the binary icosahedral group generated by
\[
-
 \begin{pmatrix}
  \varepsilon^3 & 0  \\
  0 & \varepsilon^2
 \end{pmatrix},
\frac{1}{\sqrt{5}}
\begin{pmatrix}
-\varepsilon+\varepsilon^4 & \varepsilon^2-\varepsilon^3 \\
\varepsilon^2-\varepsilon^3 & \varepsilon-\varepsilon^4 
\end{pmatrix}, \varepsilon=e^{2 \pi i/5}.
\]
\end{enumerate}

Clearly $G_2/P_2=\left(\SO{3} \times \SO{3}\right)/S_2$
for the diagonal circle subgroup $S_2 \subset \SO{3} \times \SO{3}$,
consisting in simultaneous rotations of the two spheres
about the vertical axis, in terms of Cartan's picture of 
spheres of radius $1$ and $3$, pictured sitting tangent
to one another, south pole to north pole. In terms
of quaternions, this is the subgroup
\[
\left(e^{\theta k}, e^{\theta k}\right)
\]
for $0 \le \theta < 2 \pi$.

Clearly $G_2/P_1=\left(\SO{3} \times \SO{3}\right)/S_1$
for the circle subgroup $S_1 \subset \SO{3} \times \SO{3}$,
consisting in rotating one sphere around a horizontal
axis, and the other sphere in the opposite direction
around the same axis at a third the speed, rolling them on each other so
that they travel in opposite directions. In terms
of quaternions, this is the subgroup
\[
\left(e^{3 \theta j}, e^{-\theta j}\right)
\]
for $0 \le \theta < 2 \pi$.

The problem is to figure out which finite subgroups
of $\SU{2} \times \SU{2}$ will turn out to give
a Haussdorff quotient of one or the other of these
5-manifolds $\tilde{G}_2/P_1$ and $\tilde{G}_2/P_2$, i.e. act freely. Those which act
freely on both 5-manifolds are precisely the flat tame
2-plane fields.

So far, the candidate groups $\Gamma$ are precisely
the finite subgroups of $\Gamma_1 \times \Gamma_2$,
where $\Gamma_1, \Gamma_2 \subset \SO{3}$ are 
finite subgroups. We have to check which such groups
$\Gamma$ act freely on $G_2/P_1$ and which act
freely on $G_2/P_2$. Let's try $G_2/P_2$ first.
In other words, we need to ensure that if an element
\[
 \gamma=\left(\gamma_1,\gamma_2\right) \in \Gamma
\]
fixes and element $\left(g_1,g_2\right)S_2 \in \left(\SU{2} \times \SU{2}\right)/S_2$,
then it fixes all such elements. Otherwise the
quotient $\Gamma \backslash \left(\SU{2} \times \SU{2}\right)/S_2$
will not be Hausdorff. Fixed points will be those satisfying
\begin{align*}
\gamma_1 g_1 &= g_1 e^{k \theta} \\
\gamma_2 g_2 &= g_2 e^{k \theta}, 
\end{align*}
so that in particular, $g_1^{-1} \gamma_1 g_1 = g_2^{-1} \gamma_2 g_2=e^{k \theta}$,
and so $\gamma_1$ and $\gamma_2$ are conjugate in $\SU{2}$.
Working only up to conjugacy, 
since every element of $\SU{2}$ is conjugate to something
of the form $e^{k \theta}$,
an element $\left(\gamma_1,\gamma_2\right) \in \Gamma$ has fixed
points in $G_2/P_2$ just when $\gamma_1$ and $\gamma_2$
are conjugate.

Elements of $\SU{2}$ are conjugate just when they 
have the same eigenvalues as $2 \times 2$ matrices. All elements of 
$\SU{2}$ are diagonalizable. 
The eigenvalues of elements of a finite subgroup
of $\SU{2}$ can be any root of unity. So finally, $\Gamma$ 
acts freely on $G_2/P_2$ just when, for every element
$\left(\gamma_1,\gamma_2\right) \in \Gamma$,
either $\gamma_1$ and $\gamma_2$ have distinct
eigenvalues, or else $\gamma=(-1,-1)$ or $\gamma=(1,1)$.

A similar analysis says that $\Gamma$ acts
freely on $G_2/P_1$ just when, for every element
$\left(\gamma_1,\gamma_2\right) \in \Gamma$,
either $\gamma_1^{-1}$ and $\gamma_2^{3}$ have distinct
eigenvalues, or else $\gamma=(-1,-1)$ or $\gamma=(1,1)$.

It seems likely that one could analyse all of the possible
finite subgroups $\Gamma$ to determine which act
freely on $G_2/P_1$ and which act freely on $G_2/P_2$. We will
leave this to the reader and simply give examples of cyclic
groups $\Gamma$ for which these actions are either free
or not free. Suppose that $\Gamma$ is generated by
an element $\gamma=\left(\gamma_1,\gamma_2\right) \in \SU{2} \times \SU{2}$
where $\gamma_1$ and $\gamma_2$ are of finite order,
say with eigenvalues $e^{2 \pi i p_1/q_1}$ and $e^{2 \pi i p_2/q_2}$
respectively, where $p_1, q_1, p_2$ and $q_2$ are
integers. Without loss of generality, $0 \le p_1 < q_1$
and $0 \le p_2 < q_2$.  

In order to ensure that $\Gamma$ acts freely on $G_2/P_2$,
we need to ensure precisely that, for an integer $k$,
if
\[
 k \left(\frac{p_1}{q_1} - \frac{p_2}{q_2}\right)
\]
is an integer (i.e. $\gamma_1^k$ is conjugate to 
$\gamma_2^k$), then 
\[
k \frac{p_1}{q_1} \text{ and } k \frac{p_2}{q_2}
\]
are equal half integers (i.e. $\gamma_1^k$ and $\gamma_2^k$
are both $\pm 1$).

There are three possibilities here:
\begin{enumerate}
\item $p_1=p_2=0$, i.e. $\gamma_1=\gamma_2=1$ or
 \item $\frac{p_1}{q_1}=\frac{p_2}{q_2}=\frac{1}{2}$ 
i.e. $\gamma_1=\gamma_2=-1$ or
\item if we write $\frac{p_1}{q_1}-\frac{p_2}{q_2}$
as $\frac{r}{s}$ in lowest terms, then 
if $s$ divides $kr$, 
we need $q_1$ to divide $2 k p_1$ and $q_2$ to divide $2 k p_2$.
\end{enumerate}
The second possibility is clearly one that leads precisely
to Cartan's model. Let's consider the second possibility
more carefully.

Since $p_1$ and $q_1$ are relatively prime, to 
have $q_1$ divide $2 k p_1$ is precisely to have
$q_1$ divide $2k$. Therefore $\Gamma$ acts
freely on $G_2/P_2$ just when for any integer $k$, if 
$s$ divides $k$, then both $q_1$ and $q_2$ divide
$2k$, i.e. we need $q_1$ and $q_2$ to both divide
$2s$. This quantity $s$ is just the denominator of
$\frac{p_1}{q_1}-\frac{p_2}{q_2}$ in lowest terms,
i.e.
\[
s = \frac{q_1 q_1}{\operatorname{gcd}\left(q_1q_2, p_1q_2 - q_1 p_2 \right)}.
\]
So $\Gamma$ acts freely precisely when
$q_1$ and $q_2$ both divide the integer
\[
\frac{2 q_1 q_1}{\operatorname{gcd}\left(q_1q_2, p_1q_2 - q_1 p_2 \right)}.
\]
This occurs just when the least common multiple $\ell=\operatorname{lcm}\left(q_1,q_2\right)$
divides
\[
\frac{2 q_1 q_1}{\operatorname{gcd}\left(q_1q_2, p_1q_2 - q_1 p_2 \right)},
\]
i.e. when 
\[
\frac{q_1 q_2}{\operatorname{gcd}\left(q_1,q_2\right)}
\text{ divides }
\frac{2 q_1 q_1}{\operatorname{gcd}\left(q_1q_2, p_1q_2 - q_1 p_2 \right)},
\]
i.e. when
\[
 \operatorname{gcd}\left(q_1q_2, p_1q_2 - q_1 p_2 \right) \text{ divides } 2 \operatorname{gcd}\left(q_1,q_2\right).
\]

So finally we find that $\Gamma$ acts freely on $G_2/P_2$
just when
\[
 \operatorname{gcd}\left(q_1q_2, p_1q_2 - q_1 p_2 \right) \text{ divides } 2 \operatorname{gcd}\left(q_1,q_2\right).
\]
and acts freely on $G_2/P_1$ just when
\[
 \operatorname{gcd}\left(q_1q_2, p_1q_2 + 3 \, q_1 p_2 \right) \text{ divides } 2 \operatorname{gcd}\left(q_1,q_2\right).
\]

For example, if we take any prime number $q$ and let $p_1/q_1=2/q$ and $p_2/q_2=1/q$,
then we will find that $\Gamma$ acts freely on both $G_2/P_1$ and $G_2/P_2$. Hence
there are infinitely many nonhomeomorphic compact 5-manifolds admitting flat tame 
nondegenerate 2-plane fields.

On the other hand, if we take any prime number $q>2$ and pick any number
$1 \le p_2 < q$ and let $p_1/q_1=1-3 p_2/q$, then we find that 
$\Gamma$ acts freely on $G_2/P_2$ but does \emph{not} act
freely on $G_2/P_1$. Hence there are infinitely many nonhomeomorphic compact 5-manifolds admitting flat complete
nontame nondegenerate 2-plane fields.

We leave to the reader the problem of classifying all
of the possible finite subgroups $\Gamma \subset \SU{2} \times \SU{2}$
which act freely on $G_2/P_1$ or on $G_2/P_2$.
\end{proof}

\begin{remark}
 It is not known which compact 5-manifolds admit
a nondegenerate 2-plane field. It is not even 
known which compact 5-manifolds admit a \emph{flat} 
nondegenerate 2-plane field. It is not
known whether every nondegenerate 2-plane field
can be deformed into every other nondegenerate
2-plane field through a connected family
of nondegenerate 2-plane fields.
\end{remark}

\begin{remark}
It is not known how to identify whether a nondegenerate 2-plane
field on a 5-manifold is constructed by rolling a surface
on another surface, or how invariants of surface geometry
appear in the curvature of the $G_2/P_2$-geometry.
\end{remark}

\begin{remark}
 The computation of the morphism curvature for
2-plane fields has never been done. It might
provide some interpretation of the curvature
of the $G_2/P_2$-geometry of nondegenerate 2-plane fields.
I conjecture that the morphism curvature of a 
nondegenerate 2-plane field vanishes
just when the curvature vanishes, i.e.
just when the 2-plane field is locally
isomorphic to the 2-plane field on $G_2/P_2$.
Roughly speaking, the curvature of a nondegenerate
2-plane field is conjecturally the infinitesimal deformation equation
for characteristics. 
\end{remark}

\begin{remark}
 It is not known if it is possible to find
pairs $\left(M_1,M_2\right)$ and $\left(N_1,N_2\right)$
of surfaces with Riemannian (or Lorentzian) metrics
so that the associated 2-plane field for
rolling $M_1$ on $M_2$ is (locally) isomorphic
to that for rolling $N_1$ on $N_2$,
but so that $M_1$ is not locally isometric
to either $N_1$ or $N_2$ up to constant scaling. 
\end{remark}
\end{example}

\section{Proof of the rolling theorems}

\begin{definition}
Suppose that $\left(H,\mathfrak{g}\right)$ is a local model.
Pick any positive definite scalar product
on $\mathfrak{g}$. This scalar product determines a Riemannian
metric on the total space $E$ of any $\left(H,\mathfrak{g}\right)$-geometry 
$E \to M$, the canonical metric of the coframing given by the Cartan connection.
Such a Riemannian metric is called a \emph{canonical metric}
of the Cartan geometry.
\end{definition}

\begin{definition}
Suppose that $\Phi : \left(H_0,\mathfrak{g}_0\right) \to \left(H_1,\mathfrak{g}_1\right)$ 
is a local model morphism. Suppose that $E_0 \to M_0$ and
$E_1 \to M_1$ are Cartan geometries with those local models. 
A $\Phi$-\emph{development} is a choice of
\begin{enumerate}
\item manifold $X$ and
\item map $f_0 : X \to M_0$ and
\item map $f_1 : X \to M_1$ and
\item $H_0$-equivariant map $F : f_0^* E_0 \to f_1^* E_1$ and
\item $F$ must satisfy $F^* \omega_1= \Phi \omega_0$.
\end{enumerate}
\end{definition}

\begin{definition}
If a given map $f_0 : X \to M_0$ has a $\Phi$-development,
we say that it \emph{$\Phi$-develops}.
\end{definition}

\begin{definition}
We will say that \emph{curves $\Phi$-develop freely}
from $M_0$ to $M_1$ if every curve $f_0 : C \to M_0$
has a $\Phi$-development.
\end{definition}

\begin{lemma}\label{lemma:CoframingToCartanGeometry}
The following result holds in the real or complex analytic categories.
 Suppose that $\Phi : \left(H_0,\mathfrak{g}_0\right) \to \left(H_1,\mathfrak{g}_1\right)$
is a local model morphism. Suppose that $E_0 \to M_0$ and
$E_1 \to M_1$ are Cartan geometries with those local models.

Suppose that we have a curve $f_0 : C \to M_0$,
and $C$ is simply connected. (In the complex
case, we add the assumption that $C \ne \CP{1}$.) 
Then the bundle
$f_0^* E_0 \to C$ is trivial. Let $s$ be a section.
Let $\iota : f_0^* E_0 \to E_0$ be the obvious map.
The following are equivalent:
\begin{enumerate}
 \item 
$f_0$ will $\Phi$-develop in the sense of Cartan geometries,
\item
the map $\iota \circ s : C \to E_0$ will $\Phi$-develop 
in the sense of coframings,
\item
there is some section $s_0$ of $f_0^* E_0 \to C_0$ so that 
the map $\iota \circ s_0 : C \to E_0$ will $\Phi$-develop 
in the sense of coframings.
\end{enumerate}
\end{lemma}
\begin{proof}
Clearly $f_0^* E_0 \to C$ must be trivial (holomorphically
trivial in the complex analytic category by classification
of Riemann surfaces and of the principal bundles on 
the disk and complex affine line).

Suppose (1): that we can develop $f_0$ to maps $f_1 : C \to M_1$
and $F : f_0^* E_0 \to f_1^* E_1$. 
Let $\iota_0 : f_0^* E_0 \to E_0$ and $\iota_1 : f_1^* E_1 \to E_1$
be the obvious maps.
Then $\iota_1 \circ F \circ \iota_0 \circ s : C \to E_1$
is a $\Phi$-development of $s$ in the sense of coframings.
Therefore (1) implies (2) and (3). Clearly (2) implies (3).

Suppose (3): that some section $s_0$ develops in the sense
of coframings. So suppose it develops to $s_1 : C \to E_1$,
with $s_1^* \omega_1 = \Phi s_0^* \omega_0$. Define a 
map $F : f_0^* E_0 \to E_1$ by $F\left(s_0 h_0\right)=s_1\Phi\left(h_0\right)$.
Denote the map $E_1 \to M_1$ as $\pi_1 : E_1 \to M_1$.
Define a map $f_1 : C \to M_1$ by $f_1 = \pi_1 \circ F \circ s_0$.
It is easy to check that the maps $f_0, f_1, F$ form a $\Phi$-development.
Therefore (3) implies (1).
\end{proof}

We now generalize and thereby prove theorem~\vref{thm:Main}.
\begin{theorem}\label{theorem:Rolling}
For a Cartan geometry $E_1 \to M_1$, with Cartan connection $\omega_1$, 
modelled on a homogeneous space, the following
are equivalent:
\begin{enumerate}
 \item\label{item:AcompleteMetric}
some canonical metric on $E_1$ is complete,
\item\label{item:AllcompleteMetrics}
every canonical metric on $E_1$ is complete,
\item\label{item:coframingComplete}
the coframing $\omega_1$ is complete,
\item\label{item:IdevelopModel}
curves $I$-develop to $M_1$ from the model,
\item\label{item:IdevelopSome}
curves $I$-develop to $M_1$ from some Cartan geometry,
\item\label{item:PhiDevelopAll}
curves $\Phi$-develop to $M_1$ from any Cartan geometry, for any model morphism $\Phi$.
\end{enumerate}
\end{theorem}
\begin{proof}
By proposition~\vref{proposition:Main}, (\ref{item:AcompleteMetric}), 
(\ref{item:AllcompleteMetrics}) and (\ref{item:coframingComplete}) are
equivalent. Each of them implies that all curves
$C \to E_0$ develop in the sense of coframings
to curves $C \to E_1$, by proposition~\ref{proposition:Main} again.
They are therefore equivalent to existence of developments
of all curves by lemma~\vref{lemma:CoframingToCartanGeometry},
i.e. to (\ref{item:PhiDevelopAll}). Clearly (\ref{item:PhiDevelopAll})
implies (\ref{item:IdevelopModel}) and (\ref{item:IdevelopSome}).
Obviously (\ref{item:IdevelopModel}) implies (\ref{item:IdevelopSome}).
So we have only to show that (\ref{item:IdevelopSome}) implies
(\ref{item:AcompleteMetric}). 
By lemma~\ref{lemma:CoframingToCartanGeometry}, (\ref{item:IdevelopSome})
implies that one can development all curves from some coframing.
By proposition~\ref{proposition:Main}, this implies (\ref{item:AcompleteMetric}).
\end{proof}

We now prove the complex analytic analogue 
of theorem~\vref{thm:Main}:

We now generalize and thereby prove theorem~\ref{theorem:Immersed}.
\begin{theorem}
For a holomorphic Cartan geometry $E_1 \to M_1$, with Cartan
connection $\omega_1$, modelled on a homogeneous space, the following
are equivalent:
\begin{enumerate}
 \item\label{item:AcompleteMetricC}
some canonical metric on $E_1$ is complete,
\item\label{item:AllcompleteMetricsC}
every canonical metric on $E_1$ is complete,
\item\label{item:coframingCompleteC}
the coframing $\omega_1$ is complete,
\item\label{item:IdevelopModelC}
curves $I$-develop to $M_1$ from the model,
\item\label{item:IdevelopSomeC}
curves $I$-develop to $M_1$ from some Cartan
geometry,
\item\label{item:PhiDevelopAllC}
curves $\Phi$-develop to $M_1$ from any Cartan geometry, for any model morphism $\Phi$.
\end{enumerate}
\end{theorem}

\begin{proof}
The proof is identical to the proof of theorem~\vref{theorem:Rolling}
with one small problem: if $f_0 : C \to M_0$ 
is a curve, and $C=\CP{1}$, then we can't directly employ
lemma~\ref{lemma:CoframingToCartanGeometry}. However,
at each stage where we need lemma~\ref{lemma:CoframingToCartanGeometry},
we can first pop off one point of $C$, and develop the
restriction to $C \backslash \text{pt}=\C{}$. The
development is given by holomorphic first order ordinary
differential equations, so is unique when it exists.
But then we can change our choice of point to pop off,
and by uniqueness we will be able to patch the
development of this other affine chart together
with the first one to give a global development of $f_0$.
\end{proof}

\section{Compact structure group}

We will slightly generalize Clifton's results \cite{Clifton:1966} on Euclidean connections.

\begin{definition}\label{definition:CanonicalMetricOnCartanGeometry}
Suppose that $\left(H,\mathfrak{g}\right)$ is a local model and that $H$ is compact. Pick an $H$-invariant
inner product on $\mathfrak{g}$ (the Lie algebra of $G$). Write
\[
\mathfrak{g} = \mathfrak{h} \oplus \mathfrak{g}/\mathfrak{h}
\]
as an orthogonal decomposition in the inner product. 

Suppose that $H \to E \to M$ is a Cartan geometry with
Cartan connection $\omega$. Let $E$ have the canonical
metric derived from the inner product on $\mathfrak{g}$.
Then take the orthogonal complement $V^{\perp} \subset TE$ to the
tangent spaces $V$ of the fibers of $E \to M$. This orthogonal complement
is preserved by $H$-action, since the metric is $H$-invariant.
The projection $\pi : E \to M$ identifies $V^{\perp} = \pi^* TM$,
and the map $\pi' : V^{\perp} \to TM$ is therefore an isometry for 
a unique Riemannian metric on $M$, to be called the \emph{canonical
metric} on $M$.
\end{definition}

\begin{lemma}
Suppose that $H \to E \to M$ is a Cartan geometry 
with local model $\left(H,\mathfrak{g}\right)$ and that $H$ is compact.
Then $E \to M$ is a Riemannian submersion, i.e. a submersion
for which the map $V^{\perp} \to TM$ on the orthogonal complement
to the fibers is an isometry of vector bundles with metric.
In particular, lengths of curves contract under the mapping 
$E \to M$.
\end{lemma}

\begin{corollary}
If $H \to E \to M$ is a $\left(H,\mathfrak{g}\right)$-geometry 
and $H$ is compact, then any Cartan geometry modelled
on $\left(H,\mathfrak{g}\right)$ rolls on 
$M$ just when $M$ is complete in the canonical Riemannian metric.
\end{corollary}
\begin{proof}
The map $E \to M$ is a Riemannian submersion, so if closed balls in $E$ are compact,
they will map to balls in $M$ of the same radius, which will then also be compact.
Conversely, if $M$ is complete in the canonical Riemannian metric,
then the ball of any given radius about a point of $E$ will map to a ball of the
same radius in $M$, a compact set, so will live in the preimage of that ball.
Because the fibers of $E \to M$ are compact, this preimage will also be compact. 
\end{proof}

\section{Monodromy}

\begin{definition}
Take any connected manifold $X$. Let $\pi : \tilde{X} \to X$ be the
universal covering map. Suppose that $\phi_0 : X \to M_0$ is any smooth map,
and $M_0$ and $M_1$ bear $\left(H,\mathfrak{g}\right)$-geometries. Let $\tilde{\phi}_0 =
\phi_0 \circ \pi : \tilde{X} \to M_0$. Assume that there is a
development $\tilde{\phi}_1 : \tilde{X} \to M_1$. Then we have an
isomorphism $\tilde{\phi}_0 ^* E_0 = \tilde{\phi}_1^* E_1$. The map
$\tilde{\phi}_0 ^* E_0 \to \phi_0 ^* E_0$ is a covering map, with
covering group $\pi_1(X)$. Moreover, $\pi_1(X)$ acts on
$\tilde{\phi}_0^* E_0=\tilde{\phi}_1^* E_1$ as bundle automorphisms
over the deck transformations of $\tilde{X}$.  We refer to this
action as the \emph{monodromy} of the development.
\end{definition}
Picking a frame $\tilde{e}_0 \in \tilde{\phi}_0 ^* E_0$ and
corresponding $\tilde{e}_1 \in \tilde{\phi}_1 ^* E_1$, and
corresponding points $e_0 \in \phi_0^* E_0$, $\tilde{x} \in
\tilde{X}$ and $x \in X$, we will examine the monodromy orbit of
$\tilde{e}_1$.

\begin{lemma}
The monodromy of the development is a free and proper
action preserving $\tilde{\phi}_0^* \omega$.
Two elements of $\pi_1(X)$ have the same monodromy action on all of
$\tilde{\phi}_1^* E_1$ just when they have the same monodromy action
on some element of $\tilde{\phi}_1^* E_1$.
\end{lemma}
\begin{proof}
If two elements have the same effect on $\tilde{e}_1$, then
composing one with the inverse of the other produces an element
$\gamma$ fixing $\tilde{e}_1$. The action of $\pi_1(X)$ commutes
with the action of $H$, so $\gamma$ fixes every element of the $H$-orbit
through $\tilde{e}_1$. Moreover, $\gamma$ fixes $\tilde{\phi}_1^* \omega$.

Take any path $p(t)$ in $\tilde{\phi}_1^* E_1$ with
$p(0)=\tilde{e}_1$. Let $A(t) = \dot p(t) \hook \tilde{\phi}_1^*
\omega$. The only solution $q(t)$ to $\dot q(t) \hook
\tilde{\phi}_1^* \omega = A(t)$ satisfying $q(0)=\tilde{e}_1$ is $p(t)$.
This differential equation and initial condition are $\gamma$
invariant, and therefore all points of $p(t)$ are fixed by $\gamma$.
Therefore $\gamma$ fixes every point in the path component of the
$H$-orbit of $\tilde{e}_1$, i.e. $\gamma$ acts trivially on
$\tilde{\phi}_1 ^* E_1$.
\end{proof}

\begin{lemma}\label{lemma:equiv}
Suppose that $f : X \to Y$ is a smooth map of manifolds, equivariant for
free and proper actions of a group $\Gamma$ on $X$ and $Y$. Then
\begin{enumerate}
\item the quotient spaces $\bar{X} = \Gamma \backslash X$ and $\bar{Y} = \Gamma \backslash Y$ are smooth manifolds,
\item the obvious maps $X \to \bar{X}$ and $Y \to \bar{Y}$ are $\Gamma$-bundles,
\item the quotient map $\bar{f} : \bar{X} \to \bar{Y}$ on $\Gamma$-orbits is a smooth map,
and
\item any smooth local sections of $X \to \bar{X}$ and $Y \to \bar{Y}$ identify $f$ and $\bar{f}$;
in particular
\[
\operatorname{rk} f'(x) = \dim \Gamma + \operatorname{rk} \bar{f}'\left(\bar{x}\right)
\]
when $x \in X$ maps to $\bar{x} \in \bar{X}$.
\end{enumerate}
\end{lemma}
\begin{proof}
It is well known that $X \to \bar{X}$ is a smooth principal bundle for a unique smooth structure on $\bar{X}$,
and by the same token for $Y \to \bar{Y}$.
Clearly $\bar{f}$ is well defined, and continuous, and lifts to $f$ under any smooth local sections
of $X \to \bar{X}$ and $Y \to \bar{Y}$, so is smooth.
\end{proof}

\begin{theorem}
Take any connected manifold $X$. Suppose that $M_0$ and $M_1$ are manifolds bearing
$G/H$-geometries. Take $\phi_0 : X \to M_0$ any smooth map. Assume
that there is a development of some covering space of $X$ to $M_1$.
Then a (possibly different) covering space $\hat{X} \to X$ develops from $M_0$ to $M_1$
just when the monodromy orbit of $\pi_1\left(\hat{X}\right)$ on some
point $\tilde{e}_1 \in \tilde{\phi}_1 ^* E_1$ maps to a single point
in $E_1$.
\end{theorem}
\begin{proof}
We can replace $X$ by $\hat{X}$ if needed to arrange that $X=\hat{X}$ without
loss of generality.
Clearly $X$ develops from $M_0$ to $M_1$ just when $\tilde{\phi}_1$
is $\pi_1\left(X\right)$-invariant. Moreover, if $X$ develops, then
the monodromy orbit of $\pi_1\left(X\right)$ on any point $\tilde{e}_1$ must
map to a single point in $E_1$. Suppose that the monodromy orbit of
$\pi_1(X)$ on $\tilde{e}_1$ maps to a single point of $E_1$. Let $x \in X$ be the point
of $X$ which is the image of $\tilde{e}_1$ under the obvious bundle map
$\tilde{\phi}_1^* E_1 \to X$. The map
$\tilde{\phi}_1^* E_1 \to E_1$ is $H$-equivariant, so every
monodromy orbit of $\pi_1(X)$ above $x$ is mapped
to a single point of $E_1$. Take a point $\tilde{e}_1 \in \tilde{\phi}_1^*
E_1$, and suppose that $\tilde{e}_1$ maps to $e_1 \in E_1$. Take any
smooth path $p(t)$ in $\tilde{\phi}_1^* E_1$ with $p(0)=\tilde{e}_1$. Let
$A(t) = \dot p(t) \hook \tilde{\phi}_1^* \omega$. The only solution
$q(t)$ to $\dot q(t) \hook \tilde{\phi}_1^* \omega = A(t)$
satisfying $q(0)=\tilde{e}_1$ is $p(t)$. Therefore in $E_1$, the
only solution $q(t)$ to $\dot q(t) \hook \omega = A(t)$ satisfying
$q(0)=e_1$ is the image in $E_1$ of $p(t)$. The monodromy group will
move $\tilde{e}_1$ around the monodromy orbit, moving $p(t)$ to
another curve, but yielding the same solution curve $q(t)$ in $E_1$.
Therefore the map $\tilde{\phi}_1^* E_1 \to E_1$ is invariant under
$\pi_1(X)$ and drops to a smooth map on the quotient: $\pi_1(X) \backslash \tilde{\phi}_1^* E_1 \to E_1$.

The development $\tilde{\phi}_0^* E_0 = \tilde{\phi}_1^* E_1$ identifies
$\phi_0^* E_0 = \pi_1(X) \backslash \tilde{\phi}_0^* E_0= \pi_1(X) \backslash \tilde{\phi}_1^* E_1$.
By equivariance under the action of $H$, our map descends to a smooth map
$\phi_1 : X=\pi_1(X) \backslash \tilde{\phi}_1^* E_1/H \to E_1/H=M_1$.
Define a map $\Phi : \phi_0^* E_0 \to \phi_1^* E_1$ by quotienting $\tilde{\Phi}$ to $\pi(X)$ orbits.
By lemma~\ref{lemma:equiv}, $\Phi$
is a local diffeomorphism, an $H$-bundle isomorphism, and satisfies $\omega_0=\omega_1$, so a development.
\end{proof}

\section{Rolling along Lipschitz CW-complexes}

\begin{definition}
A locally Lipschitz or $C^k$ or smooth or analytic (etc. CW-complex is a
CW-complex whose attaching maps are locally Lipschitz or 
$C^k$ or smooth or analytic (etc.). A map between CW-complexes
is locally Lipschitz or $C^k$ or smooth or analytic (etc.) just when all of
its restrictions to each simplex are.
\end{definition}
For example, a manifold or real or complex analytic variety,
possibly with boundary and corners, is a locally Lipschitz
CW-complex. Clearly the proofs above only require that $X$ be
connected by Lipschitz curves, so hold equally well with $X$ any
connected locally Lipschitz CW-complex, rather than a manifold, the
developed map $\phi_0$ locally Lipschitz, and the isomorphism $\Phi$
locally essentially bounded. In particular, in
theorem~\vref{thm:Main}, the curves we roll on can have arbitrary
analytic singularities. One could instead try to uniformize the
curve first, but then it would not be so clear that the development
could be ``deuniformized''. Locally Lipschitz development is likely
to be useful in developing calibrated cycles in studying Cartan
geometries modelled on symmetric spaces.

\section{Applications}

\begin{theorem}
Suppose that $M_0$ and $M_1$ are manifolds and bear real/complex
analytic $G/H$-geometries. Suppose that the model $G/H$ rolls
freely on $M_1$. Every local development of a real/complex analytic map $\phi_0 : X \to
M_0$ from a simply connected analytic variety $X$ extends to a
real/complex analytic development $\phi_1 : X \to M_1$. Moreover
$\phi_0^* TM_0 = \phi_1^* TM_1$ are isomorphic vector bundles on
$X$.
\end{theorem}
\begin{proof}
Clearly it is enough to extend the development along all real curves in $X$. The
local obstructions vanish by analyticity. The vector bundle
isomorphism follows from lemma~\vref{lemma:TgtBundle}.
\end{proof}

\begin{corollary}
Let $M$ be a complex manifold bearing a complex analytic
Cartan geometry modelled on a homogeneous space $G/H$, and that curves in the model roll on $M$.
Each point of $M$ lies in the image of a holomorphic map $\C{} \to M$.
Moreover $M$ contains a rational curve through every point just if
the model $G/H$ contains a rational curve, and otherwise $M$
contains no rational curves.
\end{corollary}
For example, Kobayashi hyperbolic complex manifolds admit no
complete holomorphic Cartan geometries.
\begin{proof}
We need only prove that $G$ has a complex subgroup whose orbit in
$G/H$ is a complex curve. This is easy to see if $G$ contains a
semisimple group, since $G$ then contains $\SL{2,\C{}}$, and we can
just examine the homogeneous spaces of $\SL{2,\C{}}$ by examining
closed subgroups. If $G$ contains no semisimple group, then $G$ is
solvable and contains a complex abelian subgroup of positive
dimension, and the result is obvious.
\end{proof}

\section{Ideas for further research}

Lebrun and Mason \cite{LeBrunMason:2002} found a twistor approach to construct
Zoll surfaces as moduli spaces of pseudoholomorphic
disks. It might be possible to generalize this for parabolic geometries
associated to split form semisimple Lie groups.

 Suppose that morphism curvature is valued in a representation $V$,
and that $V \subset \Sym{2}{W}$ or $V \subset W^* \otimes W$ for some
representation $W$. Then we can consider the signs of eigenvalues
of morphism curvature.
This might lead to Toponogov triangle theorems for Cartan geometries.

One can define morphism graphs, a directed graph with homogeneous spaces at
vertices and morphisms at edges, with any two paths with the same
end points representing the same composition of morphisms.
This is a natural way to look at development, in
the context of generalized Cartan geometries.
A development is really a morphism between morphisms
of generalized Cartan geometries.

In the complex analytic setting, perhaps morphisms
from a rational homogeneous variety would have unobstructed deformation
theory, and perhaps all holomorphic map deformations are
morphism deformations.

It would be nice to know which 2-plane fields on 5-manifolds come
from rolling one surface on another. More generally,
one naturally wants to understand how much surface geometry is encoded in
the 2-plane field associated to a pair of surfaces.

It would help to be able to say which morphisms $\Phi$
of homogeneous spaces have the property that
every morphism modelled on $\Phi$ will have
vanishing morphism curvature, in particular, for
$\Phi$ a morphism of rational homogeneous varieties.

Suppose that $G$ is the split real form of a semisimple Lie group,
and $P \subset G$ is a parabolic subgroup.
Suppose that $G$ has Lie algebra $\mathfrak{g}$, and suppose that $\mathfrak{g}$
is filtered, with $\mathfrak{p}=\mathfrak{g}^0$ the $0$-degree part. 
(Every parabolic Lie subalgebra arises in this way.)
Take any $G/P$-geometry, say $P \to E \to M$. 
The filtration induces a $P$-module filtration of 
of $\mathfrak{g}/\mathfrak{p}$, 
\[
0 \subset V_1 \subset V_2 \subset \dots \subset V_N = \mathfrak{g}/\mathfrak{p}
\]
for which $V_1$ is a direct sum of irreducible representations;
see \cite{Baston/Eastwood:1989,Cap:2006}. This imposes
a vector bundle $E \times_P V_1 \subset TM$.
Let $P^c$ be the maximal compact subgroup of $P$. Then $E \to M$ has
a principal $P^c$-subbundle, say $E^c \to M$. (There are many
such subbundles, so pick one.) Then $E^c \to M$ 
imposes a canonical Riemannian metric on $M$ 
as in definition~\vref{definition:CanonicalMetricOnCartanGeometry}.
\begin{conjecture}
If this Riemannian metric on $M$ has positive Ricci curvature along $E \times_P V_1$, then
the $G/P$-geometry is complete.
\end{conjecture} 
If proven, this conjecture will apply to nondegenerate 2-plane fields.

There is nothing known about Cartan geometries
on Banach manifolds. Some of these theorems
might hold in that wider context.

% ----------------------------------------------------------------
\bibliographystyle{amsplain}
\bibliography{CompleteCartanConnections}
\end{document}